\documentclass[12pt,reqno]{amsart}
\usepackage{graphicx} 
\numberwithin{equation}{section}

\usepackage{esint}
\usepackage{amsmath,amssymb,mathrsfs}
\usepackage{makecell}

\usepackage[margin=3cm]{geometry}
\allowdisplaybreaks

\newtheorem{theorem}{Theorem}[section]
\newtheorem{lemma}[theorem]{Lemma}

\theoremstyle{remark}
\newtheorem{definition}[theorem]{Definition}
\newtheorem{corollary}[theorem]{Corollary}
\newtheorem{remark}[theorem]{Remark}

\usepackage{tikz}
\usetikzlibrary{arrows.meta, calc}
\usepackage{caption}

\usepackage[colorlinks=true, linkcolor=blue, citecolor=blue, urlcolor=blue]{hyperref}

\title[Maximal function on spaces of homogeneous type]
{The maximal function on spaces of homogeneous type, or adjacent dyadic cubes do good}

\author{Alina Shalukhina}  
\address{
Center for Mathematics and Applications and
Department of Mathematics,
NOVA School of Science and Technology,
NOVA University Lisbon,
Quinta da Torre,
2829-516 Caparica, Portugal}
\email{a.shalukhina@campus.fct.unl.pt}

\begin{document}

\begin{abstract}
We prove that the Hardy--Littlewood maximal operator $M$ is bounded on the variable 
Lebesgue space $L^{p(\cdot)}(X,d,\mu)$, with $1<p_-\le p_+<\infty$, over an unbounded space
of homogeneous type $(X,d,\mu)$ with a Borel-semiregular measure $\mu$, if and only if the
averaging operators $T_\mathcal{Q}$ are bounded on $L^{p(\cdot)}(X,d,\mu)$ uniformly over 
all families $\mathcal{Q}$ of pairwise disjoint ``cubes'' from a Hyt\"onen--Kairema
dyadic system on $X$. This extends Diening's well-known characterization of the boundedness 
of $M$ on $L^{p(\cdot)}(\mathbb{R}^n)$ to the setting of
spaces of homogeneous type, while also providing a slight refinement of the original result.



\end{abstract}

\maketitle

\section{Introduction}
Spaces of homogeneous type are, in brief, quasi-metric spaces with doubling 
measure. Compact Riemannian manifolds with the Riemannian metric and volume measure,
quasi-metric spaces with Ahlfors regular measure (e.g., Lipschitz surfaces and fractal
sets), the Heisenberg group with the Carnot--Carath\'eodory metric and Haar 
measure---to give just a few examples---can all 
be represented as a triple $(X,d,\mu)$, where $d$ is a quasi-metric and there is 
a constant $A\ge1$ such that
\[
\mu(B(x,r))\le A\,\mu(B(x,r/2))
\]
for all balls $B(x,r):=\{y\in X:d(x,y)<r\}$ with $x\in X$ and $r>0$. 
The classical Hardy--Littlewood maximal function naturally appears on a space of homogeneous
type $(X,d,\mu)$. For a measurable function $f$ on $X$, it is defined by 
\[
Mf(x):=\sup_{B\ni x}\frac{1}{\mu(B)}\int_B|f(y)|\,d\mu(y),
\quad x\in X,
\]
where the supremum is taken over all balls $B\subseteq X$ containing the point $x$.
A number of recent works have been devoted to the study of the maximal operator 
$M:f\mapsto Mf$ on Banach function spaces over spaces of homogeneous type, with one
particularly interesting setting being the variable Lebesgue spaces 
$L^{p(\cdot)}(X,d,\mu)$. Among others, Cruz-Uribe and 
Cummings~\cite{CUC22} obtained weighted norm inequalities for
the maximal operator on $L^{p(\cdot)}(X,d,\mu)$; Karlovich~\cite{AK19} investigated the 
boundedness of $M$ on the associate space of an arbitrary Banach function space over 
$(X,d,\mu)$; Kokilashvili, Meskhi, Rafeiro and Samko~\cite{KMRS16} introduced grand
variable Lebesgue spaces over spaces of homogeneous type and established boundedness 
theorems for the maximal operator, singular integrals, and potentials in these new spaces.

At the same time, a well-developed theory exists for the maximal operator $M$ on variable
Lebesgue spaces $L^{p(\cdot)}(\mathbb{R}^n)$ over the Euclidean space 
$\mathbb{R}^n$---itself a prototypical example of a space of homogeneous type.
The books by Diening et al. \cite{DHHR11} and Cruz-Uribe and Fiorenza \cite{CUF13} provide a
comprehensive account of the early results in this area, while more recent research 
includes Lerner~\cite{AL25}, who established a boundedness criterion for $M$ on the space
$L^{p(\cdot)}(\mathbb{R}^n)$ with a bounded exponent $p(\cdot)$, and Adamadze, 
Diening and Kopaliani~\cite{ADK25}, who obtained a sufficient condition for the 
boundedness of $M$ on the space $L^{p(\cdot)}(\mathbb{R}^n)$ with an unbounded exponent
$p(\cdot)$. The central problem has always been, and remains, the
characterization of the spaces on which $M$ is bounded.

In 2005, Diening established a curious characterization of the boundedness of $M$
on variable Lebesgue spaces $L^{p(\cdot)}(\mathbb{R}^n)$,
with exponents $p(\cdot)$ bounded away from one and infinity, 
in terms of the uniform boundedness of a family of averaging operators on 
those $L^{p(\cdot)}(\mathbb{R}^n)$. Namely, he proved 
in~\cite[Theorem~8.1, (i)$\Leftrightarrow$(ii)]{LD05} the following result:

\begin{theorem}[Diening]\label{th:Diening}
The maximal operator $M$ is bounded on a variable Lebesgue space 
$L^{p(\cdot)}(\mathbb{R}^n)$ with an exponent $p(\cdot)$ satisfying 
\[
1<\operatornamewithlimits{ess\,inf}_{x\in\mathbb{R}^n}p(x)
\le\operatornamewithlimits{ess\,sup}_{x\in\mathbb{R}^n}p(x)<\infty
\]
if and only if there exists a constant $C>0$ such that
\begin{equation}\label{intro:averages}
\bigg\|\sum_{Q\in\mathcal{Q}}\bigg(\frac{1}{|Q|}\int_Q|f(x)|\,dx\bigg)
\,\chi_Q\bigg\|_{p(\cdot)}\le C\,\|f\|_{p(\cdot)}
\end{equation}
for all families $\mathcal{Q}$ of pairwise disjoint open cubes in $\mathbb{R}^n$ 
and all $f\in L^{p(\cdot)}(\mathbb{R}^n)$.
\end{theorem}

The condition of uniform boundedness of the averaging operators 
is examined in detail in the
recent paper by Nieraeth~\cite{ZN26}, where it is called the strong Muckenhoupt
condition for $L^{p(\cdot)}(\mathbb{R}^n)$. Vinogradov~\cite{OV25} studies the 
related problem of uniform boundedness of Steklov averaging 
operators on weighted variable Lebesgue spaces $L^{p(\cdot)}_w(\mathbb{T}^n)$ of periodic
functions.

Though Diening's criterion is not constructive---and hence, is difficult to apply in 
practice---it had important theoretical consequences. Among other things, it 
guaranteed~\cite{CUFMP06} that the theory of Rubio de Francia extrapolation could be 
extended to the variable Lebesgue
spaces. As noted in~\cite[p.~8]{CUF13}, this ``allowed the theory of weighted norm 
inequalities to be used to prove that a multitude of operators (such as singular integrals)
are bounded on $L^{p(\cdot)}(\mathbb{R}^n)$ whenever the maximal operator is.'' Another 
application of Diening's characterization recently appeared in Lerner’s 
criterion~\cite{AL25} for the boundedness
of $M$ on $L^{p(\cdot)}(\mathbb{R}^n)$, formulated in terms of the exponent $p(\cdot)$.

\emph{We generalize Diening's characterization from the Euclidean space $\mathbb{R}^n$ to 
spaces of homogeneous type, and slightly refine his result in $\mathbb{R}^n$.} 
But what could play the role of cubes in 
a space of homogeneous type $(X,d,\mu)$? In fact, there are certain ``cubes'' in this
abstract setting; they are Borel sets on $X$ organized into a system $\mathscr{D}$ of 
adjacent dyadic grids, constructed by Hyt\"onen and Kairema in~\cite{HK12}. An example of
such a system in $\mathbb{R}^n$ is the classical collection of $3^n$ adjacent dyadic
grids
\[
\mathscr{D}^{\mathbf{t}}:=\left\{
2^{-k} \big([0,1)^n+\mathbf{m}+(-1)^k\mathbf{t}\big)
\ :\ 
k\in\mathbb{Z},\,\mathbf{m}\in\mathbb{Z}^n\right\},
\quad\mathbf{t}\in\left\{0,\frac13,\frac23\right\}^n,
\] 
cleverly shifted by one-thirds. It is these dyadic ``cubes''---taken from a
Hyt\"onen--Kairema system defined in Theorem~\ref{th:HK}---that will do us good, as the 
title suggests. Building on their key properties, we prove our main result 
(restated Theorem~\ref{th:main}):

\begin{theorem}[Main result]\label{th:intromain}
Let $(X,d,\mu)$ be an unbounded space of homogeneous type with a Borel-semiregular measure
$\mu$, and let $\mathscr{D}$ be a Hyt\"onen--Kairema dyadic system on $X$. 
Then the maximal operator $M$ is bounded on a variable Lebesgue 
space $L^{p(\cdot)}(X,d,\mu)$ with an exponent $p(\cdot)$ satisfying
\[
1<\operatornamewithlimits{ess\,inf}_{x\in X}p(x)
\le\operatornamewithlimits{ess\,sup}_{x\in X}p(x)<\infty
\]
if and only if there is a constant $C>0$ such that
\[
\bigg\|\sum_{Q\in\mathcal{Q}}\bigg(\frac{1}{\mu(Q)}\int_Q|f(x)|\,d\mu(x)\bigg)
\,\chi_Q\bigg\|_{p(\cdot)}\le C\,\|f\|_{p(\cdot)}
\]
for all families $\mathcal{Q}$ of pairwise disjoint ``cubes'' $Q\in\mathscr{D}$ and
all $f\in L^{p(\cdot)}(X,d,\mu)$.
\end{theorem}

This general result, when applied to $\mathbb{R}^n$, yields a slight refinement of 
Theorem~\ref{th:Diening}: it suffices to verify
condition~\eqref{intro:averages} for all families $\mathcal{Q}$ of pairwise disjoint 
\emph{dyadic} cubes in $\mathbb{R}^n$---selected from the $3^n$ adjacent dyadic grids
$\mathscr{D}^\mathbf{t}$, $\mathbf{t}\in\{0,\frac13,\frac23\}^n$---rather than \emph{all} 
cubes in $\mathbb{R}^n$, to ensure that the maximal operator $M$ is bounded on 
$L^{p(\cdot)}(\mathbb{R}^n)$. This noteworthy consequence of Theorem~\ref{th:intromain}
is stated as Corollary~\ref{cor:main}.

Let the reader not be discouraged by the length of the proof of Theorem~\ref{th:intromain}.
Our proof largely follows the scheme developed by Diening for his original 
characterization; however, beyond merely proving our main result,
we set ourselves the additional aim of presenting the proof as clear as 
possible---filling in gaps from the original arguments in~\cite[Chapter~5]{DHHR11}, 
correcting missing or superfluous conditions, and providing additional intuition. We hope 
the text now reads easily, like a book.

We faced two main challenges on the way to generalizing Diening's result, each arising
from a sticking point where the original argument would not go through analogously. The 
first of these was to obtain Lemma~\ref{le:perturbation} about perturbation of a 
cube-constant function on dyadic ``cubes,'' and this became possible by using a version of 
the Calder\'on--Zygmund decomposition (Lemma~\ref{le:CZ}) for spaces of 
homogeneous type. The second one was to prove Lemma~\ref{le:aux-to-suff} whose prototype
in $\mathbb{R}^n$ relied on the Besicovitch covering theorem---a tool specific to Euclidean
spaces and no longer available when we pass from $\mathbb{R}^n$ to spaces of homogeneous
type. We designed an original proof for Lemma~\ref{le:aux-to-suff}, combining the use
of Hyt\"onen and Kairema's system of adjacent dyadic grids with the application of the 
Calder\'on--Zygmund-type Lemma~\ref{le:CZ} on each grid.

Following Diening's approach, we begin in the general functional setting of 
Musielak--Orlicz spaces and later restrict to variable Lebesgue spaces, which form
a special case of those more general spaces. Thus, Section~\ref{se:MO} is devoted to
Musielak--Orlicz spaces; Section~\ref{sec:SHT} presents the dyadic system
$\mathscr{D}$ on spaces of homogeneous type; and Section~\ref{sec:PHI} discusses 
generalized $\Phi$-functions on ``cubes,'' which lead to Musielak--Orlicz 
sequence spaces. These three sections provide the necessary preliminaries for the proof 
of the main result. Then, Section~\ref{sec:A} introduces the self-improving class 
$\mathscr{A^D}$ of generalized $\Phi$-functions; Section~\ref{sec:sufficient} restricts it
to the class $\mathscr{A}^{\mathscr{D}}_\textnormal{strong}$ and gives a sufficient
condition for the boundedness of $M$ on Musielak--Orlicz spaces in terms of this new,
strong class. Finally, it is proved in Section~\ref{sec:V-L} that the classes 
$\mathscr{A^D}$ and $\mathscr{A}^{\mathscr{D}}_\textnormal{strong}$ coincide for those
generalized $\Phi$-functions which generate the variable Lebesgue spaces with the
exponent function bounded away from one and infinity, and we prove our main result in
Section~\ref{sec:main}.

\section{Musielak--Orlicz Spaces}\label{se:MO}
We start from giving a brief overview of Musielak--Orlicz spaces and their basic
properties; our summary mainly follows~\cite[Sections~2.3, 2.6 and 2.7]{DHHR11}. 
Then we define variable Lebesgue spaces as Musielak--Orlicz spaces, in the spirit 
of~\cite[Chapter~3]{DHHR11}.

\subsection{Semimodular spaces induced by generalized $\Phi$-functions}

Musielak--Orlicz spaces are defined via the so-called generalized $\Phi$-functions on a
measure space, which extend the classical concept of a $\Phi$-function.
\begin{definition}
A convex, left-continuous function $\varphi:[0,\infty)\to[0,\infty]$ with $\varphi(0)=0$,
$\lim_{t\to0^+}\varphi(t)=0$ and $\lim_{t\to\infty}\varphi(t)=\infty$ is called a 
{\it $\Phi$-function}. It is said to be {\it positive} if $\varphi(t)>0$ when $t>0$.
\end{definition}
Note that the convexity of $\varphi$ and condition $\varphi(0)=0$ together give two 
useful rules for extracting a scaling factor from the function's argument: for all $t\ge0$, 
one has
\begin{align*}
&\varphi(\lambda t)\le\lambda\,\varphi(t)
\quad\text{if }\lambda\in[0,1],
\\
&\varphi(\lambda t)\ge\lambda\,\varphi(t)
\quad\text{if }\lambda\ge1;
\end{align*}
see~\cite[p.~6]{HH19}. Moreover, since $\varphi$ takes only 
nonnegative values, it follows that $\varphi$ is \emph{non-decreasing} on $[0,\infty)$: if 
$0\le t_1\le t_2$, then $t_1/t_2\in[0,1]$ and thus
\[
\varphi(t_1)=\varphi\left(\frac{t_1}{t_2}\,t_2\right)\le
\frac{t_1}{t_2}\,\varphi(t_2)\le\varphi(t_2).
\]

Let $(X,\mu)$ be an arbitrary $\sigma$-finite, complete measure space. By $L^0(X,\mu)$ we 
denote the space of all complex-valued, $\mu$-measurable functions on $X$.
\begin{definition}\label{def:gen-Phi}
We say that $\varphi:X\times[0,\infty)\to[0,\infty]$ is a 
{\it generalized $\Phi$-function} on $(X,\mu)$, and write $\varphi\in\Phi(X,\mu)$, if:
\begin{enumerate}
    \item[(i)] $t\mapsto\varphi(x,t)$ is a $\Phi$-function for every $x\in X$;
    \item[(ii)] $x\mapsto\varphi(x,t)$ is $\mu$-measurable for every $t\ge0$.
\end{enumerate}
\end{definition}
A function $\varphi\in\Phi(X,\mu)$ induces the Musielak--Orlicz space $L^\varphi(X,\mu)$
as a semimodular space. That is, first $\varphi$ generates the functional $\rho_\varphi$
on $L^0(X,\mu)$ defined by 
\[
\rho_\varphi(f):=\int_X\varphi(x,|f(x)|)\,d\mu(x).
\]
This functional is a {\it semimodular}, and in particular it is convex and left-continuous 
(in the sense that the mapping $\lambda\mapsto\rho_\varphi(\lambda f)$ is left-continuous
on $(0,\infty)$ for every measurable $f$), 
see~\cite[Lemma~2.3.10 and Definition~2.1.1]{DHHR11}. Then the semimodular $\rho_\varphi$
induces the space
\[
L^\varphi(X,\mu):=\left\{f\in L^0(X,\mu)\ :\ \rho_\varphi(\lambda f)<\infty
\text{ for some }\lambda>0 \right\}
\]
called the {\it Musielak--Orlicz space}. This is a Banach space with respect to the 
Luxemburg norm
\[
\|f\|_\varphi:=\inf\{\lambda>0\ :\ \rho_\varphi(f/\lambda)\le1\},
\]
as shown in~\cite[Theorem~2.3.13]{DHHR11}. Note that $f\in L^\varphi(X,\mu)$ if and 
only if $\|f\|_\varphi<\infty$. Here are some fundamental properties of Musielak--Orlicz
spaces.
\begin{lemma}\label{le:spaces}
Let $\varphi\in\Phi(X,\mu)$. The space $L^\varphi(X,\mu)$ has the following properties:
\begin{enumerate}
    \item[(a)] \emph{(circularity)} $\|f\|_\varphi=\|\,|f|\,\|_\varphi$
    for all $f\in L^\varphi(X,\mu)$.
    \item[(b)] \emph{(lattice property)} If $f\in L^\varphi(X,\mu)$, $g\in L^0(X,\mu)$,
    and $0\le|g|\le|f|$ holds $\mu$-almost everywhere,
    then $g\in L^\varphi(X,\mu)$ and $\|g\|_\varphi\le\|f\|_\varphi$.
    \item[(c)] \emph{(the Fatou property)} If $0\le f_k\nearrow f$ holds $\mu$-almost 
    everywhere for a sequence $\{f_k\}_{k=0}^\infty\subseteq L^\varphi(X,\mu)$ satisfying
    $\sup_{k\ge0}\|f_k\|_\varphi<\infty$, then $f\in L^\varphi(X,\mu)$
    and $\|f_k\|_\varphi\nearrow\|f\|_\varphi$.
    \item[(d)] \emph{(the Riesz--Fischer property)} For any sequence 
    $\{f_k\}_{k=0}^\infty\subseteq L^0(X,\mu)$, one has 
    \[
    \bigg\|\sum_{k=0}^\infty f_k\bigg\|_\varphi\le\sum_{k=0}^\infty\|f_k\|_\varphi.
    \]
    \item[(e)] \emph{(unit ball property)} $\|f\|_\varphi\le1$ if and only if 
    $\rho_\varphi(f)\le1$.
\end{enumerate}
\end{lemma}
\begin{proof}
Property~(a) immediately follows from $\rho_\varphi(f)=\rho_\varphi(|f|)$. The lattice 
property~(b) results from the fact that $t\mapsto\varphi(x,t)$ is non-decreasing on 
$[0,\infty)$ for all $x\in X$, and thus $\rho_\varphi(g/\lambda)\le\rho_\varphi(f/\lambda)$
for all $\lambda>0$. For the Fatou property~(c), see~\cite[Theorem~2.3.17(d)]{DHHR11}.
The unit ball property~(e) is proved in~\cite[Lemma~2.1.14]{DHHR11}.

Property~(d) obviously holds if 
$\sum_{k=0}^\infty\|f_k\|_\varphi=\infty$. If, otherwise, 
$\{f_k\}_{k=0}^\infty\subseteq L^\varphi(X,\mu)$ and 
$\sum_{k=0}^\infty\|f_k\|_\varphi<\infty$, the proof of the Riesz--Fischer property 
repeats verbatim the argument in~\cite[Chapter~1, Theorem~1.6]{BS88}, since the space
$L^\varphi(X,\mu)$ satisfies the Fatou property~(c).
\end{proof}

The next lemma shows that pointwise relations between generalized $\Phi$-functions
allow comparison between the norms of the corresponding Musielak--Orlicz spaces. To state 
it, let us establish the notation used throughout. Given arbitrary sets $X$ and
$Y$ and mappings $\xi,\eta:X\times Y\to[0,\infty]$, we will write
``$\xi(x,y)\approx\eta(x,y)$ uniformly in $x\in X$ and $y\in Y$'' if there exists a constant
$c>1$ such that 
\[
\tfrac1c\,\xi(x,y)\le\eta(x,y)\le c\,\xi(x,y)
\]
for all $x\in X$ and $y\in Y$. This notation extends
analogously to any number of variables. 
\begin{lemma}\label{le:comparison}
Let $\varphi,\psi\in\Phi(X,\mu)$. The following hold:
\begin{enumerate}
    \item[(a)] \emph{(monotonicity)} If $\psi(x,t)\le\varphi(x,t)$ for all $x\in X$ and 
    $t\ge0$, then $\|f\|_\psi\le\|f\|_\varphi$ for every $f\in L^\varphi(X,\mu)$.
    \item[(b)] \emph{(power-type rescaling)} Suppose that $\psi(x,t)\approx\varphi(x,t^s)$, 
uniformly in $x\in X$ and $t\ge0$, for some $s>0$.
Then $\|f\|_\psi\approx\|\,|f|^s\|^{1/s}_\varphi$ uniformly in 
$f\in L^\psi(X,\mu)$.
\end{enumerate}
\end{lemma}
\begin{proof}
For~(a), take an $f\in L^\varphi(X,\mu)$. For all $x\in X$ and $\lambda>0$, we have by 
assumption that $\psi(x,|f(x)/\lambda|)\le\varphi(x,|f(x)/\lambda|)$. This implies the 
set inclusion 
\[
\{\lambda>0:\rho_\varphi(f/\lambda)\le1\}\subseteq
\{\lambda>0:\rho_\psi(f/\lambda)\le1\},
\]
which leads to $\|f\|_\varphi\ge\|f\|_\psi$.

Let us prove~(b). By assumption, there exists $c>1$ such that 
\begin{equation}\label{ineq:rescaling}
\tfrac1c\,\psi(x,t)\le\varphi(x,t^s)\le c\,\psi(x,t)
\end{equation}
for all $x\in X$ and $t\ge0$. On the one hand, by the convexity of $\varphi$ in the second 
variable, this implies $\psi(x,t)\le c\,\varphi(x,t^s)\le\varphi(x,ct^s)$, which leads to
\begin{align*}
\|f\|_\psi^s&=\inf\left\{\lambda^s>0\ :\ \int_X
\psi\left(x,\frac{|f(x)|}{\lambda}\right)\,d\mu(x)\le1\right\}
\\
&\le\inf\left\{\lambda^s>0\ :\ \int_X
\varphi\left(x,\frac{c\,|f(x)|^s}{\lambda^s}\right)d\mu(x)\le1\right\}
=c\,\|\,|f|^s\|_\varphi
\end{align*}
for any $f\in L^\psi(X,\mu)$. Another consequence of~\eqref{ineq:rescaling} due to the 
convexity of $\varphi$ is the inequality 
$\psi(x,t)\ge\frac1c\,\varphi(x,t^s)\ge\varphi(x,\frac{t^s}c)$ for all $x\in X$ and $t\ge0$,
which implies for an $f\in L^\psi(X,\mu)$ that 
\begin{align*}
\|f\|_{\psi}^s&=\inf\left\{\lambda^s>0\ :\ \int_X
\psi\left(x,\frac{|f(x)|}{\lambda}\right)d\mu(x)\le1\right\}
\\
&\ge\inf\left\{\lambda^s>0\ :\ \int_X
\varphi\left(x,\frac{|f(x)|^s}{c\,\lambda^s}\right)d\mu(x)\le1\right\}
=\frac1c\,\|\,|f|^s\|_\varphi.
\end{align*}
The obtained norm estimates give $\|f\|_\psi\approx\|\,|f|^s\|_\varphi^{1/s}$
uniformly in $f\in L^\psi(X,\mu)$.
\end{proof}

Many desirable properties of the spaces $L^\varphi(X,\mu)$ follow when the growth of 
$\varphi$ in the second variable is suitably restricted. The most common way
to impose such control is via the $\Delta_2$-condition, which provides a uniform bound
on the growth of $\varphi$.
\begin{definition}
A function $\varphi\in\Phi(X,\mu)$ is said to satisfy the \emph{$\Delta_2$-condition} if 
there exists a constant $D\ge2$ such that 
\begin{equation}\label{ineq:delta2}
\varphi(x,2t)\le D\,\varphi(x,t)
\end{equation}
for all $x\in X$ and $t\ge0$. The smallest such $D$ is referred to as the 
\emph{$\Delta_2$-constant} of $\varphi$.
\end{definition}
Note that the lower bound $D\ge2$ is due to the convexity of $\varphi$ in the 
second variable, since $2\,\varphi(x,t)\le\varphi(x,2t)\le D\,\varphi(x,t)$.
As a consequence of~\eqref{ineq:delta2}, for all $f\in L^0(X,\mu)$ one has 
$\rho_\varphi(2f)\le D\,\rho_\varphi(f)$,
so the semimodular $\rho_\varphi$ is said to satisfy the $\Delta_2$-condition with 
the same constant as $\varphi$. Here are several properties of the spaces 
$L^\varphi(X,\mu)$ under the $\Delta_2$-condition on $\varphi$. 
\begin{lemma}\label{le:delta2}
If $\varphi\in\Phi(X,\mu)$ satisfies the $\Delta_2$-condition with constant $D$, the
following properties hold:
\begin{enumerate}
    \item[(a)] \emph{(finite semimodular)} $\rho_\varphi(f)<\infty$ for all 
    $f\in L^\varphi(X,\mu)$;
    \item[(b)] \emph{(positive definiteness)} $\rho_\varphi(f)=0$ if and only if 
    $f=0$ $\mu$-almost everywhere;
    \item[(c)] \emph{(unit sphere property)} $\|f\|_\varphi=1$ if and only if
    $\rho_\varphi(f)=1$;
    \item[(d)] \emph{(norm-to-semimodular lower bound)} for every $\varepsilon>0$ there 
    exists a small number $\delta=\delta(\varepsilon,D)\in(0,1)$ such that, for all 
    $f\in L^0(X,\mu)$,
    \[
    \|f\|_\varphi\ge\varepsilon\implies\rho_\varphi(f)\ge\delta;
    \]
    \item[(e)] \emph{(semimodular-to-norm upper bound)} for every small number 
    $\varepsilon\in(0,1)$ there exists $\delta=\delta(\varepsilon,D)\in(0,1)$ such that, for 
    all $f\in L^\varphi(X,\mu)$, 
    \[
    \rho_\varphi(f)\le\varepsilon\implies\|f\|_\varphi\le\delta.
    \]
\end{enumerate}
\end{lemma}
\begin{proof}
For claim~(a), see the discussion in~\cite[above Remark~2.5.4]{DHHR11}. 
Claim~(b) is proved in~\cite[Lemma~2.4.3]{DHHR11} in view 
of~\cite[Definition~2.1.1(f)]{DHHR11}. Property~(c) follows
from~\cite[Lemma~2.1.14 and Lemma~2.4.3]{DHHR11}. Claim~(d) is a restatement
of~\cite[Lemma~2.4.2]{DHHR11}, and claim~(e) is proved in~\cite[Lemma~2.4.3]{DHHR11}.
\end{proof}

\subsection{Conjugates of generalized $\Phi$-functions}

Within the theory of Musielak--Orlicz spaces, each generalized $\Phi$-function
is paired with its conjugate, defined in the spirit of the Legendre transform.
\begin{definition}\label{def:conjug}
Let $\varphi\in\Phi(X,\mu)$. The {\it conjugate function} of $\varphi$, denoted by 
$\varphi^*$, is defined for all $x\in X$ and $u\ge0$ by 
\[
\varphi^*(x,u):=\sup_{t\ge0}(tu-\varphi(x,t)).
\]
\end{definition}
The conjugate function of $\varphi\in\Phi(X,\mu)$ is again a generalized $\Phi$-function
on $(X,\mu)$: this is implicitly stated in~\cite[see discussion after formula~(2.6.2), 
Theorem~2.2.4, and Lemma~2.3.2]{DHHR11}. Moreover, $\varphi^*\in\Phi(X,\mu)$ if and only if
$\varphi\in\Phi(X,\mu)$, and there holds the involutive property
\[
(\varphi^*)^*=\varphi,
\]
see \cite[Corollary~2.6.3]{DHHR11}, which we will frequently use. 
An immediate consequence of Definition~\ref{def:conjug} is {\it Young's inequality}
\begin{equation}\label{eq:Young}
tu\le\varphi(x,t)+\varphi^*(x,u),
\end{equation}
which holds for all $x\in X$ and $t,u\ge0$.
The next lemma follows from~\cite[Lemma~2.6.4]{DHHR11}.
\begin{lemma}\label{le:conjugation}
Let $\varphi,\psi\in\Phi(X,\mu)$. The following hold: 
\begin{enumerate}
    \item[(a)] If $\varphi(x,t)\le\psi(x,t)$ for all $x\in X$ and $t\ge0$, then 
    $\psi^*(x,u)\le\varphi^*(x,u)$ for all $x\in X$ and $u\ge0$.
    \item[(b)] Suppose that $a,b>0$ and $\psi(x,t)=a\varphi(x,bt)$ for all $x\in X$ and
    $t\ge0$. Then $\psi^*(x,u)=a\varphi^*(x,\frac{u}{ab})$ for all $x\in X$ and $u\ge0$.
\end{enumerate}
\end{lemma}
Finally, there is a version of H\"older's inequality for the conjugates $\varphi$ and
$\varphi^*$ \cite[Lemma~2.6.5]{DHHR11}. In what follows, let the angle brackets 
$\langle f,g\rangle$ denote the integral of the product of two nonnegative functions 
$f,g\in L^0(X,\mu)$ over the entire space $X$, i.e.,
\[
\langle f,g\rangle:=\int_Xf(x)g(x)\,d\mu(x).
\]
\begin{lemma}[H\"older's inequality]\label{le:Holder}
Let $\varphi\in\Phi(X,\mu)$. Then 
$\langle |f|,|g|\rangle\le2\,\|f\|_\varphi\|g\|_{\varphi^*}$
for all $f\in L^\varphi(X,\mu)$ and $g\in L^{\varphi^*}(X,\mu)$.
\end{lemma}

\subsection{Proper generalized $\Phi$-functions}

The notion of a proper function $\varphi\in\Phi(X,\mu)$, introduced
in~\cite{DHHR11}, was originally designed to characterize the necessary and
sufficient condition for $L^\varphi(X,\mu)$ to be a Banach function space
in the sense of Bennett and Sharpley~\cite[Chapter~1]{BS88}. 
Our definition of a proper function is equivalent to that 
in~\cite[Definition~2.7.8]{DHHR11} due to~\cite[Corollary~2.7.9]{DHHR11}.
Hereafter, the symbol $\chi_E$ will denote the characteristic function of a set $E$.
\begin{definition}
We say that $\varphi\in\Phi(X,\mu)$ is {\it proper} if 
$\chi_E\in L^\varphi(X,\mu)\cap L^{\varphi^*}(X,\mu)$ 
for all $E\subseteq X$ with 
$\mu(E)<\infty$.
\end{definition}
Since $\varphi=(\varphi^*)^*$, it is clear from the definition that $\varphi$ is proper if 
and only if $\varphi^*$ is proper.
The following lemma gives part of the norm conjugate formula for the spaces 
$L^\varphi(X,\mu)$ with proper $\varphi$; it is a special case 
of~\cite[Corollary~2.7.5]{DHHR11}.
\begin{lemma}\label{le:norm-conj}
Let $\varphi\in\Phi(X,\mu)$ be proper. Then for every $f\in L^0(X,\mu)$, one has
\[
\|f\|_\varphi\le\sup_{\substack{g\in L^{\varphi^*}(X,\mu):\\
\|g\|_{\varphi^*}\le1}}
\langle |f|,|g|\rangle.
\]   
\end{lemma}
The property of being proper conveniently extends to all generalized $\Phi$-function that
are pointwise equivalent to a proper function.
\begin{lemma}\label{le:approper}
Let $\varphi,\psi\in\Phi(X,\mu)$ be such that $\varphi(x,t)\approx\psi(x,t)$
uniformly in $x\in X$ and $t\ge0$. If $\varphi$ is proper, then $\psi$ is also proper.
\end{lemma}
\begin{proof}
By Lemma~\ref{le:comparison}(b) with $s=1$, we deduce from the pointwise equivalence
of $\varphi$ and $\psi$ that  
$\|f\|_\varphi\approx\|f\|_\psi$ uniformly in $f\in L^\varphi(X,\mu)=L^\psi(X,\mu)$.
Fix a set $E\subseteq X$ with $\mu(E)<\infty$. Since $\varphi$ is proper, it follows that 
$\chi_E\in L^\varphi(X,\mu)$ and thus $\chi_E\in L^\psi(X,\mu)$.

To prove that $\chi_E\in L^{\psi^*}(X,\mu)$, note that by assumption there exists $c>1$ 
such that $\frac1c\,\varphi(x,t)\le\psi(x,t)\le c\,\varphi(x,t)$ for all $x\in X$ and 
$t\ge0$. Due to the convexity of $\varphi$ in the second variable, this implies 
$\varphi(x,\frac{t}{c})\le\psi(x,t)\le\varphi(x,ct)$ for all $x\in X$ and $t\ge0$, 
which is equivalent to 
\[
\varphi^*\left(x,\frac{t}{c}\right)\le\psi^*(x,t)\le\varphi^*(x,ct),
\]
for all $x\in X$ and $t\ge0$, in view of Lemma~\ref{le:conjugation}. By a similar argument
as in the proof of Lemma~\ref{le:comparison}(b), we deduce from the last inequality that 
$\frac1c\,\|f\|_{\varphi^*}\le\|f\|_{\psi^*}\le c\,\|f\|_{\varphi^*}$ for any 
measurable function $f$, and hence $L^{\varphi^*}(X,\mu)=L^{\psi^*}(X,\mu)$. Since $\varphi$
is proper, we have $\chi_E\in L^{\varphi^*}(X,\mu)$ and therefore  
$\chi_E\in L^{\psi^*}(X,\mu)$. Thus $\psi$ is proper.
\end{proof}

The next, purely technical, lemma is best understood in light of  
Definition~\ref{def:conjug}. It offers a way to represent an integral of 
the conjugate function by a similar supremum construction as in the definition of the
conjugate function itself.
\begin{lemma}\label{le:int-phi*}
Let $\varphi\in\Phi(X,\mu)$ be proper. Then for all $u\ge0$ and $E\subseteq X$ with 
$\mu(E)<\infty$, one has
\begin{equation}\label{eq:int-phi*}
\int_E \varphi^*(x,u)\,d\mu(x)=\sup_{f\in L^1(E,\mu)}
\int_E\left(u\,|f(x)|-\varphi(x,|f(x)|)\right)\,d\mu(x).
\end{equation}
\end{lemma}
\begin{proof}
Since $\varphi$ is proper, the conclusion of~\cite[Theorem~2.7.4]{DHHR11} holds
for the function $g=u\chi_E\in L^{\varphi^*}(X,\mu)$. Written out using the definitions from
\cite[Definition~2.2.2 and formula~(2.7.3)]{DHHR11}, this theorem states that 
\begin{equation}\label{eq:int-phi*-0}
\int_E \varphi^*(x,u)\,d\mu(x)=\sup_{f\in L^\varphi(X,\mu)}
\left(u\left|\int_E f(x)\,d\mu(x)\right|-\rho_\varphi(f)\right).
\end{equation}
We need to check that the right-hand side of~\eqref{eq:int-phi*-0} coincides with the 
right-hand side of~\eqref{eq:int-phi*}. 
Take a function $f\in L^\varphi(X,\mu)$. Since $\varphi$ is proper, we have that 
$\chi_E\in L^{\varphi^*}(X,\mu)$. Then it follows by
H\"older's inequality in Lemma~\ref{le:Holder} that
\[
\int_E|f(x)|\,d\mu(x)=\langle|f|,\chi_E\rangle
\le2\,\|f\|_\varphi\|\chi_E\|_{\varphi^*}<\infty,
\]
so $f\in L^1(E,\mu)$. Then for any $f\in L^\varphi(X,\mu)$, there obviously holds
\begin{align*}
u\left|\int_E f(x)\,d\mu(x)\right|-\rho_\varphi(f)
&\le u\int_E |f(x)|\,d\mu(x)-\rho_\varphi(f\chi_E)
\\
&\le\sup_{f\in L^1(E,\mu)}
\int_E\left(u\,|f(x)|-\varphi(x,|f(x)|)\right)\,d\mu(x).
\end{align*}
Passing to the supremum over all $f\in L^\varphi(X,\mu)$ in the left-hand side of the
above inequality, we prove that ``\eqref{eq:int-phi*-0} $\le$ \eqref{eq:int-phi*}.'' 

To prove the converse, let us take a function $f\in L^1(E,\mu)$ and extend it from $E$ 
to $X$ by setting $f(x)=0$ for all $x\in X\setminus E$. For this extended $f\in L^0(X,\mu)$,
find a sequence of simple functions $\{g_k\}_{k=1}^\infty$, which all belong 
to $L^\varphi(X,\mu)$ since $\varphi$ is proper, such that
\[
0\le g_1(x)\le g_2(x)\le\ldots\le |f(x)|
\quad\text{and}\quad
g_k(x)\to|f(x)|=|f(x)|\,\chi_E(x)
\]
for $\mu$-almost every $x\in X$; see, e.g., \cite[Theorem~2.10(a)]{GF99}. 
By the monotonicity 
and left-continuity of $t\mapsto\varphi(x,t)$, we consequently have that
$0\le\varphi(\cdot,g_k(\cdot))\nearrow\varphi(\cdot,|f(\cdot)|\,\chi_E(\cdot))$. Then it
follows from the monotone convergence theorem~\cite[Theorem~2.14]{GF99}, applied to the 
non-decreasing sequences 
$\{g_k\}_{k=1}^\infty$ and $\{\varphi(\cdot,g_k(\cdot))\}_{k=1}^\infty$, that
\begin{align*}
\int_E\left(u\,|f(x)|-\varphi(x,|f(x)|)\right)\,d\mu(x)
&=\int_E u\,\lim_{k\to\infty}g_k(x)\,d\mu(x)
-\int_X\lim_{k\to\infty}\varphi(x,g_k(x))\,d\mu(x)
\\
&=\lim_{k\to\infty}\left(\int_E u\,g_k(x)\,d\mu(x)-\rho_\varphi(g_k)\right)
\\
&\le\sup_{k\in\mathbb{N}}
\left(u\left|\int_E g_k(x)\,d\mu(x)\right|-\rho_\varphi(g_k)\right)
\\
&\le\sup_{f\in L^\varphi(X,\mu)}
\left(u\left|\int_E f(x)\,d\mu(x)\right|-\rho_\varphi(f)\right).
\end{align*}
Therefore, passing to the supremum over all $f\in L^1(E,\mu)$ in 
the left-hand side above, we obtain the desired converse inequality 
``\eqref{eq:int-phi*} $\le$ \eqref{eq:int-phi*-0}.''
\end{proof}

\subsection{Generalized $N$-functions: the nicest among generalized $\Phi$-functions}

It is common practice in the study of Musielak--Orlicz spaces to restrict one's attention to
the so-called generalized $N$-functions within the class $\Phi(X,\mu)$. This special subclass
of $\Phi(X,\mu)$ consists of generalized $\Phi$-functions characterized by ``normal'' growth
in the second variable at zero and infinity. We adapt the definition of a generalized
$N$-function from~\cite[Definition~2.4.4, discussion on p.~53]{DHHR11}.
\begin{definition}\label{def:N}
We say that a function $\varphi\in\Phi(X,\mu)$ is a \emph{generalized $N$-function} 
on $(X,\mu)$, and write $\varphi\in N(X,\mu)$, if for every $x\in X$, the $\Phi$-function
$t\mapsto\varphi(x,t)$ is positive, finite-valued and satisfies
\[
\lim_{t\to0^+}\frac{\varphi(x,t)}{t}=0
\quad\text{and}\quad
\lim_{t\to\infty}\frac{\varphi(x,t)}{t}=\infty.
\]
\end{definition}
\begin{remark}\label{re:N-function}
A function $\varphi:X\times[0,\infty)\to[0,\infty)$ belongs to the class $N(X,\mu)$ if and 
only if it admits the representation 
\[
\varphi(x,t)=\int_0^t\varphi'(x,\tau)\,d\tau
\]
for all $x\in X$ and all $t\ge0$, where 
$\varphi':X\times[0,\infty)\to[0,\infty)$
is the right derivative of $\varphi$ with respect to the second variable, 
which satisfies, for any $x\in X$, the following properties:
\begin{enumerate}
    \item[(a)] $t\mapsto\varphi'(x,t)$ is non-decreasing and right-continuous;
    \item[(b)] $\varphi'(x,0)=0$;
    \item[(c)] $\varphi'(x,t)>0$ for $t>0$;
    \item[(d)] $\varphi'(x,t)\to\infty$ as $t\to\infty$.
\end{enumerate}
The necessity part of this statement follows 
from~\cite[Chapter~I, \S1, Subsection~5]{KR61}. The sufficiency part can be 
established by repeating the arguments 
from~\cite[Chapter~I, \S1, Subsection 4]{KR61}. 
\end{remark}

The integral representation of a function $\varphi\in N(X,\mu)$ offers a convenient way 
to express its conjugate $\varphi^*$. For this, define  
$(\varphi')^{-1}:X\times[0,\infty)\to[0,\infty)$ by 
\[
(\varphi')^{-1}(x,t):=\inf\{u\ge0\,:\,\varphi'(x,u)>t\}
\]
for all $x\in X$ and $t\ge0$. Logically, when $x\in X$ is fixed, the function
$t\mapsto(\varphi')^{-1}(x,t)$ is the \emph{right-continuous inverse} of the right
derivative $t\mapsto\varphi'(x,t)$. It is shown in~\cite[Theorem~2.6.8]{DHHR11} that if
$\varphi\in N(X,\mu)$, then $\varphi^*\in N(X,\mu)$ and $(\varphi^*)'=(\varphi')^{-1}$. Thus
for all $x\in X$ and $t\ge0$, the conjugate function $\varphi^*$ can be 
written as
\[
\varphi^*(x,t)=\int_0^t (\varphi')^{-1}(x,\tau)\,d\tau.
\]

Here are several useful inequalities for a generalized $N$-function, involving its
right derivative with respect to the second variable and its conjugate functions.
\begin{lemma}\label{le:N-ineq}
Let $\varphi\in N(X,\mu)$. The following inequalities hold:
\begin{enumerate}
    \item[(a)] $\frac{t}2\,\varphi'(x,\frac{t}2)\le\varphi(x,t)\le t\,\varphi'(x,t)$
    for all $x\in X$ and $t\ge0$;
    \item[(b)] $\varphi(x,\frac{\varphi^*(x,t)}{t})<\varphi^*(x,t)
    \le\varphi(x,\frac{2\varphi^*(x,t)}{t})$
    for all $x\in X$ and $t>0$.
\end{enumerate}
If, moreover, $\varphi$ satisfies the $\Delta_2$-condition, then
\begin{enumerate}
    \item[(c)] $\varphi(x,t)\approx t\,\varphi'(x,t)$
    uniformly in $x\in X$ and $t\ge0$;
    \item[(d)] $\varphi(x,\frac{\varphi^*(x,t)}{t})\approx\varphi^*(x,t)$
    uniformly in $x\in X$ and $t>0$.
\end{enumerate}
\end{lemma}
\begin{proof}
Inequality~(a) is proved in~\cite[Lemma~2.6.6]{DHHR11}. Statement~(c) is a consequence of 
part~(a) when $\varphi$ satisfies the $\Delta_2$-condition, as noted 
in~\cite[Remark~2.6.7]{DHHR11}. Indeed, if $D\ge2$ is the $\Delta_2$-constant
of $\varphi$, then it follows from inequality~(a) that 
\[
\varphi(x,t)\le t\,\varphi'(x,t)\le\varphi(x,2t)\le D\,\varphi(x,t)
\]
holds for all $x\in X$ and $t\ge0$, hence $\varphi(x,t)\approx t\,\varphi'(x,t)$
uniformly in $x\in X$ and $t\ge0$.

Inequality~(b) is a corollary of the well-known inequality
\begin{equation}\label{eq:Krasn-Rut}
\tau<\varphi^{-1}(x,\tau)\,(\varphi^*)^{-1}(x,\tau)\le2\tau
\end{equation}
for all $x\in X$ and $\tau>0$, see~\cite[Formula~(2.10)]{KR61}. Here 
$\varphi^{-1}(x,\tau)$ denotes, for all $x\in X$, the inverse of 
the function $\tau\mapsto\varphi(x,\tau)$, which is increasing due to $\varphi\in N(X,\mu)$. 
Making the change of variables $\tau=\varphi^*(x,t)$ with $x\in X$ and $t>0$ 
in~\eqref{eq:Krasn-Rut}, we obtain
\[
\frac{\varphi^*(x,t)}{t}<\varphi^{-1}(x,\varphi^*(x,t))
\le\frac{2\varphi^*(x,t)}{t}.
\]
Applying $\varphi(x,\cdot)$ to all three sides of the above inequality proves part~(b). 
Finally, the $\Delta_2$-condition for $\varphi$ allows us to extend statement~(b) to
\[
\varphi\left(x,\frac{\varphi^*(x,t)}{t}\right)<\varphi^*(x,t)
\le\varphi\left(x,\frac{2\varphi^*(x,t)}{t}\right)\le 
D\,\varphi\left(x,\frac{\varphi^*(x,t)}{t}\right),
\]
for all $x\in X$ and $t>0$, and thus equivalence~(d) follows.
\end{proof}

\subsection{Variable Lebesgue spaces as Musielak--Orlicz spaces}
\label{subseq:VLS}

Variable Lebesgue spaces $L^{p(\cdot)}(X,\mu)$ generalize 
the classical Lebesgue spaces $L^p(X,\mu)$ by allowing the exponent $p$ to vary from point 
to point. Instead a fixed exponent $p\in[1,\infty]$, one
considers a \emph{variable exponent}, that is, a $\mu$-measurable function 
\[
p(\cdot)\ :\ X\to[1,\infty].
\]
The set of all variable exponents on $(X,\mu)$ is traditionally denoted by
$\mathcal{P}(X,\mu)$. Given $p(\cdot)\in\mathcal{P}(X,\mu)$, the 
\emph{dual variable exponent}
$p'(\cdot)\in\mathcal{P}(X,\mu)$ is defined by the formula
\[
\frac{1}{p(x)}+\frac{1}{p'(x)}=1,\quad x\in X,
\]
with the convention that $1/\infty=0$.

Each exponent $p(\cdot)\in\mathcal{P}(X,\mu)$ determines two generalized 
$\Phi$-functions on $(X,\mu)$, which induce the space $L^{p(\cdot)}(X,\mu)$ as a 
Musielak--Orlicz space. These are the power-type functions $\bar{\varphi}_{p(\cdot)}$ and 
$\widetilde{\varphi}_{p(\cdot)}$, defined for all $x\in X$ and $t\ge0$ by
\[
\bar{\varphi}_{p(\cdot)}(x,t):=
\left\{\begin{array}{ll}
     t^{p(x)} & \text{if }p(x)<\infty,  \\
     0 & \text{if }p(x)=\infty\text{ and }t\in[0,1], \\
     \infty & \text{if }p(x)=\infty\text{ and }t\in(1,\infty)
\end{array}
\right.
\]
and
\[
\widetilde{\varphi}_{p(\cdot)}(x,t):=
\left\{\begin{array}{ll}
     \frac{t^{p(x)}}{p(x)} & \text{if }p(x)<\infty,  \\
     0 & \text{if }p(x)=\infty\text{ and }t\in[0,1], \\
     \infty & \text{if }p(x)=\infty\text{ and }t\in(1,\infty).
\end{array}
\right.
\]
Indeed, one easily checks by Definition~\ref{def:gen-Phi} that
$\bar{\varphi}_{p(\cdot)},\widetilde{\varphi}_{p(\cdot)}\in\Phi(X,\mu)$. The corresponding
Musielak--Orlicz spaces $L^{\bar{\varphi}_{p(\cdot)}}(X,\mu)$ 
and $L^{\widetilde{\varphi}_{p(\cdot)}}(X,\mu)$
coincide up to the equivalence of norms: for all 
$f\in L^{\bar{\varphi}_{p(\cdot)}}(X,\mu)=L^{\widetilde{\varphi}_{p(\cdot)}}(X,\mu)$, 
there holds
\[
\|f\|_{\widetilde{\varphi}_{p(\cdot)}}\le\|f\|_{\bar{\varphi}_{p(\cdot)}}
\le2\,\|f\|_{\widetilde{\varphi}_{p(\cdot)}},
\]
see~\cite[formula~(3.2.2)]{DHHR11}. In view of this remarkable fact, we may use the
notation $\varphi_{p(\cdot)}$ to refer to either $\bar{\varphi}_{p(\cdot)}$
or $\widetilde{\varphi}_{p(\cdot)}$, and define the \emph{variable Lebesgue space}
$L^{p(\cdot)}(X,\mu)$ as the Musielak--Orlicz space $L^{\varphi_{p(\cdot)}}(X,\mu)$
with the norm $\|\cdot\|_{p(\cdot)}:=\|\cdot\|_{\varphi_{p(\cdot)}}$. Since the aim of our
work is to characterize the boundedness of the maximal operator $M$ on 
$L^{p(\cdot)}(X,\mu)$, the exact norm of the space $L^{p(\cdot)}(X,\mu)$ is not important
and may be taken to be either of the above equivalent norms.

Note that both $\varphi_{p(\cdot)}:=\bar{\varphi}_{p(\cdot)}$ and 
$\varphi_{p(\cdot)}:=\widetilde{\varphi}_{p(\cdot)}$
are \emph{proper} if $p(\cdot)\in\mathcal{P}(X,\mu)$; 
see~\cite[Theorem~3.2.13]{DHHR11} for the proof. To ensure other desirable properties of 
$\varphi_{p(\cdot)}$, we will hereafter restrict ourselves to the case when the 
variable exponent $p(\cdot)$ is bounded away from one and infinity, that is, when
\[
p_-:=\operatornamewithlimits{ess\,inf}_{x\in X}p(x)>1
\quad\text{and}\quad
p_+:=\operatornamewithlimits{ess\,sup}_{x\in X}p(x)<\infty.
\]
Since two variable exponents that differ only on a set of measure zero generate the same
variable Lebesgue space, for our purposes we may, without loss of generality, assume that 
every exponent $p(\cdot)\in\mathcal{P}(X,\mu)$ with $1<p_-\le p_+<\infty$ satisfies 
$1<p(x)<\infty$ 
for all $x\in X$. Therefore, in what follows, let the condition $1<p_-\le p_+<\infty$ 
implicitly include this convention. Under the assumptions made, the functions 
$\bar{\varphi}_{p(\cdot)}$ and $\widetilde{\varphi}_{p(\cdot)}$
behave well in the following respects.
\begin{lemma}\label{le:phi-bar-tilde}
Suppose that $p(\cdot)\in\mathcal{P}(X,\mu)$ satisfies $1<p_-\le p_+<\infty$. Then uniformly in 
$x\in X$ and $t\ge0$, there hold:
\begin{enumerate}
    \item[(a)] $\bar{\varphi}_{p(\cdot)}(x,t)\approx\widetilde{\varphi}_{p(\cdot)}(x,t)$;
    \item[(b)] $\widetilde{\varphi}^*_{p(\cdot)}(x,t)
    =\widetilde{\varphi}_{p'(\cdot)}(x,t)\approx\bar{\varphi}^*_{p(\cdot)}(x,t)$.
\end{enumerate}
Moreover, $\varphi_{p(\cdot)}\in N(X,\mu)$ and both $\varphi_{p(\cdot)}$ and
$\varphi^*_{p(\cdot)}$ satisfy the $\Delta_2$-condition.
\end{lemma}
\begin{proof}
Since $p_+<\infty$, we have $\bar{\varphi}_{p(\cdot)}(x,t)=t^{p(x)}$ and 
$\widetilde{\varphi}_{p(\cdot)}(x,t)=t^{p(x)}/p(x)$
for all $x\in X$ and $t\ge0$. By checking the conditions of Definition~\ref{def:N}, 
one readily sees that $\bar{\varphi}_{p(\cdot)}$ and $\widetilde{\varphi}_{p(\cdot)}$ are 
generalized $N$-functions on $(X,\mu)$ due to $1<p(x)<\infty$ for all $x\in X$. It is also
immediate that both $\bar{\varphi}_{p(\cdot)}$ and $\widetilde{\varphi}_{p(\cdot)}$ satisfy
the $\Delta_2$-condition with constant $2^{p_+}$. Part~(a) follows since 
for all $x\in X$ and $t\ge0$, there holds
\[
\frac{1}{p_+}\,\bar{\varphi}_{p(\cdot)}(x,t)
\le\widetilde{\varphi}_{p(\cdot)}(x,t)\le
\bar{\varphi}_{p(\cdot)}(x,t).
\]

It is shown in~\cite[Example~(b) on p.~56]{DHHR11} that 
$\widetilde{\varphi}^*_{p(\cdot)}=\widetilde{\varphi}_{p'(\cdot)}$, where 
$p'(\cdot)\in\mathcal{P}(X,\mu)$ is the dual variable exponent. Note that the condition 
$1<p_-\le p_+<\infty$ implies, via direct computation, that
\[
(p')_+=\frac{p_-}{p_--1}<\infty
\quad\text{and}\quad
(p')_-=\frac{p_+}{p_+-1}>1.
\]
Then calculations in~\cite[Example~(c) on p.~56]{DHHR11},
conditions $1<p_-\le p_+<\infty$ and $1<(p')_-\le(p')_+<\infty$, and part~(a) together yield
\[
\bar{\varphi}^*_{p(\cdot)}(x,t)=\frac{p(x)-1}{p(x)^{p'(x)}}\,
\bar{\varphi}_{p'(\cdot)}(x,t)\approx \bar{\varphi}_{p'(\cdot)}(x,t)
\approx\widetilde{\varphi}_{p'(\cdot)}(x,t)
\]
uniformly in $x\in X$ and $t\ge0$, which completes the proof of part~(b).
Since $(p')_+<\infty$, the function $\varphi_{p'(\cdot)}$ satisfies the $\Delta_2$-condition.
Thus $\varphi^*_{p(\cdot)}$ also satisfies the $\Delta_2$-condition due to the relations in 
part~(b).
\end{proof}

\section{Spaces of Homogeneous Type}\label{sec:SHT}
Let us now pass from abstract measure spaces to the specific spaces that form the foundation
of our work: spaces of homogeneous type. As a main analytic tool on these spaces, we present
Hyt\"onen and Kairema's dyadic system with a distinguished center point 
(Theorem~\ref{th:HK})---precisely the system of cubes which will ``do us good,'' true to
the title. Then the maximal, dyadic maximal, and averaging operators are defined
(Definitions~\ref{def:MO} and \ref{def:averaging-operators}), and we give a version of
the Calder\'on--Zygmund lemma in spaces of homogeneous type (Lemma~\ref{le:CZ}).

\subsection{Definition of spaces of homogeneous type}
Spaces of homogeneous type, if briefly introduced, are quasi-metric spaces with 
a doubling measure. This notion first appeared in the works of Coifman and Weiss
\cite{CW71,CW77}. 

A more precise definition of spaces of homogeneous type is given below.
We take as our starting point the definition of Coifman and Weiss~\cite[pp.~587--588]{CW77},
while also adopting the useful approach of Alvarado and Mitrea, who interpret
a Borel measure as a measure defined on any $\sigma$-algebra containing the 
Borel $\sigma$-algebra; we refer the reader to~\cite[Definitions~2.9 and 3.2]{AM15}. 
The latter helps us to include the Euclidean space $\mathbb{R}^n$ with the Lebesgue measure 
as an example of a space of homogeneous type. Note also that, unlike other common 
definitions of spaces of homogeneous type, our Definition~\ref{def:SHT} requires the 
measure $\mu$ to be complete---so that we can consider Musielak--Orlicz spaces 
over spaces of homogeneous type and make use of the properties given in 
Section~\ref{se:MO}.

\begin{definition}\label{def:SHT}
Given a set $X$, we say that a triple $(X,d,\mu)$ is a {\it space of homogeneous
type} if
\begin{enumerate}
    \item[(i)] $d:X\times X\to[0,\infty)$ is a quasi-metric on $X$, i.e., the following hold:
    \begin{enumerate}
        \item[(a)] $d(x,y)=0$ if and only if $x=y$;
        \item[(b)] $d(x,y)=d(y,x)$ for all $x,y\in X$;
        \item[(c)] (\textit{quasi-triangle inequality}) there exists a constant $A_0\ge1$ such that 
        \[
        d(x,z)\le A_0\,(d(x,y)+d(y,z))
        \quad\text{for all }x,y,z\in X.
        \]
    \end{enumerate}
    \item[(ii)] $\mu$ is a complete Borel measure with respect to the topology
    $\tau_d$ induced by the quasi-metric $d$, or more precisely,
    \begin{enumerate}
        \item[(a)] a set $\Omega\subseteq X$ is open, $\Omega\in\tau_d$, if for each
        $x\in\Omega$ there exists $r>0$ such that the {\it quasi-metric ball} 
        \[
        B(x,r):=\{y\in X\ :\ d(x,y)<r\},
        \]
        centered at $x$ and of radius $r$, is contained in $\Omega$;
        \item[(b)] $\mu$ is a nonnegative, complete measure defined on a $\sigma$-algebra
        containing the Borel $\sigma$-algebra on $(X,\tau_d)$.
    \end{enumerate}
    \item[(iii)] Quasi-metric balls are $\mu$-measurable, and there exists a constant
    $A=A_\mu\ge1$ such that, for any $x\in X$ and $r>0$, one has
    \[
    0<\mu(B(x,r))\le A\,\mu(B(x,r/2))<\infty,
    \]
    i.e., $\mu$ is a \emph{doubling measure}. The smallest such constant $A$ is 
    called the \emph{doubling constant}.
\end{enumerate}
\end{definition}

Notice that condition~(iii) implies that 
there exists a positive constant $C=C(A,A_0)$ such that for all $x\in X$, $0<r<R$ and
$y\in B(x,R)$,
\begin{equation}\label{eq:LMB}
\frac{\mu(B(y,r))}{\mu(B(x,R))}\ge C\left(\frac rR\right)^{\log_2A}.
\end{equation}
This {\it lower mass bound} is well-known and stated, e.g., in~\cite[Lemma~2.3]{CUS18}. 
The condition that balls have finite measure ensures that $\mu$ is $\sigma$-finite, since 
for any $x\in X$, we can represent $X$ as the countable union of the balls
$\bigcup_{n=1}^\infty B(x,n)$. The assumption that balls have positive measure helps to 
avoid trivial measures.

Remarkably, quasi-metric balls need {\it not} be open in the topology induced by the 
quasi-metric, as shown in a simple example by Hyt\"onen and Kairema~\cite[p.~5]{HK12}.
Thus, the requirement in~(iii) that balls are $\mu$-measurable is not redundant. An 
alternative way to ``measure balls with a Borel measure,'' however, would be by using the 
classical result by Mac\'{i}as and Segovia~\cite[Theorem~2]{MS79} who proved that there is
always an equivalent quasi-metric in which all balls are open. Some works display this 
approach~\cite{HPR12,KS25}.

Another well-known fact is that a space of homogeneous type $(X,d,\mu)$ has finite
measure $\mu(X)<\infty$ if and only if $X=B(x,r)$ for some $x\in X$ and $r>0$, 
see~\cite[Lemma~1.9]{BC96}. In this paper, we will mainly work with \emph{unbounded} 
spaces of homogeneous type, which is equivalent to assuming that $\mu(X)=\infty$.

On a space of homogeneous type $(X,d,\mu)$, we naturally define the space
$L^s_\textnormal{loc}(X,d,\mu)$ for $s\in[1,\infty)$ by
\[
L^s_\textnormal{loc}(X,d,\mu):=\left\{f\in L^0(X,d,\mu)\ :\ \int_B|f(y)|^s\,d\mu(y)<\infty
\;\text{ for all balls }B\subseteq X\right\}.
\]
Simply put, locally integrable functions on $X$ are those with
the finite integrals over quasi-metric balls.

\subsection{The Hyt\"onen--Kairema dyadic system}
A fundamental role in harmonic analysis on the Euclidean space $\mathbb{R}^n$
belongs to the standard system of dyadic cubes 
\[
\{2^{-k}([0,1)^n+\mathbf{m})\ :\ k\in\mathbb{Z},\,\mathbf{m}\in\mathbb{Z}^n\}. 
\]
In particular, dyadic cubes allow one to define dyadic versions of classical operators,
such as the Hardy--Littlewood maximal operator, and employ powerful decomposition techniques,
such as the Calder\'on--Zygmund decomposition. A refined dyadic approach suggests the use of
\emph{several adjacent dyadic grids} instead of just one fixed system. Typically, 
one works with $3^n$ grids obtained by shifting the standard dyadic system in each 
coordinate direction by $0$, $1/3$ or $2/3$ of the sidelength---with the shifts depending on 
the scale so that the hierarchical structure of the cubes is preserved.
This collection of $1/3$-shifted dyadic grids has fine covering 
properties and is used, for example, to dominate the maximal 
function by a sum of dyadic maximal functions (see~\cite[Section~3]{LN19}
and~\cite[Section~7.3.2]{CP19} for a quick introduction to the one-third trick). 

In the influential 2012 paper~\cite{HK12}, Hyt\"onen and Kairema continued 
earlier work of Christ~\cite[Section~3]{MC90} and Hyt\"onen and 
Martikainen~\cite[Section~4]{HM12}---and constructed a system of 
finitely many adjacent grids of ``dyadic cubes'' with useful covering properties on
geometrically doubling quasi-metric spaces. Since spaces of homogeneous type
in the sense of Definition~\ref{def:SHT} are geometrically doubling 
(see~\cite[p.~67]{CW71} and \cite[p.~12]{AS25}), Hyt\"onen 
and Kairema's result made the dyadic tool---originally developed in the Euclidean 
setting---available for analysis on quasi-metric spaces with doubling measure.

We present a version of the Hyt\"onen--Kairema dyadic 
construction, which originates in \cite[Theorems~2.2 and 4.1, 
Corollary~7.4]{HK12}, 
for spaces of homogeneous type. In this formulation, we omit details which are not 
necessary for our work.
\begin{theorem}\label{th:HK}
Given a space of homogeneous type $(X,d,\mu)$,
there exists a {\rm system $\mathscr{D}$ of dyadic cubes on $X$} with parameters 
$0<\varepsilon,\eta<1$ and $0<r_0\le R_0<\infty$ 
depending on the space. The system 
\[
\mathscr{D}:=\bigcup_{t=1}^K\mathscr{D}^t
\]
consists of finitely many families $\mathscr{D}^t$, $t\in\{1,\ldots,K\}$, called 
{\rm the adjacent dyadic grids}. Each $\mathscr{D}^t$ is a countable collection of 
Borel sets $Q\subseteq X$ (with respect to the topology $\tau_d$), called {\rm dyadic cubes},
with the following properties:
\begin{itemize}
    \item[(a)] \textnormal{(generations)} cubes in $\mathscr{D}^t$ are organized into 
    generations, i.e., $\mathscr{D}^t:=\bigcup_{k\in\mathbb{Z}}\mathscr{D}_k^t$, where
    subcollection $\mathscr{D}_k^t$ is \emph{the $k$-th generation of dyadic cubes in the 
    grid $\mathscr{D}^t$}. Each generation $k\in\mathbb{Z}$ forms a partition of $X$:
    \[
    X=\bigcup_{Q\in\mathscr{D}_k^t}Q
    \quad\text{(disjoint union)};
    \]
    \item[(b)] \textnormal{(nestedness)} there is no partial overlap across generations:
    if $k>l$, then for any $Q_1\in\mathscr{D}^t_l$ and $Q_2\in\mathscr{D}^t_k$, either
    $Q_2\subseteq Q_1$ (in which case we say that $Q_1$ is \emph{an ancestor} of $Q_2$) or 
    $Q_1\cap Q_2=\emptyset$;
    \item[(c)] \textnormal{(children and parent)} for any $Q_1\in\mathscr{D}_k^t$,
    there exists at least one cube $Q_2\in\mathscr{D}^t_{k+1}$ 
    (called \emph{a child} of $Q_1$) such that $Q_2\subseteq Q_1$, 
    and there exists exactly one cube $Q_3\in\mathscr{D}^t_{k-1}$
    (called \emph{the parent} of $Q_1$) such that 
    $Q_3\supseteq Q_1$;
    \item[(d)] \textnormal{(nondegeneracy of children)} if $Q_2\in\mathscr{D}^t_k$ is a 
    child of $Q_1\in\mathscr{D}^t_{k-1}$, then 
    \[
    \mu(Q_2)\ge\varepsilon\,\mu(Q_1);
    \]
    \item[(e)] \textnormal{(inner/outer ball)} for any $k\in\mathbb{Z}$
    and $Q\in\mathscr{D}^t_k$, there is a point $x_c(Q)\in Q$, called 
    \emph{the center} of $Q$, such that 
    \[
    B(x_c(Q),r_0\eta^k)\subseteq Q\subseteq B(x_c(Q),R_0\eta^k)=:B_Q.
    \]
    The number $\eta^k$ is interpreted as the ``sidelength'' of 
    every cube $Q$ in the $k$-th generation of any grid $\mathscr{D}^t$.
\end{itemize}
Due to several adjacent grids in its structure, the dyadic system $\mathscr{D}$ has a 
fine covering property: for every ball $B:=B(x,r)\subseteq X$,
there exists a dyadic cube $Q_B\in\mathscr{D}$ such that 
\[
B\subseteq Q_B
\quad\text{and}\quad
\mu(Q_B)\le C\,\mu(B),
\]
where $C<\infty$ is a constant independent of $x$ and $r$.
\end{theorem}
In the above version of Hyt\"onen and Kairema's result, we have also included two direct
corollaries of the canonical statement in~\cite[Theorem~4.1]{HK12}. First, we 
explicitly indicated that cubes $Q\in\mathscr{D}$ are Borel sets~\cite[Remark~4.2]{AHT17}.
The second addition is property~(d), which logically extends the doubling 
property of the measure from balls to dyadic cubes (see \cite[Corollary~2.9]{AW18} or
\cite[Lemma~8]{AK19}).

Note that the construction in Theorem~\ref{th:HK} allows the situation when a cube has 
just one child, itself, and this can last for several generations. In other words, it is
possible that for a cube $Q\in\mathscr{D}^t$, there exist more than one $k\in\mathbb{Z}$ so
that $Q\in\mathscr{D}^t_k$. To avoid ambiguity, we may only speak of a child or the parent 
of $Q\in\mathscr{D}^t$ with respect to a specific $k\in\mathbb{Z}$ such that 
$Q\in\mathscr{D}^t_k$, as we do this, for example, when stating property~(d).
\begin{remark}\label{re:Euclidean}
The Euclidean space $\mathbb{R}^n$ with the Euclidean distance and the 
$n$-dimensional Lebesgue measure $|\cdot|$ is a space of homogeneous type in the sense of 
Definition~\ref{def:SHT}. Consequently, Theorem~\ref{th:HK} guarantees the existence of 
numerous dyadic systems on $\mathbb{R}^n$. These include the system
$\mathscr{D}:=\bigcup_{\mathbf{t}\in\mathcal{K}}\mathscr{D}^{\mathbf{t}}$ formed by the
classical collection of 1/3-shifted dyadic grids
\[
\mathscr{D}^{\mathbf{t}}:=\left\{
2^{-k} \big([0,1)^n+\mathbf{m}+(-1)^k\mathbf{t}\big)
\ :\ 
k\in\mathbb{Z},\,\mathbf{m}\in\mathbb{Z}^n\right\},
\quad\mathbf{t}\in\mathcal{K},
\]
where $\mathcal{K}:=\{0,\frac13,\frac23\}^n$ denotes the set of $3^n$ translation vectors
$\mathbf{t}:=(t_1,\ldots,t_n)$ with
$t_i\in\{0,\frac13,\frac23\}$ for all $i\in\{1,\ldots,n\}$. Such adjacent 
dyadic grids have been exploited, e.g., in the works of Muscalu, Tao and 
Thiele~\cite{MTT02}, Hyt\"onen, Lacey and P\'erez~\cite{HLP13}, and Lerner~\cite{AL17}.

To verify that the above system $\mathscr{D}$ fits into the construction of 
Theorem~\ref{th:HK}, let us check first that each grid $\mathscr{D}^\mathbf{t}$,
$\mathbf{t}\in\mathcal{K}$, satisfies individual properties~(a)--(e) listed in 
Theorem~\ref{th:HK}:

\medskip
(a) (\emph{generations}) Cubes in $\mathscr{D}^{\mathbf{t}}$ are organized into 
generations, i.e., $\mathscr{D}^\mathbf{t}:=
\bigcup_{k\in\mathbb{Z}}\mathscr{D}^\mathbf{t}_k$,
where the subcollection 
\[
\mathscr{D}^\mathbf{t}_k:=
\left\{Q_{k,\mathbf{m}}^{\mathbf{t}}:=
2^{-k} \big([0,1)^n+\mathbf{m}+(-1)^k\mathbf{t}\big)
\ :\ \mathbf{m}\in\mathbb{Z}^n\right\}
\]
is the $k$-th generation of half-open cubes in the grid $\mathscr{D}^\mathbf{t}$. 
Observe that a cube $Q_{k,\mathbf{m}}^{\mathbf{t}}$ has the center 
$\mathbf{x}_c(Q_{k,\mathbf{m}}^{\mathbf{t}}):=2^{-k}(\frac{\mathbf{e}}2+\mathbf{m}
+(-1)^k\mathbf{t})$, where $\mathbf{e}:=(1,\ldots,1)\in\mathbb{R}^n$ is the all-ones
vector, and the sidelength $\ell(Q_{k,\mathbf{m}}^{\mathbf{t}}):=2^{-k}$. It easily follows
that the cubes in $\mathscr{D}^\mathbf{t}_k$ are pairwise disjoint. Indeed, take two vectors
$\mathbf{m}^1=(m_1^1,\ldots,m_n^1)$ and $\mathbf{m}^2=(m_1^2,\ldots,m_n^2)$ such that
$\mathbf{m}^1\ne\mathbf{m}^2$. Then $m_i^1\ne m_i^2$ for some $i\in\{1,\ldots,n\}$, 
so the absolute value of the difference between the $i$-th coordinates of the centers 
$\mathbf{x}_c(Q_{k,\mathbf{m}^1}^{\mathbf{t}})$ and 
$\mathbf{x}_c(Q_{k,\mathbf{m}^2}^{\mathbf{t}})$ is at 
least $2^{-k}$. Since also $\ell(Q_{k,\mathbf{m}^1}^{\mathbf{t}})=
\ell(Q_{k,\mathbf{m}^2}^{\mathbf{t}})=2^{-k}$ and the cubes are half-open, we conclude
that $Q_{k,\mathbf{m}^1}^{\mathbf{t}}\cap Q_{k,\mathbf{m}^2}^{\mathbf{t}}=\emptyset$.
Furthermore,
\[
\mathbb{R}^n=\bigcup_{Q\in\mathscr{D}^\mathbf{t}_k}Q=\bigcup_{\mathbf{m}\in\mathbb{Z}^n}
Q^\mathbf{t}_{k,\mathbf{m}}
\quad\textnormal{(disjoint union)},
\]
that is, each generation $k\in\mathbb{Z}$ forms a partition of $\mathbb{R}^n$. To check 
this, take $\mathbf{x}=(x_1,\ldots,x_n)\in\mathbb{R}^n$, and for each $i\in\{1,\ldots,n\}$, 
find a unique number $m_i\in\mathbb{Z}$ such that 
\[
x_i\in2^{-k}\big([0,1)+m_i+(-1)^kt_i\big).
\]
Then $\mathbf{x}\in Q^\mathbf{t}_{k,\mathbf{m}}$ with 
$\mathbf{m}:=(m_1,\ldots,m_n)\in\mathbb{Z}^n$, which proves our claim.

\medskip
(b) (\emph{nestedness}) There is no partial overlap between generations:
any two cubes in the grid $\mathscr{D}^\mathbf{t}$ are either disjoint or 
one is contained in the other. Indeed, let 
$k>l$, and let cubes $Q^\mathbf{t}_{k,\mathbf{m}}\in\mathscr{D}^{\mathbf{t}}_k$
and $Q^\mathbf{t}_{l,\mathbf{m}'}\in\mathscr{D}^{\mathbf{t}}_l$ be such that 
$Q^\mathbf{t}_{k,\mathbf{m}}\cap Q^\mathbf{t}_{l,\mathbf{m}'}\ne\emptyset$. We will prove
that then $Q^\mathbf{t}_{k,\mathbf{m}}\subseteq Q^\mathbf{t}_{l,\mathbf{m}'}$. 
First, let us check that for a cube 
$Q^\mathbf{t}_{k,\mathbf{m}}\in\mathscr{D}^{\mathbf{t}}_k$, we 
can always find a larger-scale cube $Q^\mathbf{t}_{k-1,\widetilde{\mathbf{m}}}\in
\mathscr{D}^\mathbf{t}_{k-1}$ covering $Q^\mathbf{t}_{k,\mathbf{m}}$. Take a point
$\mathbf{x}\in Q^\mathbf{t}_{k,\mathbf{m}}$. Since the cubes from 
$\mathscr{D}^\mathbf{t}_{k-1}$ partition $\mathbb{R}^n$, there exists a unique vector
$\widetilde{\mathbf{m}}\in\mathbb{Z}^n$ such that 
$\mathbf{x}\in Q^\mathbf{t}_{k-1,\widetilde{\mathbf{m}}}$. Decompose
\begin{align}
Q^\mathbf{t}_{k-1,\widetilde{\mathbf{m}}}&=
2^{-(k-1)}\big([0,1)^n+\widetilde{\mathbf{m}}+(-1)^{k-1}\mathbf{t}\big)
\nonumber
\\
&=2^{-k}\big([0,2)^n+2\widetilde{\mathbf{m}}+2(-1)^{k-1}\mathbf{t}\big)
\nonumber
\\
&=2^{-k}\Bigg(
\bigcup_{\mathbf{s}\in\{0,1\}^n}\big([0,1)^n+\mathbf{s}\big)
+2\widetilde{\mathbf{m}}+3(-1)^{k-1}\mathbf{t}+(-1)^k\mathbf{t}
\Bigg)
\nonumber
\\
&=\bigcup_{\mathbf{s}\in\{0,1\}^n}
2^{-k}\Big([0,1)^n+\big(2\widetilde{\mathbf{m}}+3(-1)^{k-1}\mathbf{t}+\mathbf{s}\big)
+(-1)^k\mathbf{t}\Big)
\nonumber
\\
&=\bigcup_{\mathbf{s}\in\{0,1\}^n}
Q^\mathbf{t}_{k,2\widetilde{\mathbf{m}}+3(-1)^{k-1}\mathbf{t}+\mathbf{s}}
\quad\text{(disjoint union)},
\label{eq:parent-decomposed}
\end{align}
taking into account that 
$(2\widetilde{\mathbf{m}}+3(-1)^{k-1}\mathbf{t}+\mathbf{s})\in\mathbb{Z}^n$
for every $\mathbf{s}\in\{0,1\}^n$, i.e., every $\mathbf{s}:=(s_1,\ldots,s_n)$ with
$s_i\in\{0,1\}$ for all $i\in\{1,\ldots,n\}$. Then 
$\mathbf{x}\in Q^\mathbf{t}_{k,2\widetilde{\mathbf{m}}+3(-1)^{k-1}\mathbf{t}+\mathbf{s}'}$
for some $\mathbf{s}'\in\{0,1\}^n$. Since also $\mathbf{x}\in Q^\mathbf{t}_{k,\mathbf{m}}$,
and different cubes in the same generation $\mathscr{D}^\mathbf{t}_k$ are disjoint by 
property~(a), we conclude that 
\[
Q^\mathbf{t}_{k,\mathbf{m}}=
Q^\mathbf{t}_{k,2\widetilde{\mathbf{m}}+3(-1)^{k-1}\mathbf{t}+\mathbf{s}'}
\subseteq Q^\mathbf{t}_{k-1,\widetilde{\mathbf{m}}}.
\]
Then, by induction, we can find $\mathbf{m}''\in\mathbb{Z}^n$ such that 
\[
Q^\mathbf{t}_{k,\mathbf{m}}
\subseteq Q^\mathbf{t}_{k-1,\widetilde{\mathbf{m}}}
\subseteq\ldots\subseteq Q^\mathbf{t}_{l,\mathbf{m}''}.
\]
By assumption $Q^\mathbf{t}_{k,\mathbf{m}}\cap Q^\mathbf{t}_{l,\mathbf{m}'}\ne\emptyset$,
therefore
$Q^\mathbf{t}_{l,\mathbf{m}''}\cap Q^\mathbf{t}_{l,\mathbf{m}'}\ne\emptyset$. Again by 
property~(a), it follows that $Q^\mathbf{t}_{l,\mathbf{m}''}= 
Q^\mathbf{t}_{l,\mathbf{m}'}$, and hence 
$Q^\mathbf{t}_{k,\mathbf{m}}\subseteq Q^\mathbf{t}_{l,\mathbf{m}'}$.

\medskip
(c) (\emph{children and parent}) In view of decomposition~\eqref{eq:parent-decomposed},
for any cube $Q_{k,\mathbf{m}}^\mathbf{t}\in\mathscr{D}_k^\mathbf{t}$, there exist exactly
$2^n$ cubes $Q_{k+1,2\mathbf{m}+3(-1)^k\mathbf{t}+\mathbf{s}}^\mathbf{t}
\in\mathscr{D}^\mathbf{t}_{k+1}$ with 
$\mathbf{s}\in\{0,1\}^n$ (the children of $Q_{k,\mathbf{m}}^\mathbf{t}$) such that
\[
Q_{k,\mathbf{m}}^\mathbf{t}:=\bigcup_{\mathbf{s}\in\{0,1\}^n}
Q_{k+1,2\mathbf{m}+3(-1)^k\mathbf{t}+\mathbf{s}}^\mathbf{t}.
\] 
Furthermore, for a cube $Q_{k,\mathbf{m}}^\mathbf{t}\in\mathscr{D}_k^\mathbf{t}$, there
exists exactly one cube 
$Q_{k-1,\widetilde{\mathbf{m}}}^\mathbf{t}\in\mathscr{D}^\mathbf{t}_{k-1}$ 
(the parent of $Q_{k,\mathbf{m}}^\mathbf{t}$) such that 
$Q_{k-1,\widetilde{\mathbf{m}}}^\mathbf{t}\supseteq Q_{k,\mathbf{m}}^\mathbf{t}$.
The existence of one such cube was proved in part~(b). The uniqueness of the parent follows 
from the pairwise disjointness of the cubes in generation $\mathscr{D}^\mathbf{t}_{k-1}$.

\medskip
(d) (\emph{nondegeneracy of children}) Let 
$Q_{k,\mathbf{m}}^\mathbf{t}\in\mathscr{D}^\mathbf{t}_k$ be a child of
$Q_{k-1,\widetilde{\mathbf{m}}}^\mathbf{t}\in\mathscr{D}^\mathbf{t}_{k-1}$. 
Since the latter cube has exactly $2^n$ children and all cubes in the same generation have
equal measure, it follows that
\[
|Q_{k,\mathbf{m}}^\mathbf{t}|=2^{-n}\,|Q_{k-1,\widetilde{\mathbf{m}}}^\mathbf{t}|.
\]

\medskip
(e) (\emph{inner/outer ball}) For any $k\in\mathbb{Z}$ and $\mathbf{m}\in\mathbb{Z}^n$,
there hold the inclusions
\[
B\left(\mathbf{x}_c(Q_{k,\mathbf{m}}^\mathbf{t}),\frac12\times2^{-k}\right)
\subseteq Q_{k,\mathbf{m}}^\mathbf{t}\subseteq
B\left(\mathbf{x}_c(Q_{k,\mathbf{m}}^\mathbf{t}),\frac{\sqrt{n}}2\times2^{-k}\right).
\]

\medskip\noindent
Then it only remains to check that the system 
$\mathscr{D}:=\bigcup_{\mathbf{t}\in\mathcal{K}}\mathscr{D}^{\mathbf{t}}$, comprising 
all the adjacent dyadic grids, has the fine covering property stated in Theorem~\ref{th:HK}.
Take a ball $B:=B(\mathbf{x},r)$ with the center $\mathbf{x}\in\mathbb{R}^n$ and 
radius $r>0$. Clearly, the ball $B$ is contained in the cube $Q\subseteq\mathbb{R}^n$
centered at $\mathbf{x}$ and with the sidelength $\ell(Q)=2r$. It is proved in
\cite[Lemma~2.5]{HLP13} that there exists $Q_B\in\mathscr{D}$ such that $Q\subseteq Q_B$
and $\ell(Q_B)\le6\,\ell(Q)$. Then $B\subseteq Q_B$ and 
\[
|Q_B|\le 6^n|Q|=6^n(2r)^n=
12^n\,\frac{\Gamma\left(\frac{n}2+1\right)}{\pi^{n/2}}|B|,
\]
where the constant in the final estimate is independent of $B$. 
\qed
\end{remark}

Here are some further properties of dyadic cubes related to their nesting.
\begin{lemma}\label{le:cones}
Let $(X,d,\mu)$ be a space of homogeneous type and 
$\mathscr{D}:=\bigcup_{t=1}^K\mathscr{D}^t$ be a dyadic system on $X$ defined in 
Theorem~\ref{th:HK}. For an $x\in X$, denote by $Q_k^t(x)\in\mathscr{D}^t_k$ 
\emph{the unique cube from the $k$-th generation of cubes in $\mathscr{D}^t$ that 
contains $x$}. Each sequence $\{Q_k^t(x)\}_{k\in\mathbb{Z}}$ is monotone, since
\[
Q_{k+1}^t(x)\subseteq Q^t_k(x) 
\quad\text{for all }k\in\mathbb{Z},
\]
and forms the \emph{$t$-dyadic cone} over $x$. The following hold for any 
$t\in\{1,\ldots,K\}$:
\begin{enumerate}
    \item[(a)] \emph{(cone expansion)} If $\mu(X)=\infty$, then for every $x\in X$, 
    one has
    \[
    \mu(Q^t_k(x))\to\infty 
    \quad\text{ as }k\to-\infty.
    \]
    \item[(b)] \emph{(measure gap)} If $Q^t_k(x)\varsubsetneq Q^t_l(x)$ for some $k>l$
    and $x\in X$, then 
    \[
    \frac{\mu(Q^t_k(x))}{\mu(Q^t_l(x))}\le1-\varepsilon,
    \]
    where $\varepsilon\in(0,1)$ is the parameter of $\mathscr{D}$ specified
    in Theorem~\ref{th:HK}(d).
    \item[(c)] \emph{(the Lebesgue differentiation theorem for cubes)} Let $\mu$ be 
    \emph{Borel-semiregular}, meaning that for every $E\subseteq X$ with $\mu(E)<\infty$, 
    there exists a Borel set $B\subseteq X$ (with respect to the topology $\tau_d$) 
    such that $\mu(E\triangle B)=0$.
    Given a function $f\in L^1_\textnormal{loc}(X,d,\mu)$, the limit
    \[
    \lim_{k\to+\infty}\frac{1}{\mu(Q^t_k(x))}\int_{Q^t_k(x)}f(y)\,d\mu(y)=f(x)
    \]
    holds for $\mu$-almost every $x\in X$. 
\end{enumerate}
\end{lemma}
\begin{proof}
The idea for the proof of property~(a) is due to Ryan Alvarado, co-author of~\cite{AM15}.
Fix an $x\in X$. By Theorem~\ref{th:HK}(e), for any $k\in\mathbb{Z}$ there holds
\[
B(x^t_k,r_0\eta^k)\subseteq Q_k^t(x)\subseteq B(x_k^t,R_0\eta^k),
\]
where $x^t_k:=x_c(Q^t_k(x))$ is the center of a cube $Q^t_k(x)$. Consequently, 
for any integer $k\le0$, we have
\[
x_0^t\in Q_0^t(x)\subseteq Q_k^t(x)\subseteq B(x_k^t,R_0\eta^k),
\]
which implies that $d(x_0^t,x_k^t)<R_0\eta^k$. For a fixed $k\le0$, consider 
$y\in B(x_0^t,\eta^k)$. Then, by the quasi-triangle inequality,
\begin{align*}
d(y,x_k^t)&\le A_0\,(d(y,x_0^t)+d(x_0^t,x_k^t))
\\
&<A_0\,(\eta^k+R_0\eta^k)
\\
&=A_0\,(R_0+1)\,\eta^k,
\end{align*}
hence $B(x_0^t,\eta^k)\subseteq B(x_k^t,A_0(R_0+1)\eta^k)$ for all $k\le0$. Note that by the
continuity of measure from below, we have
\[
\mu(B(x_0^t,\eta^k))\to\mu\bigg(\bigcup_{k\in\mathbb{Z}}B(x_0^t,\eta^k)\bigg)
=\mu(X)=\infty
\]
as $k\to-\infty$. Then this, the lower mass bound~\eqref{eq:LMB}, and the set inclusions 
established above, yield the desired limit
\begin{align*}
\mu(Q_k^t(x))&\ge\mu(B(x_k^t,r_0\eta^k))
\\
&\ge C\left(\frac{r_0}{A_0(R_0+1)}\right)^{\log_2 A}
\mu(B(x_k^t,A_0(R_0+1)\eta^k))
\\
&\ge C\left(\frac{r_0}{A_0(R_0+1)}\right)^{\log_2 A}
\mu(B(x_0^t,\eta^k))
\to\infty
\quad\text{as }k\to-\infty.
\end{align*}

To prove~(b), we may assume without loss of generality that 
$Q^t_k(x)\varsubsetneq Q^t_{k-1}(x)$. Otherwise, if $Q^t_k(x)$ could not be assigned to
such a generation $k\in\mathbb{Z}$, it would coincide with all its ancestors, contradicting
the assumption of claim~(b). Consequently, $k-1\ge l$ implies that 
$Q^t_{k-1}(x)\subseteq Q^t_l(x)$, so overall we have
\[
Q_k^t(x)\varsubsetneq Q^t_{k-1}(x)\subseteq Q^t_l(x).
\]
Due to $Q^t_k(x)\varsubsetneq Q^t_{k-1}(x)$, the cube $Q^t_{k-1}(x)\in\mathscr{D}^t_{k-1}$
is the parent of $Q^t_k(x)\in\mathscr{D}^t_k$ and of another child $Q\in\mathscr{D}^t_k$. By
Theorem~\ref{th:HK}(d), we have the estimate $\mu(Q)\ge\varepsilon\,\mu(Q^t_{k-1}(x))$. 
This gives the inequality
\[
\mu(Q^t_{k-1}(x))\ge\mu(Q^t_k(x))+\mu(Q)\ge\mu(Q^t_k(x))+\varepsilon\,\mu(Q^t_{k-1}(x)),
\]
which proves our claim that
\[
\mu(Q^t_k(x))\le(1-\varepsilon)\,\mu(Q^t_{k-1}(x))\le(1-\varepsilon)\,\mu(Q^t_l(x)).
\]

Property~(c) is an extension of the well-known fact that the Lebesgue differentiation 
theorem in its classical form, for concentric balls, holds in spaces of homogeneous type
with a Borel-semiregular measure~\cite[Theorem~3.14]{AM15}. The version for cubes has 
appeared in the literature (see~\cite[Corollary~2.6]{ACUM15} or 
\cite[Lemma~4.18]{AV14}), but without any regularity condition on the measure $\mu$---that
is, for the case when $\mu$ is defined \emph{only} on the $\sigma$-algebra of 
Borel sets, which is a less general setting than ours. We provide a slightly more explicit
proof here.

Fix an $x\in X$. Denoting again $x^t_k:=x_c(Q^t_k(x))$, we have from Theorem~\ref{th:HK}(e)
that for any $k\in\mathbb{Z}$, there holds
\[
B(x^t_k,r_0\eta^k)\subseteq Q_k^t(x)\subseteq B(x_k^t,R_0\eta^k).
\]
Since $x\in B(x^t_k,R_0\eta^k)$, it follows that 
$B(x_k^t,R_0\eta^k)\subseteq B(x,2A_0R_0\eta^k)$.
Thus
\[
B(x^t_k,r_0\eta^k)\subseteq Q_k^t(x)\subseteq B(x,2A_0R_0\eta^k),
\]
and the lower mass bound~\eqref{eq:LMB} gives the estimate
\[
\mu(Q^t_k(x))\ge\mu(B(x_k^t,r_0\eta^k))\ge 
\frac{1}{c}\,\mu(B(x,2A_0R_0\eta^k)),
\]
where $c:=\frac{1}{C}(\frac{2A_0R_0}{r_0})^{\log_2A}$. This and the Lebesgue 
differentiation theorem for concentric balls~\cite[Theorem~3.14(2)]{AM15} allow us 
to conclude that for $\mu$-almost every $x\in X$,
\begin{align*}
\frac{1}{\mu(Q^t_k(x))}&\int_{Q^t_k(x)}|f(y)-f(x)|\,d\mu(y)\\
&\le\frac{c}{\mu(B(x,2A_0R_0\eta^k))}\int_{B(x,2A_0R_0\eta^k)}|f(y)-f(x)|\,d\mu(y)
\to0
\end{align*}
as $k\to+\infty$, since the radii $2A_0R_0\eta^k\to0^+$ as $k\to+\infty$. This immediately
implies the claim of the theorem.
\end{proof}

From now on, when speaking of a space of homogeneous type $(X,d,\mu)$, we will always 
assume that there is a fixed dyadic system $\mathscr{D}$, chosen arbitrarily among 
the various systems provided by Theorem~\ref{th:HK}, associated with this space.

\subsection{Maximal, dyadic maximal and averaging operators}
The central object of our study is the Hardy--Littlewood maximal function on a space of
homogeneous type. Many of its key properties, however, can be derived by examining the
dyadic maximal function---a version of the maximal function associated with
Hyt\"onen and Kairema's dyadic cubes. To define these functions, we start with the
mean integral values (or simply, the means) of a measurable function.
\begin{definition}\label{def:means}
The {\it means} of a function $f\in L^0(X,d,\mu)$ over a $\mu$-measurable set
$E\subseteq X$ with $0<\mu(E)<\infty$ are defined by
\[
M_{s,E}f:=\left(\fint_E|f(y)|^s\,d\mu(y)\right)^{1/s}=
\left(\frac{1}{\mu(E)}\int_E|f(y)|^s\,d\mu(y)\right)^{1/s}
\]
for every $s\in[1,\infty)$.
The ``barred'' integral sign denotes the integral divided by the measure of its domain; 
we will frequently use this simplified notation. We also write $M_Ef$ instead of 
$M_{1,E}f$.
\end{definition}
Note that the means $M_{s,E}f$ are non-decreasing in $s$ as a consequence
of Jensen's inequality. Namely, one has $M_{r,E}f\le M_{s,E}f$ for any $1\le r<s$ due to
\[
(M_{r,E}f)^s=\left(\fint_E|f(y)|^r\,d\mu(y)\right)^{s/r}
\le\fint_E|f(y)|^s\,d\mu(y)=(M_{s,E}f)^s.
\]
The common idea behind maximal functions of a measurable function $f$ is to find the
``maximum'' value among the smallest ($s=1$) means of $f$ over some neighborhoods of 
a fixed point $x\in X$. If these neighborhoods of $x$ are balls, we speak of the classical
maximal function; if they are dyadic cubes, we obtain dyadic maximal functions.
\begin{definition}\label{def:MO}
Let $(X,d,\mu)$ be a space of homogeneous type and 
$\mathscr{D}:=\bigcup_{t=1}^K\mathscr{D}^t$ be a dyadic system 
associated with it. Given a function $f\in L^0(X,d,\mu)$, 
\begin{enumerate}
    \item[(i)] the {\it Hardy--Littlewood maximal function} $Mf$ is defined by 
    \[
    Mf(x):=\sup_{B\ni x}M_Bf,
    \quad x\in X,
    \]
    where the supremum is taken over all quasi-metric balls $B\subseteq X$ containing $x$;
    \item[(ii)] the \emph{dyadic maximal function} $M^\mathscr{D}f$ is defined by
    \[
    M^{\mathscr{D}}f(x):=\sup_{\substack{Q\ni x,\\Q\in\mathscr{D}}}M_Qf,
    \quad x\in X,
    \]
    where the supremum is taken over all dyadic cubes in $\mathscr{D}$ containing $x$;
    \item[(iii)] the \emph{$t$-dyadic maximal function} $M^{\mathscr{D}^t}f$
    is defined for any $t\in\{1,\ldots,K\}$ by using dyadic cubes of the grid 
    $\mathscr{D}^t$ only, that is,
    \[
    M^{\mathscr{D}^t}f(x):=\sup_{\substack{Q\ni x,\\Q\in\mathscr{D}^t}}M_{Q}f,
    \quad x\in X.
    \]
\end{enumerate} 
The {\it maximal operator} $M$, the \emph{dyadic maximal operator} $M^\mathscr{D}$
and the \emph{$t$-dyadic maximal operator} $M^{\mathscr{D}^t}$ are 
sublinear operators acting by the rule $f\mapsto Mf$, $f\mapsto M^\mathscr{D}f$
and $f\mapsto M^{\mathscr{D}^t}f$, respectively.
\end{definition}
Due to the available comparison between each dyadic cube and its inner/outer ball,
as well as the fine covering property of a dyadic system $\mathscr{D}$ stated in
Theorem~\ref{th:HK}, we have pointwise equivalence between the classical maximal 
function $Mf$ and the dyadic maximal function $M^\mathscr{D}f$ 
(see~\cite[Proposition~7.9]{HK12} and \cite[Theorem~3.5]{AS25}).
\begin{lemma}\label{le:M-and-dyadicM}
Let $f\in L^0(X,d,\mu)$. There holds the pointwise estimate
\[
\frac{1}{c}\,Mf(x)\le M^\mathscr{D}f(x)\le c\,Mf(x),
\]
for all $x\in X$, with a constant $c\ge 1$ independent of $f$.
\end{lemma}
We will also need the following corollary of the Lebesgue differentiation theorem for cubes
in Lemma~\ref{le:cones}(c). It shows that the $t$-dyadic maximal function of a function
dominates the function itself, giving thus one more justification for the name ``maximal.''
\begin{lemma}\label{le:maxfun}
Let $(X,d,\mu)$ be a space of homogeneous type \emph{with a Borel-semiregular measure $\mu$},
and $\mathscr{D}:=\bigcup_{t=1}^K\mathscr{D}^t$ be a dyadic system associated with the space.
For all $t\in\{1,\ldots,K\}$ and $f\in L^1_\textnormal{loc}(X,d,\mu)$, there holds
\[
|f(x)|\le M^{\mathscr{D}^t}f(x)
\quad\text{for $\mu$-almost every }x\in X.
\]
\end{lemma}
\begin{proof}
By Lemma~\ref{le:cones}(c), for $\mu$-almost every $x\in X$ we have
\[
|f(x)|=\lim_{k\to+\infty}\left|\fint_{Q^t_k(x)}f(y)\,d\mu(y)\right|
\le\sup_{\substack{Q\ni x,\\Q\in\mathscr{D}^t}}\fint_Q|f(y)|\,d\mu(y)
=M^{\mathscr{D}^t}f(x),
\]
where $Q^t_k(x)$ denotes, as in Lemma~\ref{le:cones}, the unique cube from 
the $k$-th generation of cubes in $\mathscr{D}^t$ that contains $x$.
\end{proof}

The main result of this paper characterizes the boundedness of the maximal
operator in terms of the uniform boundedness of a bunch of \emph{averaging operators}
over families of pairwise disjoint dyadic cubes. We next introduce these auxiliary dyadic
operators. 
\begin{definition}\label{def:averaging-operators}
Denote by $\mathfrak{S}$ the \emph{set of all families $\mathcal{Q}$ of pairwise disjoint
cubes} from a dyadic system $\mathscr{D}$ associated with $(X,d,\mu)$.
For $s\in[1,\infty)$, the \emph{$s$-averaging operator} 
$T_{s,\mathcal{Q}}:L^s_\text{loc}(X,d,\mu)\to L^0(X,d,\mu)$ 
over a family $\mathcal{Q}\in\mathfrak{S}$ is defined by 
\[
T_{s,\mathcal{Q}}f(x):=\sum_{Q\in\mathcal{Q}}\chi_Q(x)\,M_{s,Q}f.
\]
For simplicity, we write $T_\mathcal{Q}$ instead of $T_{1,\mathcal{Q}}$ and call it 
the {\it averaging operator} over $\mathcal{Q}$.
\end{definition}
To understand why locally-$L^s$ functions form the domain for $s$-averaging operators,
observe that a function 
$f\in L^s_\textnormal{loc}(X,d,\mu)$ has finite means
$M_{s,Q}f<\infty$ over all cubes $Q\in\mathscr{D}$: since by construction in
Theorem~\ref{th:HK} every $Q\in\mathscr{D}$ has the outer ball $B_Q\supseteq Q$, 
we easily estimate that
\[
M_{s,Q}f\le\left(\frac1{\mu(Q)}\int_{B_Q}|f(y)|^s\,d\mu(y)\right)^{1/s}<\infty.
\]

Two properties of averaging operators are worth mentioning. First, for any 
$\mathcal{Q\in\mathfrak{S}}$, the operator $T_\mathcal{Q}$ is self-dual with respect 
to nonnegative functions, i.e.,
\begin{equation}\label{eq:self-dual}
\langle T_\mathcal{Q}f,|g|\rangle
=\sum_{Q\in\mathcal{Q}}\mu(Q)\,M_Qf\,M_Qg
=\langle|f|,T_\mathcal{Q}g\rangle
\end{equation}
for all $f,g\in L^1_{\text{loc}}(X,d,\mu)$. Second, the dyadic maximal function dominates
the image
of the averaging operator over a family of pairwise disjoint dyadic cubes: for any 
$f\in L^1_\textnormal{loc}(X,d,\mu)$ and $\mathcal{Q}\in\mathfrak{S}$, one has
\begin{equation}\label{eq:aver-max}
T_\mathcal{Q}f(x)\le M^\mathscr{D}f(x),\quad x\in X.
\end{equation}
To check this inequality, fix a family $\mathcal{Q}\in\mathfrak{S}$ and 
an $x\in X$. If $x$ is outside every cube from $\mathcal{Q}$, then 
$0=T_\mathcal{Q}f(x)\le M^\mathscr{D}f(x)$. Otherwise $x\in Q'\in\mathcal{Q}$, which 
implies that $T_\mathcal{Q}f(x)=M_{Q'}f\le
\sup_{Q\ni x,\,Q\in\mathscr{D}}M_Qf=M^\mathscr{D}f(x)$.

\subsection{The Calder\'on--Zygmund lemma}
In an unbounded space of homogeneous type with a dyadic system 
$\mathscr{D}:=\bigcup_{t=1}^K\mathscr{D}^t$, the superlevel set of each $t$-dyadic maximal 
function $M^{\mathscr{D}^t}f$ admits a Calder\'on--Zygmund-type
decomposition into a family of pairwise disjoint dyadic cubes from $\mathscr{D}^t$.
Versions of the Calder\'on--Zygmund lemma in spaces of homogeneous type are definitely
known: see~\cite[Theorem~2.8]{ACUM15}, \cite[Theorem~9]{AK19}, and
\cite[Lemma~4.5]{CUC22}. What distinguishes our Lemma~\ref{le:CZ} is the maximality 
property~(iii). The proof of this property is tied to the principle by which one chooses 
the Calder\'on--Zygmund cubes, and for this reason, we present the proof of the entire lemma.
\begin{lemma}\label{le:CZ}
Suppose that $(X,d,\mu)$ is an \emph{unbounded} space of homogeneous type and 
$\mathscr{D}:=\bigcup_{t=1}^K\mathscr{D}^t$ is a dyadic system associated with it.
Let $f\in L^1(X,d,\mu)$ and $\lambda>0$. Then for any $t\in\{1,\ldots,K\}$, 
there exists a family of dyadic cubes $\{Q_j\}\subseteq\mathscr{D}^t$ such that
\[
E^t_\lambda:=\left\{x\in X\,:\,M^{\mathscr{D}^t}f(x)>\lambda\right\}=\bigcup_j Q_j,
\]
and the following properties hold:
\begin{enumerate}
    \item[(i)] \emph{(disjointness)} the cubes in $\{Q_j\}$ are pairwise disjoint and are referred to as 
    \emph{the Calder\'on--Zygmund cubes} of the superlevel set $E^t_\lambda$;
    \item[(ii)] \emph{(inequality for means)} for all $j$, one has
    \[
    \lambda<M_{Q_j}f\le\frac{\lambda}{\varepsilon},
    \]
    where $\varepsilon\in(0,1)$ is the parameter of the system $\mathscr{D}$
    specified in Theorem~\ref{th:HK}(d);
    \item[(iii)] \emph{(maximality)} none of the cubes $Q'\in\mathscr{D}^t$ satisfying 
    $Q'\varsupsetneq Q_j$ for some $j$ is contained in $E^t_\lambda$.
\end{enumerate}
We let the family $\{Q_j\}$ be empty if $E^t_\lambda=\emptyset$.
\end{lemma}
\begin{remark}\label{re:CZ}
The hypothesis that $f\in L^1(X,d,\mu)$ can be replaced by the weaker requirements
that $f\in L^1_\text{loc}(X,d,\mu)$ and, for all $x\in X$ and $t\in\{1,\ldots,K\}$, 
the means $M_{Q^t_k(x)}f\to0$ as $k\to-\infty$ (the cubes $Q^t_k(x)$ are defined in 
Lemma~\ref{le:cones}).
\end{remark}
\begin{proof}[Proof of Lemma~\ref{le:CZ}]
Let, as in Lemma~\ref{le:cones}, $Q^t_k(x)\in\mathscr{D}^t_k$ denote the unique cube from
the $k$-th generation of cubes in $\mathscr{D}^t$ that contains $x\in X$. Note first that
by Lemma~\ref{le:cones}(a), the limit $\mu(Q^t_k(x))\to\infty$ as $k\to-\infty$
holds for all $x\in X$ and $t\in\{1,\ldots,K\}$ since the space $(X,d,\mu)$ is unbounded. 
Then for $f\in L^1(X,d,\mu)$ and all $x\in X$ and $t\in\{1,\ldots,K\}$, we have 
\begin{equation}\label{eq:small-means}
M_{Q^t_k(x)}f\le\frac{\|f\|_1}{\mu(Q^t_k(x))}\to0
\quad\text{as }k\to-\infty;
\end{equation}
alternatively, we can take~\eqref{eq:small-means} as a hypothesis, see Remark~\ref{re:CZ}.

Assume that $E^t_\lambda\ne\emptyset$. In view of~\eqref{eq:small-means}, 
for each $x\in E^t_\lambda$ 
we can find the unique cube $Q^t_k(x)\in\mathscr{D}^t_k$ such that 
$M_{Q^t_k(x)}f>\lambda$ but $M_{Q^t_l(x)}f\le\lambda$ for all $l<k$. This choice of 
$k\in\mathbb{Z}$ implies by Theorem~\ref{th:HK}(d) that
\begin{align*}
M_{Q^t_k(x)}f&=\frac{1}{\mu(Q^t_k(x))}\int_{Q^t_k(x)}|f(y)|\,d\mu(y)
\\
&\le\frac{1}{\varepsilon\,\mu(Q^t_{k-1}(x))}\int_{Q^t_{k-1}(x)}|f(y)|\,d\mu(y)
\\
&=\frac{1}{\varepsilon}\,M_{Q^t_{k-1}(x)}f\le\frac{\lambda}{\varepsilon}.
\end{align*}
Among all cubes $Q\in\mathscr{D}^t$ which contain the point $x$ and have the property
$M_Qf>\lambda$, the cube $Q^t_k(x)$ is the largest. Let us define 
$\{Q_j\}\subseteq\mathscr{D}^t$ as the family of such largest cubes for each 
point $x\in E^t_\lambda$, with the convention that duplicates are omitted and 
each cube appears in this family only once. By our definition, 
$E^t_\lambda\subseteq\bigcup_j Q_j$ and property~(ii) holds straightaway.
On the other hand, if $x\in\bigcup_jQ_j$, then 
$x\in Q_j\in\mathscr{D}^t$ for some $j$ and there holds $M_{Q_j}f>\lambda$, which implies 
that $M^{\mathscr{D}^t}f(x)>\lambda$ and so $x\in E^t_\lambda$. Overall, 
$E^t_\lambda=\bigcup_jQ_j$.

To check property~(i), the pairwise disjointness of cubes in 
$\{Q_j\}\subseteq\mathscr{D}^t$, suppose to the
contrary that $Q_i\cap Q_j\ne\emptyset$ for some distinct indices $i$ and $j$
in the family. Since all cubes in $\{Q_j\}$ are different, we may assume without loss of 
generality that $Q_i\varsupsetneq Q_j$. By construction, for any $x\in Q_j$, the cube $Q_j$ 
is the largest cube in $\mathscr{D}^t$ which contains $x$ and satisfies $M_{Q_j}f>\lambda$,
therefore $M_{Q_i}f\le\lambda$. This is a contradiction with $Q_i\in\{Q_j\}$.

Finally, we prove the maximality property~(iii) by contradiction. Assume that a cube 
$Q'\in\mathscr{D}^t$, such that $Q'\varsupsetneq Q_j$ for some Calder\'on--Zygmund cube
$Q_j$ of $E^t_\lambda$, is a subset of $E^t_\lambda$. Then $Q'$ intersects more than one 
Calder\'on--Zygmund cube of $E^t_\lambda$. Namely, $Q'\subseteq\bigcup_{i\in I}Q_i$ 
for some index set $I$ with cardinality $|I|\ge2$, and
$Q'\cap Q_i\ne\emptyset$ for all $i\in I$. The situation $Q'\subseteq Q_i$ is not possible
for any $i\in I$, since the intersection of $Q'$ with the rest of the cubes would be empty;
hence $\bigcup_{i\in I}Q_i\subseteq Q'$ and eventually $Q'=\bigcup_{i\in I}Q_i$. The
inequalities $M_{Q_i}f>\lambda$, $i\in I$, are equivalent to
\[
\int_{Q_i}|f(y)|\,d\mu(y)>\lambda\,\mu(Q_i),\quad i\in I. 
\]
Summing over $i\in I$ gives $M_{Q'}f>\lambda$. This contradict the fact that for every 
$x\in Q_j$, the cube $Q_j\varsubsetneq Q'$ is the largest among the cubes 
$Q\in\mathscr{D}^t$ which contain $x$ and satisfy $M_Qf>\lambda$.
\end{proof}

\section{Generalized $\Phi$-Functions on Cubes}\label{sec:PHI}
In what follows, we work with a fixed \emph{unbounded} space of homogeneous type $(X,d,\mu)$,
where $\mu$ is a \emph{Borel-semiregular measure}. We denote it simply by $X$, writing in
particular $L^0(X):=L^0(X,d,\mu)$, $L^1(X):=L^1(X,d,\mu)$, and 
$L^1_\textnormal{loc}(X):=L^1_\textnormal{loc}(X,d,\mu)$. As before, there is a dyadic
system $\mathscr{D}:=\bigcup_{t=1}^K\mathscr{D}^t$ provided by Theorem~\ref{th:HK}, 
associated with $X$. Accordingly, the set $\mathfrak{S}$ comprises all families 
$\mathcal{Q}$ of pairwise disjoint cubes from $\mathscr{D}$. 

This section presents two important examples of generalized $\Phi$-functions on 
dyadic cubes, which are the means $M_{s,Q}\varphi$ and mirror means $(M_{s,Q}\varphi^*)^*$
of a generalized $\Phi$-function $\varphi$ on $X$. Then relations of domination and 
strong domination are defined for generalized $\Phi$-functions on cubes, 
following~\cite[Definitions~5.2.17 and 5.5.1]{DHHR11}.

\subsection{Three types of generalized $\Phi$-functions}

Let us introduce notation for the three types of generalized $\Phi$-functions
used hereafter. By $\Phi(X)$ we denote the set of generalized $\Phi$-functions on 
$(X,d,\mu)$. Then any function $\varphi\in\Phi(X)$ induces the
corresponding Musielak--Orlicz space $L^\varphi(X):=L^\varphi(X,d,\mu)$.

Next, regard the dyadic system $\mathscr{D}$ 
associated with $X$ as a countable collection of cubes ``as objects.'' We can introduce the
atomic measure $\mu_{\rm a}$ on $\mathscr{D}$ defined by $\mu_{\rm a}(Q):=\mu(Q)$ for
all $Q\in\mathscr{D}$, where $\mu(Q)$ is the measure of a cube $Q$ as a set in $X$. Every set
$\mathcal{U}\subset\mathscr{D}$ is $\mu_{\rm a}$-measurable, since
\[
\mu_{\rm a}(\mathcal{U})=\sum_{Q\in\mathcal{U}}\mu_{\rm a}(Q),
\]
therefore all functions on $(\mathscr{D},\mu_{\rm a})$ are measurable. Thus
the generalized $\Phi$-functions on $(\mathscr{D},\mu_{\rm a})$ can be defined 
using Definition~\ref{def:gen-Phi} with condition~(ii) excluded, because it is 
always satisfied. We denote the set of all such functions by $\Phi(\mathscr{D})$.

Finally, we will make extensive use of restrictions of generalized $\Phi$-functions 
on $\mathscr{D}$ to 
families of pairwise disjoint dyadic cubes, which are nice subfamilies of $\mathscr{D}$.
For any $\mathcal{Q}\in\mathfrak{S}$, we can speak about generalized $\Phi$-functions on 
$(\mathcal{Q},\mu_{\rm a})$ and denote their set by $\Phi(\mathcal{Q})$. 
Just like above, every function on $\mathcal{Q}$ is $\mu_{\rm a}$-measurable, hence
\[
L^0(\mathcal{Q},\mu_{\rm a})=\big\{
\{t_Q\}_{Q\in\mathcal{Q}}\ :\ t_Q\in\mathbb{C}
\big\}=\mathbb{C}^\mathcal{Q}.
\]
If $\varphi\in\Phi(\mathcal{Q})$, then the semimodular induced by $\varphi$ on 
$L^0(\mathcal{Q},\mu_{\rm a})$ is given by
\[
\rho_{\varphi(Q)}\big(\{t_Q\}_{Q\in\mathcal{Q}}\big):=
\sum_{Q\in\mathcal{Q}}\mu(Q)\,\varphi(Q,|t_Q|).
\]
Notice that we write $\rho_{\varphi(Q)}$ instead of $\rho_\varphi$ to highlight 
that $\varphi$ is a generalized $\Phi$-function on cubes. This helps us immediately 
realize that the semimodular $\rho_{\varphi(Q)}$ is defined on sequences 
$\{t_Q\}_{Q\in\mathcal{Q}}\in\mathbb{C}^\mathcal{Q}$, whereas the notation $\rho_\varphi$ 
will be reserved for the case when $\varphi\in\Phi(X)$ and the semimodular
$\rho_\varphi$ is defined on functions $f\in L^0(X)$.
Finally, we denote by  
$l^{\varphi(Q)}(\mathcal{Q}):=L^\varphi(\mathcal{Q},\mu_{\rm a})$
the Musielak--Orlicz {\it sequence} space induced by the semimodular $\rho_{\varphi(Q)}$.
The Luxemburg norm on this space is
\[
\big\|\{t_Q\}_{Q\in\mathcal{Q}}\big\|_{l^{\varphi(Q)}(\mathcal{Q})}:=
\inf\left\{\lambda>0\ :\ 
\rho_{\varphi(Q)}\big(\{t_Q/\lambda\}_{Q\in\mathcal{Q}}\big)\le1\right\}.
\]

Overall, the reader will meet the following three types of generalized $\Phi$-functions
in the text:
\vspace{4pt}

\begin{center}
\begin{tabular}{|c|c|}
\hline
$\varphi\in\Phi(X)$ &
\makecell{\vspace{-10pt}\\
$\varphi:X\times[0,\infty)\to[0,\infty]$ \\
is a generalized $\Phi$-function on $(X,d,\mu)$. \\[4pt]
Semimodular $\rho_\varphi(f)=\int_X\varphi(x,|f(x)|)\,d\mu(x)$ \\
induces the space $L^\varphi(X)$ with the norm $\|\cdot\|_\varphi$.\\
\vspace{-10pt}} \\
\hline
$\varphi\in\Phi(\mathscr{D})$ &
\makecell{\vspace{-10pt}\\
$\varphi:\mathscr{D}\times[0,\infty)\to[0,\infty]$, \\
$t\mapsto\varphi(Q,t)$ is a $\Phi$-function for all $Q\in\mathscr{D}$. \\
\vspace{-10pt}}\\
\hline
\makecell{$\varphi\in\Phi(\mathcal{Q})$, \\
$\mathcal{Q}\in\mathfrak{S}$} &
\makecell{\vspace{-10pt}\\
$\varphi:\mathcal{Q}\times[0,\infty)\to[0,\infty]$,  \\
$t\mapsto\varphi(Q,t)$ is a $\Phi$-function for all $Q\in\mathcal{Q}$. \\[4pt]
Semimodular $\rho_{\varphi(Q)}(\{t_Q\}_{Q\in\mathcal{Q}})=
\sum_{Q\in\mathcal{Q}}\mu(Q)\,\varphi(Q,|t_Q|)$ \\
induces the space $l^{\varphi(Q)}(\mathcal{Q})$ with the norm 
$\|\cdot\|_{l^{\varphi(Q)}(\mathcal{Q})}$.\\
\vspace{-10pt}} \\
\hline
\end{tabular}
\end{center}
\vspace{4pt}

\noindent
Note also that the classes of generalized $N$-functions within $\Phi(X)$ and
$\Phi(\mathscr{D})$ will be denoted by $N(X)$ and $N(\mathscr{D})$, respectively.

\subsection{Means of $\varphi\in\Phi(X)$ over cubes as 
generalized $\Phi$-functions on cubes}

A natural way to construct a generalized $\Phi$-function on the dyadic system $\mathscr{D}$
is as follows: we take a suitable function $\varphi\in\Phi(X)$ and compute the means of 
the measurable functions $x\mapsto\varphi(x,t)$ over all cubes $Q\in\mathscr{D}$,
in the sense of Definition~\ref{def:means}, for each value of $t\ge0$.
This gives us a function of two variables $Q$ and $t$, defined below.
\begin{definition}\label{def:means-phi}
For any $s\in[1,\infty)$, the function of \emph{means of a function $\varphi\in\Phi(X)$} 
over dyadic cubes from $\mathscr{D}$ is defined by
\begin{align*}
M_{s,Q}\varphi\ :\ &\mathscr{D}\times[0,\infty)\to[0,\infty],
\\
&(Q,t)\mapsto
\left(
\fint_Q\big(\varphi(x,t)\big)^sd\mu(x)
\right)^{1/s}=:(M_{s,Q}\varphi)(t).
\end{align*}
For simplicity, we write $M_Q\varphi$ instead of $M_{1,Q}\varphi$.
\end{definition}
The notation $M_{s,Q}\varphi$ is slightly abusive, since the name of the
function contains the symbol of its first variable $Q$. However, this notation 
conveniently highlights that we are dealing with the means of 
$\varphi\in\Phi(X)$ taken over \emph{cubes}. We adopt the name $M_{s,Q}\varphi$ 
for its intuitive appeal, viewing the symbol $Q$ as a general reference
to dyadic cubes rather than a specific cube. At the same time, we will often work with
subfamilies $\mathcal{Q}\in\mathfrak{S}$ of the dyadic system $\mathscr{D}$ and consider 
restrictions of $M_{s,Q}\varphi$ to the corresponding subdomains 
$\mathcal{Q}\times[0,\infty)$. We use the same notation $M_{s,Q}\varphi$ for these
restrictions, always ensuring that the domain is clear from the context.

Just as we remarked after Definition~\ref{def:means}, Jensen's inequality implies 
that the means are non-decreasing in parameter $s$. In the context of the above 
definition, we can say that for any function $\varphi\in\Phi(X)$, cube $Q\in\mathscr{D}$,
and $t\ge0$, there holds
\[
(M_{r,Q}\varphi)(t)\le(M_{s,Q}\varphi)(t)
\quad\text{if }1\le r<s.
\]

The properties of $M_{s,Q}\varphi$ as generalized $\Phi$-functions on the dyadic
system $\mathscr{D}$ or a family $\mathcal{Q}\in\mathfrak{S}$ are determined
by the properties of $\varphi$ as a generalized $\Phi$-function on $X$. The following 
Lemmas~\ref{le:mean-Phi}--\ref{le:mean-N} describe the transference of nice properties 
from a function $\varphi\in\Phi(X)$ to its means $M_{s,Q}\varphi$. Schematically, the main
implications here are these ($s\ge1$):
\begin{align*}
\varphi^s\text{ proper }&\implies M_{s,Q}\varphi\in\Phi(\mathscr{D});
\\
\varphi\in\Delta_2&\implies M_{s,Q}\varphi\in\Delta_2;
\\
\varphi\text{ proper, }\mathcal{Q}\in\mathfrak{S}&\implies 
M_Q\varphi\in\Phi(\mathcal{Q})\text{ proper};
\\
\left.
\begin{array}{rr}
\varphi\in N(X) \\
\varphi^s\in\Delta_2+\text{proper}
\end{array}
\right\}
&\implies M_{s,Q}\varphi\in N(\mathscr{D}).
\end{align*}

We begin by showing what are these ``suitable'' functions $\varphi\in\Phi(X)$ for which
$M_{s,Q}\varphi$ are generalized $\Phi$-function on $\mathscr{D}$. The statement below
simply extends~\cite[Corollary~5.2.12]{DHHR11} from the setting of $\mathbb{R}^n$ to 
spaces of homogeneous type, but we provide a more detailed proof. 
\begin{lemma}\label{le:mean-Phi}
Let $\varphi\in\Phi(X)$ and $s\in[1,\infty)$. If $\varphi^s$ is proper, then 
$M_{s,Q}\varphi\in\Phi(\mathscr{D})$.
\end{lemma}
\begin{proof}
Fix a cube $Q\in\mathscr{D}$ and check that $t\mapsto(M_{s,Q}\varphi)(t)$ is a 
$\Phi$-function. Convexity follows from the convexity of $\varphi$ in the second variable
and the Minkowski inequality, which imply for any 
$t_1,t_2\in\{t\ge0:(M_{s,Q}\varphi)(t)<\infty\}$ and $\alpha\in[0,1]$ that
\begin{align*}
(M_{s,Q}\varphi)(\alpha t_1+(1-\alpha)t_2)
&=\|\varphi(\cdot,\alpha t_1+(1-\alpha)t_2)\|_{L^s(Q,\frac{d\mu(x)}{\mu(Q)})}
\\
&\le\|\alpha\,\varphi(\cdot,t_1)+(1-\alpha)\,\varphi(\cdot,t_2)\|
_{L^s(Q,\frac{d\mu(x)}{\mu(Q)})}
\\
&\le\alpha\,\|\varphi(\cdot,t_1)\|_{L^s(Q,\frac{d\mu(x)}{\mu(Q)})}
+(1-\alpha)\,\|\varphi(\cdot,t_2)\|_{L^s(Q,\frac{d\mu(x)}{\mu(Q)})}
\\
&=\alpha\,(M_{s,Q}\varphi)(t_1)+(1-\alpha)\,(M_{s,Q}\varphi)(t_2).
\end{align*}
For the left-continuity, note that $\varphi\in\Phi(X)$ implies $\varphi^s\in\Phi(X)$, 
therefore $\rho_{\varphi^s}$ is a semimodular and hence is left-continuous. Then for any
$t_0\in(0,\infty)$, 
\[
\lim_{t\to t_0^-}(M_{s,Q}\varphi)(t)=\lim_{t\to t_0^-}
\left(\frac{\rho_{\varphi^s}(t\chi_Q)}{\mu(Q)}\right)^{1/s}
=\left(\frac{\rho_{\varphi^s}(t_0\chi_Q)}{\mu(Q)}\right)^{1/s}
=(M_{s,Q}\varphi)(t_0).
\]
Clearly, $(M_{s,Q}\varphi)(0)=0$ because $\varphi(x,0)=0$ for all $x\in X$. Since
$\varphi^s$ is proper, it follows that $\chi_Q\in L^{\varphi^s}(X)$ and thus 
$\rho_{\varphi^s}(t_1\chi_Q)<\infty$ for some $t_1\in(0,\infty)$. Then for every 
$0\le t\le t_1$ we have, by the convexity of $M_{s,Q}\varphi$, that
\[
(M_{s,Q}\varphi)(t)\le\frac{t}{t_1}\,(M_{s,Q}\varphi)(t_1)
=\frac{t}{t_1}\left(\frac{\rho_{\varphi^s}(t_1\chi_Q)}{\mu(Q)}\right)^{1/s},
\]
which implies the limit $(M_{s,Q}\varphi)(t)\to0$ as $t\to0^+$.

Finally, we prove that $(M_{s,Q}\varphi)(t)\to\infty$ as $t\to\infty$. Since 
$0\le(\varphi(\cdot,k))^s\nearrow\infty$ as $k\to\infty$, the monotone convergence 
theorem~\cite[Theorem~2.14]{GF99} gives
\[
\lim_{k\to\infty}(M_{s,Q}\varphi)(k)
=\lim_{k\to\infty}\left(
\fint_Q\big(\varphi(x,k)\big)^sd\mu(x)
\right)^{1/s}=\infty.
\]
This and the monotonicity of $t\mapsto(M_{s,Q}\varphi)(t)$ imply that for any $N>0$, 
there is $k\in\mathbb{N}$ such that for all $t\ge k$, one has
$(M_{s,Q}\varphi)(t)\ge(M_{s,Q}\varphi)(k)>N$. By definition, this means that 
$(M_{s,Q}\varphi)(t)\to\infty$ as $t\to\infty$.
\end{proof}

It is immediate from Definition~\ref{def:means-phi} that if $\varphi\in\Phi(X)$ satisfies
the $\Delta_2$-condition, then $M_{s,Q}\varphi$ satisfies the $\Delta_2$-condition with 
the same constant as $\varphi$ for any $s\ge1$. The inheritance of the property of being 
``proper,'' however, is more delicate. We reproduce here part of~\cite[Lemma~5.2.11]{DHHR11}
for our setting, supplementing its proof with more details for the reader's convenience.
\begin{lemma}\label{le:means-proper}
Let $\varphi\in\Phi(X)$ be proper and $\mathcal{Q}\in\mathfrak{S}$. Then 
$M_Q\varphi\in\Phi(\mathcal{Q})$ is proper. 
\end{lemma}
\begin{proof}
By Lemma~\ref{le:mean-Phi}, we have $M_Q\varphi\in\Phi(\mathscr{D})$, and
thus its restrictions $M_Q\varphi\in\Phi(\mathcal{Q})$ for any
$\mathcal{Q}\in\mathfrak{S}$.
Take a set $\mathcal{U}\subseteq\mathcal{Q}$ with $\mu_{\rm a}(\mathcal{U})<\infty$; 
we want to show that 
\[
\chi_\mathcal{U}
\in l^{M_Q\varphi}(\mathcal{Q})\cap
l^{(M_Q\varphi)^*}(\mathcal{Q}).
\]
Let $E:=\bigcup_{Q\in\mathcal{U}}Q\subseteq X$. Then 
$\mu(E)=\sum_{Q\in\mathcal{U}}\mu(Q)=\mu_{\rm a}(\mathcal{U})<\infty$. Since $\varphi$ is
proper, we have $\chi_E\in L^\varphi(X)$. Hence for some $t>0$, there holds
\[
\rho_{M_Q\varphi}(t\chi_\mathcal{U})=\sum_{Q\in\mathcal{U}}\mu(Q)\,(M_Q\varphi)(t)
=\rho_{\varphi}(t\chi_E)<\infty
\]
and thus $\chi_\mathcal{U}\in l^{M_Q\varphi}(\mathcal{Q})$. The conjugate
$\varphi^*$ is proper together with $\varphi$, so we conclude similarly that 
$\chi_\mathcal{U}\in l^{M_Q\varphi^*}(\mathcal{Q})$. By Young's inequality~\eqref{eq:Young},
there holds
\[
tu=M_Q(tu)\le(M_Q\varphi)(t)+(M_Q\varphi^*)(u)
\]
for all $t,u\ge0$ and $Q\in\mathscr{D}$, which implies by Definition~\ref{def:conjug} that
\[
(M_Q\varphi)^*(u)=\sup_{t\ge0}(tu-(M_Q\varphi)(t))
\le(M_Q\varphi^*)(u).
\]
The relation $(M_Q\varphi)^*(u)\le(M_Q\varphi^*)(u)$ for all $Q\in\mathcal{Q}$ and $u\ge0$
yields the inclusion 
$l^{M_Q\varphi^*}(\mathcal{Q})\subseteq l^{(M_Q\varphi)^*}(\mathcal{Q})$ due to 
Lemma~\ref{le:comparison}(a), thus $\chi_\mathcal{U}\in l^{(M_Q\varphi)^*}(\mathcal{Q})$.
\end{proof}

The remaining part of~\cite[Lemma~5.2.11]{DHHR11} claims, up to a change of notation, that if a function 
$\varphi\in N(X)$ is proper, then $M_Q\varphi\in N(\mathscr{D})$. This claim, 
however, is lacking some additional requirements on $\varphi$, because the fact 
that $\varphi$ is a proper generalized $N$-function on $X$ does not guarantee that 
$(M_Q\varphi)(t)<\infty$ for all $Q\in\mathscr{D}$ and $t\ge0$. Without knowing that
$t\mapsto(M_Q\varphi)(t)$ is finite-valued for all $Q\in\mathscr{D}$, one cannot conclude
that $M_Q\varphi\in N(\mathscr{D})$. We correct this statement and generalize it to the 
means $M_{s,Q}\varphi$ for all $s\ge1$.
\begin{lemma}\label{le:mean-N}
Let $\varphi\in N(X)$ and $s\in[1,\infty)$. If $\varphi^s$ is proper and satisfies the
$\Delta_2$-condition, then $M_{s,Q}\varphi\in N(\mathscr{D})$.
\end{lemma}
\begin{proof}
By Lemma~\ref{le:mean-Phi} we immediately have that $M_{s,Q}\varphi\in\Phi(\mathscr{D})$.
Fix a cube $Q\in\mathscr{D}$. Note that $\varphi\in\Phi(X)$ implies 
$\varphi^s\in\Phi(X)$; since $\varphi^s$ is proper, it follows that 
$t\chi_Q\in L^{\varphi^s}(X)$ for all $t\ge0$. Then the $\Delta_2$-condition for
$\varphi^s$ and the finite semimodular property from Lemma~\ref{le:delta2}(a) imply
\[
(M_{s,Q}\varphi)(t)=\left(\frac{\rho_{\varphi^s}(t\chi_Q)}{\mu(Q)}\right)^{1/s}<\infty
\]
for all $t\ge0$, hence the function $t\mapsto(M_{s,Q}\varphi)(t)$ is finite-valued.
Due to $\varphi\in N(X)$, the $\Phi$-functions $t\mapsto\varphi(x,t)$ are positive 
for every $x\in X$, therefore the $\Phi$-function $t\mapsto(M_{s,Q}\varphi)(t)$ is positive.

It remains to show the limits
\begin{equation}\label{eq:limits}
\lim_{t\to0^+}\frac{(M_{s,Q}\varphi)(t)}{t}=0
\quad\text{and}\quad
\lim_{t\to\infty}\frac{(M_{s,Q}\varphi)(t)}{t}=\infty
\end{equation}
to conclude that $M_{s,Q}\varphi\in N(\mathscr{D})$. Notice that for any $x\in X$, 
the function $t\mapsto\frac{\varphi(x,t)}{t}$ is strictly increasing for $t>0$, 
see~\cite[formula~(1.18)]{KR61}, and we know that
\[
\lim_{t\to0^+}\frac{\varphi(x,t)}{t}=0
\quad\text{and}\quad
\lim_{t\to\infty}\frac{\varphi(x,t)}{t}=\infty
\]
since $\varphi\in N(X)$. Then $0<(\frac{\varphi(\cdot,k)}{k})^s\nearrow\infty$
as $k\to\infty$, so the monotone convergence theorem~\cite[Theorem~2.14]{GF99} gives
\[
\lim_{k\to\infty}\frac{(M_{s,Q}\varphi)(k)}{k}
=\lim_{k\to\infty}\left(
\fint_Q\left(\frac{\varphi(x,k)}{k}\right)^sd\mu(x)
\right)^{1/s}=\infty.
\]
Consequently, for any $N>0$, there exists $k\in\mathbb{N}$ such that for all $t\ge k$, 
we have
\[
\frac{(M_{s,Q}\varphi)(t)}{t}\ge\frac{(M_{s,Q}\varphi)(k)}{k}>N.
\]
This proves the second limit in~\eqref{eq:limits}. For the first limit in~\eqref{eq:limits},
take a sequence $0<t_k\to0$ as $k\to\infty$.
For a number $\lambda>0$ such that all $t_k\le\lambda$, there holds
\[
\left(\frac{\varphi(\cdot,t_k)}{t_k}\right)^s
\le\left(\frac{\varphi(\cdot,\lambda)}{\lambda}\right)^s\in L^1(Q)
\]
with the integrability due to $(M_{s,Q}\varphi)(\lambda)<\infty$. Then it follows by the
dominated convergence theorem~\cite[Theorem~2.24]{GF99} that 
\[
\lim_{k\to\infty}\frac{(M_{s,Q}\varphi)(t_k)}{t_k}=
\left(\fint_Q\lim_{k\to\infty}\left(\frac{\varphi(x,t_k)}{t_k}\right)^s d\mu(x)\right)^{1/s}
=0,
\]
which concludes the proof of the lemma.
\end{proof}

\subsection{Mirror means $(M_{s,Q}\varphi^*)^*$ and their interpretation}

As shown in Lemma~\ref{le:mean-Phi}, if a function $\varphi\in\Phi(X)$ is such
that $\varphi^s$ is proper for some $s\ge1$, then the means of $\varphi$
over dyadic cubes $Q\in\mathscr{D}$ form the generalized $\Phi$-function $M_{s,Q}\varphi$ on 
$\mathscr{D}$. Assuming additionally that $(\varphi^*)^s$ is proper, we can construct another
example of a generalized $\Phi$-function on cubes as follows.

\begin{figure}[h!]
  \centering
\begin{tikzpicture}[>=Stealth, every node/.style={scale=1.1}]

  \coordinate (O) at (0,0); 

  \node at ($(O)+(-4.5,1)$) {$\Phi(X)$};
  \node (phi)  at ($(O)+(-2,1)$) {$\varphi$};
  \node (phi*) at ($(O)+(2,1)$) {$\varphi^*$};

  \node at ($(O)+(-4.5,-1)$) {$\Phi(\mathscr{D})$};
  \node (M-phi) at ($(O)+(-2,-1)$) {$M_{s,Q}\varphi$};
  \node (M-phi*) at ($(O)+(2,-1)$) {$M_{s,Q}\varphi^*$};
  \node (M-phi**) at ($(O)+(-2,-2)$) {$(M_{s,Q}\varphi^*)^*$};
  \node[rotate=90] at (-2,-1.5) {$\ne$};

  \draw[very thick] ($(O)+(0,2)$) -- ($(O)+(0,-3.1)$);
  \node at ($(O)+(0,2.3)$) {\small $*$};
  \node at ($(O)+(0,-3.6)$) {\tiny \shortstack{conjugation\\ ``mirror''}};
  \fill[gray, opacity=0.05] (0,2) rectangle (3.5,-3.1);

  \draw[dashed, thick] ($(O)+(-5.3,0)$) -- ($(O)+(3.5,0)$);
  \node at ($(O)+(4.2,0)$) {\tiny \shortstack{integral\\ means}};

  \draw[->, very thin] (phi) to[out=30, in=150] (phi*);
  \draw[->, very thin] (phi*) to (M-phi*);
  \draw[->, very thin] ($(M-phi*.south)+(-0,0)$) to[out=-90, in=-30] 
  ($(M-phi**.south)+(0.35,0)$);
  \draw[->, very thin, gray] (phi) to (M-phi);

\end{tikzpicture}
  \caption{Construction of the mirror means $(M_{s,Q}\varphi^*)^*$.}
  \label{fig:conjug}
\end{figure}
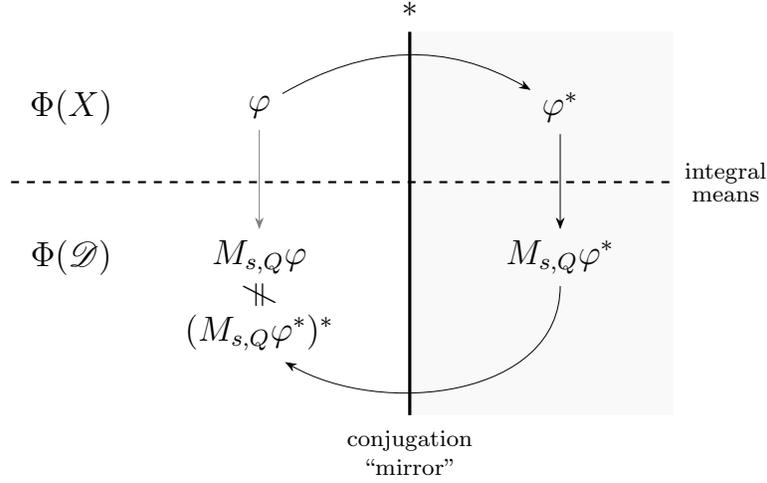

Consider $\varphi\in\Phi(X)$ and pass to its conjugate $\varphi^*$: in
Figure~\ref{fig:conjug}, this passage is symbolically depicted as a transition through 
the vertical conjugation ``mirror.'' Then form the function
$M_{s,Q}\varphi^*$ by taking the means of $\varphi^*$ over cubes $Q\in\mathscr{D}$ 
according to the rule in Definition~\ref{def:means-phi}. This function belongs to the class
$\Phi(\mathscr{D})$ when $(\varphi^*)^s$ is proper. Finally, 
step back from behind the conjugation ``mirror'' to obtain the conjugate of 
$M_{s,Q}\varphi^*$, which belongs to $\Phi(\mathscr{D})$ and will be denoted by
\[
(M_{s,Q}\varphi^*)^*\ :\ \mathscr{D}\times[0,\infty)\to[0,\infty].
\]
We refer to $(M_{s,Q}\varphi^*)^*$ as the function of {\it mirror means} of 
a function $\varphi\in\Phi(X)$ over dyadic cubes from $\mathscr{D}$.

Generally speaking, there is a gap between the values of mirror means $(M_{s,Q}\varphi^*)^*$ 
and standard means $M_{s,Q}\varphi$. As noted in~\cite[p.~152]{DHHR11}, in the translation
invariant case---when $\varphi(x,t)$ is independent of $x$---the means coincide, since
for all $Q\in\mathscr{D}$ and $t\ge0$,
\[
(M_{s,Q}\varphi^*)^*(t)=(\varphi^*)^*(t)=\varphi(t)=(M_{s,Q}\varphi)(t).
\]
This fails in general when 
$\varphi(x,t)$ is $x$-dependent, so the authors suggest interpreting the difference
between $(M_{s,Q}\varphi^*)^*$ and $M_{s,Q}\varphi$ as a measure of the failure of
translation invariance of the space $L^\varphi(X)$. The next lemma shows that the values 
of mirror means {\it do not exceed} the corresponding values of standard means
(cf. \cite[Lemma~5.2.10]{DHHR11} and \cite[Lemma~3.7]{LD05}). Hereafter, for any 
$s\in[1,\infty)$ and $f\in L^0(X)$, we denote the ordinary means of the function 
$\varphi(f):=\varphi(\cdot,|f(\cdot)|)$ over a cube $Q\in\mathscr{D}$ by 
\[
M_{s,Q}(\varphi(f)):=\left(
\fint_Q\big(\varphi(x,|f(x)|)\big)^sd\mu(x)
\right)^{1/s}.
\]
\begin{lemma}\label{le:mirror-inequality}
Let $\varphi\in\Phi(X)$ be such that $\varphi^s$ and $(\varphi^*)^s$ are proper for some
$s\ge1$. Then for all cubes $Q\in\mathscr{D}$ and functions $f\in L^0(X)$ with 
$M_{s,Q}f<\infty$, there holds
\begin{equation}\label{ineq:for-fav-lemma}
(M_{s,Q}\varphi^*)^*(M_{s,Q}f)\le M_{s,Q}(\varphi(f)).
\end{equation}
In particular, $(M_{s,Q}\varphi^*)^*(t)\le(M_{s,Q}\varphi)(t)$ for all
$Q\in\mathscr{D}$ and $t\ge0$.
\end{lemma}
\begin{proof}
Note first that $M_{s,Q}\varphi,(M_{s,Q}\varphi^*)^*\in\Phi(\mathscr{D})$ by 
Lemma~\ref{le:mean-Phi}, since $\varphi^s$ and $(\varphi^*)^s$ are proper. Fix a cube 
$Q\in\mathscr{D}$ and a function $f$ satisfying $M_{s,Q}f<\infty$. 
By Young's inequality~\eqref{eq:Young}, we have 
$u\,|f(x)|\le\varphi(x,|f(x)|)+\varphi^*(x,u)$ for all $x\in X$ and $u\ge0$.
This implies by the Minkowski inequality for the space $L^s(Q,\frac{d\mu(x)}{\mu(Q)})$ that
\[
u\,M_{s,Q}f= M_{s,Q}(uf)\le M_{s,Q}(\varphi(f))+(M_{s,Q}\varphi^*)(u).
\]
Applying Definition~\ref{def:conjug} of the conjugate function and this inequality, we obtain
\begin{align*}
(M_{s,Q}\varphi^*)^*(M_{s,Q}f)&=\sup_{u\ge0}\left(u\,M_{s,Q}f-(M_{s,Q}\varphi^*)(u)\right)
\\
&\le\sup_{u\ge0}\left(M_{s,Q}(\varphi(f))+(M_{s,Q}\varphi^*)(u)-(M_{s,Q}\varphi^*)(u)\right)
\\
&=M_{s,Q}(\varphi(f)),
\end{align*}
which is exactly~\eqref{ineq:for-fav-lemma}. Then the choice $f=t\chi_Q$ 
yields $(M_{s,Q}\varphi^*)^*(t)\le(M_{s,Q}\varphi)(t)$.
\end{proof}

A convenient feature of the mirror means $(M_{s,Q}\varphi^*)^*$ is that they inherit the
$\Delta_2$-property directly from $\varphi$. For completeness, we reproduce the proof of this
simple fact from~\cite[Lemma~5.7.10]{DHHR11}.
\begin{lemma}\label{le:delta2-means}
Let $\varphi\in\Phi(X)$ be such that $(M_{s,Q}\varphi^*)^*\in\Phi(\mathscr{D})$ for
some $s\ge1$. If $\varphi$ satisfies the $\Delta_2$-condition, then 
$(M_{s,Q}\varphi^*)^*$ satisfies the $\Delta_2$-condition with the same $\Delta_2$-constant 
as $\varphi$.
\end{lemma}
\begin{proof}
Let $D$ denote the $\Delta_2$-constant of $\varphi$. The inequality 
$\varphi(x,2t)\le D\,\varphi(x,t)$ implies that  
$\varphi^*(x,\frac{t}{2})\ge D\,\varphi^*(x,\frac{t}{D})$ for all $x\in X$ and $t\ge 0$
by Lemma~\ref{le:conjugation}. 
It follows for any cube $Q\in\mathscr{D}$ that 
$(M_{s,Q}\varphi^*)(\frac{t}{2})\ge D\,(M_{s,Q}\varphi^*)(\frac{t}{D})$. Since
$M_{s,Q}\varphi^*\in\Phi(\mathscr{D})$, we apply Lemma~\ref{le:conjugation} once again
to conclude that $(M_{s,Q}\varphi^*)^*(2t)\le D\,(M_{s,Q}\varphi^*)^*(t)$ for all
$Q\in\mathscr{D}$ and $t\ge0$.
\end{proof}

It is difficult, even in the simplest case of $(M_Q\varphi^*)^*$, to obtain an
explicit expression for the mirror means. Nevertheless, for each cube 
$Q\in\mathscr{D}$, the function $t\mapsto(M_Q\varphi^*)^*(t)$ can be represented as the
\emph{infimal convolution} of $\varphi(x,\cdot)$ with respect to the measure 
$d\mu(x)/\mu(Q)$, speaking in terms of~\cite[Section~2]{TK07}. The following lemma 
formalizes this statement and provides a variant of~\cite[Lemma~5.2.8]{DHHR11}, which itself
is a counterpart of~\cite[Proposition~2.1]{TK07}. 

Let us note that the proof of~\cite[Lemma~5.2.8]{DHHR11} contains an unexplained step, 
which we were not able to reconstruct. It is not clarified why the function $\psi$ (the 
same as in the proof of Lemma~\ref{le:mirror-mean} below) is left-continuous---a property
necessary to conclude that $\psi$ is a $\Phi$-function. However, the issue of 
left-continuity is easily resolved by adding the $\Delta_2$-condition to the assumption of
the lemma. This makes the statement less general, but still adequate for our purposes. 
Thus, Lemma~\ref{le:mirror-mean} is a variant of~\cite[Lemma~5.2.8]{DHHR11} formulated 
under the $\Delta_2$-condition and accompanied by a more detailed proof.
\begin{lemma}\label{le:mirror-mean}
Let $\varphi\in\Phi(X)$ be proper and satisfy the $\Delta_2$-condition. Then for all cubes
$Q\in\mathscr{D}$ and every $t\ge0$, one has
\[
(M_Q\varphi^*)^*(t)=\inf_{\substack{f\in L^0(X):\\M_Qf\ge t}}M_Q(\varphi(f))
=\inf_{\substack{f\in L^0(X):\\M_Qf=t}}M_Q(\varphi(f)).
\]
\end{lemma}
\begin{proof}
Note first that since $\varphi^*$ is proper together with $\varphi\in\Phi(X)$,
by Lemma~\ref{le:mean-Phi} we have $(M_Q\varphi^*)^*\in\Phi(\mathscr{D})$.
Fix a cube $Q\in\mathscr{D}$ and $t\ge0$. For any $f\in L^0(X)$ satisfying 
$t\le M_Qf<\infty$, the monotonicity of $(M_Q\varphi^*)^*$ in the second variable
and inequality~\eqref{ineq:for-fav-lemma} for $s=1$ yield
\begin{equation}\label{eq:fav-1}
(M_Q\varphi^*)^*(t)\le(M_Q\varphi^*)^*(M_Qf)\le M_Q(\varphi(f)).
\end{equation}
Alternatively, if $M_Qf=\infty$, it follows from H\"older's inequality in
Lemma~\ref{le:Holder} that
\[
\infty=\int_Q|f(x)|\,d\mu(x)=\langle|f|\chi_Q,\chi_Q\rangle\le2\,\|f\chi_Q\|_\varphi
\|\chi_Q\|_{\varphi^*}.
\]
Since $\varphi$ is proper and thus $\|\chi_Q\|_{\varphi^*}<\infty$, we conclude that 
$f\chi_Q\notin L^\varphi(X)$. Then the inequality
\begin{equation}\label{eq:fav-2}
(M_Q\varphi^*)^*(t)\le M_Q(\varphi(f))=\frac{\rho_\varphi(f\chi_Q)}{\mu(Q)}=\infty
\end{equation}
holds trivially. Combining inequality~\eqref{eq:fav-1} for measurable 
$f$ satisfying $t\le M_Qf<\infty$ with inequality~\eqref{eq:fav-2} for measurable $f$
satisfying $M_Qf=\infty$, we conclude that
\[
(M_Q\varphi^*)^*(t)\le\inf_{\substack{f\in L^0(X):\\M_Qf\ge t}}M_Q(\varphi(f))
\le\inf_{\substack{f\in L^0(X):\\M_Qf=t}}M_Q(\varphi(f))=:\psi(t).
\]
The second inequality here is obvious, as it compares the infima over a “wider” and
a “narrower” set. To complete the proof, it suffices to show that 
$(M_Q\varphi^*)^*(t)=\psi(t)$. 

We first show that $t\mapsto\psi(t)$ is a $\Phi$-function.
Since $f=t\chi_Q$ satisfies $M_Qf=t$, it follows that 
$\psi(t)\le M_Q(\varphi(f))=(M_Q\varphi)(t)$.
Together with the estimate proved above, this gives
\begin{equation}\label{eq:fav-3}
(M_Q\varphi^*)^*(t)\le\psi(t)\le(M_Q\varphi)(t)
\end{equation}
for all $t\ge0$. By Lemma~\ref{le:mean-Phi}, both 
$t\mapsto(M_Q\varphi^*)^*(t)$ and $t\mapsto(M_Q\varphi)(t)$ are $\Phi$-functions. As a 
function squeezed between two $\Phi$-functions, $\psi$ inherits their boundary 
properties that $\psi(0)=0$, $\psi(t)\to0$ as $t\to0^+$, and $\psi(t)\to\infty$ as 
$t\to\infty$. 
Next, let us check that $\psi$ is convex: fix $t_1,t_2\ge0$ and $\alpha\in[0,1]$. By 
applying the well-known equality $\inf(A+B)=\inf A +\inf B$
for the Minkowski sum $A+B=\{x+y:x\in A,\,y\in B\}$
together with the convexity of $\varphi$ in the second variable, we obtain that
\begin{align*}
\alpha\,\psi(t_1)+(1-\alpha)\,\psi(t_2)
&=\alpha\inf_{\substack{f_1\in L^0(X):\\M_Qf_1=t_1}}M_Q(\varphi(f_1))
+(1-\alpha)\inf_{\substack{f_2\in L^0(X):\\M_Qf_2=t_2}}M_Q(\varphi(f_2))
\\
&=\inf_{\substack{f_1,f_2\in L^0(X):\\M_Qf_1=t_1,\,M_Qf_2=t_2}}
M_Q(\alpha\,\varphi(f_1)+(1-\alpha)\,\varphi(f_2))
\\
&\ge\inf_{\substack{f_1,f_2\in L^0(X):\\M_Qf_1=t_1,\,M_Qf_2=t_2}}
M_Q(\varphi(\alpha |f_1|+(1-\alpha)|f_2|))
\\
&\ge\inf_{\substack{f\in L^0(X):\\M_Qf=\alpha t_1+(1-\alpha)t_2}}
M_Q(\varphi(f))
=\psi(\alpha t_1+(1-\alpha)t_2).
\end{align*}
Finally, notice that for all $t\ge0$, the functions $t\chi_Q\in L^\varphi(X)$ since 
$\varphi$ is proper; therefore, it follows from the $\Delta_2$-condition for $\varphi$
that $(M_Q\varphi)(t)=\rho_\varphi(t\chi_Q)/\mu(Q)<\infty$ by Lemma~\ref{le:delta2}(a). Then 
estimate~\eqref{eq:fav-3} yields that $\psi$ is finite-valued, so it is continuous by 
the convexity. Thus we have shown that $\psi$ is a $\Phi$-function.

Now, we can calculate the conjugate of $\psi$. By Definition~\ref{def:conjug}
for the case without the $x$-dependence and Lemma~\ref{le:int-phi*},
we obtain for any $u\ge0$ that
\begin{align*}
\psi^*(u)&=\sup_{t\ge0}\bigg(tu-\inf_{\substack{f\in L^0(X):\\M_Qf=t}}M_Q(\varphi(f))\bigg)
\\
&=\sup_{t\ge0}\,\sup_{\substack{f\in L^0(X):\\M_Qf=t}}
\big(tu-M_Q(\varphi(f))\big)
\\
&=\sup_{f\in L^1(Q)}\big(u\,M_Qf-M_Q(\varphi(f))\big)
\\
&=\frac1{\mu(Q)}\sup_{f\in L^1(Q)}\int_Q
(u|f(x)|-\varphi(x,|f(x)|))\,d\mu(x)
\\
&=\frac{1}{\mu(Q)}\int_Q\varphi^*(x,u)\,d\mu(x)=(M_Q\varphi^*)(u).
\end{align*}
Thus the equality $\psi(t)=(\psi^*)^*(t)=(M_Q\varphi^*)^*(t)$ holds for all $t\ge0$. 
This concludes the proof of the lemma.
\end{proof}

Equivalently, Lemma~\ref{le:mirror-mean} gives a pointwise characterization
of $(M_Q\varphi^*)^*$ as a function from the class $\Phi(\mathcal{Q})$ for an arbitrary
$\mathcal{Q}\in\mathfrak{S}$. Interpreted this way, it has a useful consequence
for the corresponding semimodular $\rho_{(M_Q\varphi^*)^*}$ on 
$\mathbb{C}^\mathcal{Q}$. The following result is 
from~\cite[Lemma~5.2.13]{DHHR11}, where it is stated without proof.
\begin{lemma}\label{le:mirror-modular}
Let $\varphi\in\Phi(X)$ be proper and satisfy the $\Delta_2$-condition, and let
$\mathcal{Q}\in\mathfrak{S}$. Then the semimodular $\rho_{(M_Q\varphi^*)^*}$ on 
$\mathbb{C}^\mathcal{Q}$ admits the representation 
\[
\rho_{(M_Q\varphi^*)^*}\big(\{t_Q\}_{Q\in\mathcal{Q}}\big)
=\inf_{\substack{f\in L^0(X):\\ M_Qf=t_Q,\,Q\in\mathcal{Q}}} \rho_\varphi(f)
\]
for all nonnegative sequences $\{t_Q\}_{Q\in\mathcal{Q}}$.  
\end{lemma}
\begin{proof}
Note first that Lemma~\ref{le:mean-Phi} implies $(M_Q\varphi^*)^*\in\Phi(\mathcal{Q})$ 
since $\varphi$ is proper.
It follows from Lemma~\ref{le:mirror-mean} that 
\begin{align*}
\rho_{(M_Q\varphi^*)^*}\big(\{t_Q\}_{Q\in\mathcal{Q}}\big)
&=\sum_{Q\in\mathcal{Q}}\mu(Q)\,(M_Q\varphi^*)^*(t_Q) 
=\sum_{Q\in\mathcal{Q}}\mu(Q)\inf_{\substack{f\in L^0(X):\\M_Qf=t_Q}} M_Q(\varphi(f)) \\
&=\sum_{Q\in\mathcal{Q}}\,\inf_{\substack{f\in L^0(X):\\M_Qf=t_Q}}\rho_\varphi(f\chi_Q) 
\overset{!}{=}\inf_{\substack{f\in L^0(X):\\M_Qf=t_Q,\,Q\in\mathcal{Q}}}\,
\sum_{Q\in\mathcal{Q}} \rho_\varphi(f\chi_Q) \\
&=\inf_{\substack{f\in L^0(X):\\M_Qf=t_Q,\,Q\in\mathcal{Q}}}\rho_\varphi(f).
\end{align*}

The penultimate equality (!) is not quite trivial. The part ``$\le$'' of it follows
simply from the classical fact that the sum of infima does not exceed the infimum of 
the sum. For the converse, let us enumerate the cubes of the family 
$\mathcal{Q}=\{Q_k\}_{k=1}^\infty$ and denote
\[
I_k:=\inf_{\substack{f\in L^0(X):\\M_{Q_k}f=t_{Q_k}}}\rho_\varphi(f\chi_{Q_k}),
\quad k\in\mathbb{N}.
\]
Fix an $\varepsilon>0$. For any $k\in\mathbb{N}$, there exists $f_k\in L^0(X)$ such that 
$M_{Q_k}f_k=t_{Q_k}$ and $\rho_\varphi(f_k\,\chi_{Q_k})<I_k+\frac{\varepsilon}{2^k}$. Define
\[
\tilde{f}:=\sum_{k=1}^\infty f_k\,\chi_{Q_k},
\]
note that $M_{Q_k}\tilde{f}=M_{Q_k}f_k=t_{Q_k}$ for all $k\in\mathbb{N}$. Therefore
\begin{align*}
\inf_{\substack{f\in L^0(X):\\M_{Q_k}f=t_{Q_k},\,k\in\mathbb{N}}}\,
\sum_{k=1}^\infty \rho_\varphi(f\chi_{Q_k}) 
&\le\sum_{k=1}^\infty \rho_{\varphi}(\tilde{f}\chi_{Q_k})
=\sum_{k=1}^\infty \rho_{\varphi}(f_k\,\chi_{Q_k}) \\
&\le \sum_{k=1}^\infty \left(I_k+\frac{\varepsilon}{2^k}\right)
=\sum_{k=1}^\infty\inf_{\substack{f\in L^0(X):\\M_{Q_k}f=t_{Q_k}}}\rho_\varphi(f\chi_{Q_k})
+\varepsilon.
\end{align*}
Since $\varepsilon>0$ is arbitrary, we get the desired part ``$\ge$'' of the non-trivial
equality~(!).
\end{proof}

\subsection{Notion of domination for generalized $\Phi$-functions on cubes}

It is common to say that the embedding of a normed space $(U,\|\cdot\|_U)$ into 
a normed space $(V,\|\cdot\|_V)$ is \emph{bounded}, and write $U\hookrightarrow V$, if 
$U\subseteq V$ and there is a constant $c>0$ such that 
\[\|u\|_V\le c\,\|u\|_U
\quad\text{for all }u\in U.
\] 
Bounded embeddings provide a convenient language for defining a domination relation 
between generalized $\Phi$-functions on dyadic cubes. We adapt this language
from~\cite{DHHR11} and present below a dyadic analogue of~\cite[Definition~5.2.17]{DHHR11}.
\begin{definition}\label{def:domination}
Let $\varphi,\psi\in\Phi(\mathscr{D})$. We say that $\psi$ is {\it dominated} by $\varphi$,
and write $\psi\preceq\varphi$, if the embeddings
$l^{\varphi(Q)}(\mathcal{Q})\hookrightarrow l^{\psi(Q)}(\mathcal{Q})$
are bounded uniformly for all families $\mathcal{Q}\in\mathfrak{S}$; that is, 
if there exists a constant $c>0$ such that 
\begin{equation}\label{ineq:domination}
\big\|\{t_Q\}_{Q\in\mathcal{Q}}\big\|_{l^{\psi(Q)}(\mathcal{Q})}
\le c\,\big\|\{t_Q\}_{Q\in\mathcal{Q}}\big\|_{l^{\varphi(Q)}(\mathcal{Q})}
\end{equation}
holds for all $\mathcal{Q}\in\mathfrak{S}$ and 
$\{t_Q\}_{Q\in\mathcal{Q}}\in l^{\varphi(Q)}(\mathcal{Q})$.
We write 
$\psi\cong\varphi$ if $\psi\preceq\varphi$ and $\varphi\preceq\psi$.
\end{definition}
This definition yields two observations.
First, the domination relation is transitive: 
if $\psi\preceq\kappa$ and $\kappa\preceq\varphi$, then $\psi\preceq\varphi$. Second,
if $\psi(Q,t)\le\varphi(Q,t)$ for all $Q\in\mathscr{D}$ and $t\ge0$, then 
$\psi\preceq\varphi$ by Lemma~\ref{le:comparison}(a). For example, we know that for any
$\varphi\in\Phi(X)$ and $s\in[1,\infty)$, Jensen's inequality implies 
$(M_Q\varphi)(t)\le(M_{s,Q}\varphi)(t)$ for all $Q\in\mathscr{D}$ and $t\ge0$. Hence 
for any $\varphi\in\Phi(X)$ such that $M_Q\varphi,M_{s,Q}\varphi\in\Phi(\mathscr{D})$, we
have
\[
M_Q\varphi\preceq M_{s,Q}\varphi
\quad\text{if }s\ge1.
\]

Under a regularity assumption on the \emph{dominant} function, the next lemma establishes 
the ``reverse'' domination between conjugate functions (cf.~\cite[Remark~5.2.19]{DHHR11}). 
\begin{lemma}\label{le:domin-conj}
Let $\varphi,\psi\in\Phi(\mathscr{D})$. Suppose that $\varphi$ is proper as an element 
of $\Phi(\mathcal{Q})$ for every $\mathcal{Q}\in\mathfrak{S}$. If 
$\psi\preceq\varphi$, then $\varphi^*\preceq\psi^*$.
\end{lemma}
\begin{proof}
The function $\varphi^*\in\Phi(\mathcal{Q})$ is proper for each $\mathcal{Q}\in\mathfrak{S}$
since $\varphi\in\Phi(\mathcal{Q})$ has the same property. Therefore, we can apply 
Lemma~\ref{le:norm-conj} to the function 
$\varphi^*\in\Phi(\mathcal{Q})$, observing that $(\varphi^*)^*=\varphi$. This, norm
inequality~\eqref{ineq:domination}, and H\"older's
inequality in Lemma~\ref{le:Holder} together yield
\begin{align*}
\big\|\{t_Q\}_{Q\in\mathcal{Q}}\big\|_{l^{\varphi^*(Q)}(\mathcal{Q})}
&\le\sup_{\|\{s_Q\}_{Q\in\mathcal{Q}}\|_{l^{\varphi(Q)}(\mathcal{Q})}\le1}
\,\sum_{Q\in\mathcal{Q}}\mu(Q)\,|t_Qs_Q|
\\
&\le\sup_{\|\{s_Q\}_{Q\in\mathcal{Q}}\|_{l^{\psi(Q)}(\mathcal{Q})}\le c}
\,\sum_{Q\in\mathcal{Q}}\mu(Q)\,|t_Qs_Q|
\\
&=c\,\sup_{\|\{s_Q\}_{Q\in\mathcal{Q}}\|_{l^{\psi(Q)}(\mathcal{Q})}\le 1}
\,\sum_{Q\in\mathcal{Q}}\mu(Q)\,|t_Qs_Q|
\\
&\le 2c\,\|\{t_Q\}_{Q\in\mathcal{Q}}\|_{l^{\psi^*(Q)}(\mathcal{Q})}
\end{align*}
for all $\mathcal{Q}\in\mathfrak{S}$ and 
$\{t_Q\}_{Q\in\mathcal{Q}}\in\mathbb{C}^\mathcal{Q}$.
Thus the embeddings
$l^{\psi^*(Q)}(\mathcal{Q})\hookrightarrow l^{\varphi^*(Q)}(\mathcal{Q})$ 
are bounded uniformly for all $\mathcal{Q}\in\mathfrak{S}$, so $\varphi^*\preceq\psi^*$.
\end{proof}

Imposing certain assumption on the \emph{dominated} function gives some benefits as well. 
In particular, Definition~\ref{def:domination} can be reformulated in terms of the
semimodulars rather than the norms if the dominated function satisfies the 
$\Delta_2$-condition. Lemma~\ref{le:domination}, which formalizes this, is inspired by the
early definition of domination in~\cite[Definition~3.5]{LD05}.
\begin{lemma}[Equivalent definition]\label{le:domination}
Let $\varphi,\psi\in\Phi(\mathscr{D})$ and $\psi$ satisfy the $\Delta_2$-condition. Then 
$\psi\preceq\varphi$ if and only if there exists a constant $A_1>1$ such that for all $\mathcal{Q}\in\mathfrak{S}$ and
$\{t_Q\}_{Q\in\mathcal{Q}}\in l^{\varphi(Q)}(\mathcal{Q})$ satisfying
\[
\rho_{\varphi(Q)}\big(\{t_Q\}_{Q\in\mathcal{Q}}\big)
=\sum_{Q\in\mathcal{Q}}\mu(Q)\,\varphi(Q,|t_Q|)\le1,
\]
there holds
\[
\rho_{\psi(Q)}\big(\{t_Q\}_{Q\in\mathcal{Q}}\big)
=\sum_{Q\in\mathcal{Q}}\mu(Q)\,\psi(Q,|t_Q|)\le A_1.
\]
\end{lemma}
\begin{proof}
For the necessity part, take a family $\mathcal{Q}\in\mathfrak{S}$ and
a sequence $\{t_Q\}_{Q\in\mathcal{Q}}\in l^{\varphi(Q)}(\mathcal{Q})$ 
with $\rho_{\varphi(Q)}(\{t_Q\}_{Q\in\mathcal{Q}})\le1$. Then
$\|\{t_Q\}_{Q\in\mathcal{Q}}\|_{l^{\varphi(Q)}(\mathcal{Q})}\le1$ by the unit ball 
property in Lemma~\ref{le:spaces}(e). 
Since $\psi\preceq\varphi$, 
it follows that $\|\{t_Q\}_{Q\in\mathcal{Q}}\|_{l^{\psi(Q)}(\mathcal{Q})}\le c$
with some $c>0$ independent of $\mathcal{Q}$ and $\{t_Q\}_{Q\in\mathcal{Q}}$.
Choose a number $k\in\mathbb{N}$ such that $c\le2^k$, and therefore
$\|2^{-k}\,\{t_Q\}_{Q\in\mathcal{Q}}\|_{l^{\psi(Q)}(\mathcal{Q})}\le 1$. By the unit ball 
property and the $\Delta_2$-condition for $\psi$ with constant $D\ge2$, we have 
\[
\rho_{\psi(Q)}\big(\{t_Q\}_{Q\in\mathcal{Q}}\big)
\le D^k\,\rho_{\psi(Q)}\big(2^{-k}\,\{t_Q\}_{Q\in\mathcal{Q}}\big)\le D^k=:A_1.
\]

For the sufficiency part, take a family $\mathcal{Q}\in\mathfrak{S}$ and a sequence 
$\{t_Q\}_{Q\in\mathcal{Q}}\in l^{\varphi(Q)}(\mathcal{Q})$ satisfying 
$\|\{t_Q\}_{Q\in\mathcal{Q}}\|_{l^{\varphi(Q)}(\mathcal{Q})}\le1$, or equivalently,
$\rho_{\varphi(Q)}(\{t_Q\}_{Q\in\mathcal{Q}})\le1$. It follows by assumption that
$\rho_{\psi(Q)}(\{t_Q\}_{Q\in\mathcal{Q}})\le A_1$ with a constant $A_1>1$ independent
of $\mathcal{Q}$ and $\{t_Q\}_{Q\in\mathcal{Q}}$. The convexity of $\rho_{\psi(Q)}$ implies
\[
\rho_{\psi(Q)}\left(\frac{1}{A_1}\,\{t_Q\}_{Q\in\mathcal{Q}}\right)
\le\frac1{A_1}\,\rho_{\psi(Q)}\big(\{t_Q\}_{Q\in\mathcal{Q}}\big)\le1.
\]
Using the unit ball property, we arrive at
$\|\{t_Q\}_{Q\in\mathcal{Q}}\|_{l^{\psi(Q)}(\mathcal{Q})}\le A_1$.
A scaling argument gives
$\|\{t_Q\}_{Q\in\mathcal{Q}}\|_{l^{\psi(Q)}(\mathcal{Q})}\le A_1\,
\|\{t_Q\}_{Q\in\mathcal{Q}}\|_{l^{\varphi(Q)}(\mathcal{Q})}$
for all $\mathcal{Q}\in\mathfrak{S}$ and 
$\{t_Q\}_{Q\in\mathcal{Q}}\in l^{\varphi(Q)}(\mathcal{Q})$, hence $\psi\preceq\varphi$
by Definition~\ref{def:domination}.
\end{proof}

Observe that if $\psi\cong\varphi$, then for any $\mathcal{Q}\in\mathfrak{S}$, the spaces 
$l^{\varphi(Q)}(\mathcal{Q})$ and $l^{\psi(Q)}(\mathcal{Q})$ coincide up to the 
equivalence of norms, and there exists a constant $c\ge1$ such that 
\begin{equation}\label{ineq:norm-equiv}
\frac1c\,\big\|\{t_Q\}_{Q\in\mathcal{Q}}\big\|_{l^{\varphi(Q)}(\mathcal{Q})}\le
\big\|\{t_Q\}_{Q\in\mathcal{Q}}\big\|_{l^{\psi(Q)}(\mathcal{Q})}
\le c\,\big\|\{t_Q\}_{Q\in\mathcal{Q}}\big\|_{l^{\varphi(Q)}(\mathcal{Q})}
\end{equation}
holds for all $\mathcal{Q}\in\mathfrak{S}$ and 
$\{t_Q\}_{Q\in\mathcal{Q}}\in l^{\varphi(Q)}(\mathcal{Q})=l^{\psi(Q)}(\mathcal{Q})$.
When both $\varphi$ and $\psi$ satisfy the $\Delta_2$-condition, their equivalence in
terms of domination has a useful consequence, stated in the following lemma. This 
is a simplified version of~\cite[Lemma~5.7.11]{DHHR11}, sufficient for our needs.
\begin{lemma}\label{le:mod-equiv-1}
Let $\varphi,\psi\in\Phi(\mathscr{D})$ satisfy the $\Delta_2$-condition, and let  
$\psi\cong\varphi$. Then 
\[
\rho_{\varphi(Q)}\big(\{t_Q\}_{Q\in\mathcal{Q}}\big)=1 \implies
\rho_{\psi(Q)}\big(\{t_Q\}_{Q\in\mathcal{Q}}\big)\approx1
\]
uniformly in $\mathcal{Q}\in\mathfrak{S}$ and 
$\{t_Q\}_{Q\in\mathcal{Q}}\in l^{\varphi(Q)}(\mathcal{Q})=l^{\psi(Q)}(\mathcal{Q})$.
\end{lemma}
\begin{proof}
Consider $\mathcal{Q}\in\mathfrak{S}$, and let
$\{t_Q\}_{Q\in\mathcal{Q}}\in l^{\varphi(Q)}(\mathcal{Q})$ satisfy 
$\rho_{\varphi(Q)}(\{t_Q\}_{Q\in\mathcal{Q}})=1$, which is equivalent to 
$\|\{t_Q\}_{Q\in\mathcal{Q}}\|_{l^{\varphi(Q)}(\mathcal{Q})}=1$ by the unit sphere
property in Lemma~\ref{le:delta2}(c) due to the $\Delta_2$-condition for $\varphi$. 
Since $\psi\cong\varphi$, norm
inequality~\eqref{ineq:norm-equiv} holds with a constant $c\ge1$ independent of 
$\mathcal{Q}$ and $\{t_Q\}_{Q\in\mathcal{Q}}$. Choose a number $k\in\mathbb{N}$ such that 
$c<2^k$. On the one hand, inequality~\eqref{ineq:norm-equiv} and the choice of $k$
imply $\|\{t_Q\}_{Q\in\mathcal{Q}}\|_{l^{\psi(Q)}(\mathcal{Q})}<2^k$, so we have
by the $\Delta_2$-condition for $\psi$ with constant $D\ge2$ and the unit ball property that
\[
\rho_{\psi(Q)}\big(\{t_Q\}_{Q\in\mathcal{Q}}\big)\le
D^k\,\rho_{\psi(Q)}\big(2^{-k}\,\{t_Q\}_{Q\in\mathcal{Q}}\big)\le D^k.
\]
On the other hand, it follows from~\eqref{ineq:norm-equiv} that 
$\|\{t_Q\}_{Q\in\mathcal{Q}}\|_{l^{\psi(Q)}(\mathcal{Q})}>2^{-k}$, which implies
\[
\rho_{\psi(Q)}\big(\{t_Q\}_{Q\in\mathcal{Q}}\big)\ge
D^{-k}\,\rho_{\psi(Q)}\big(2^k\,\{t_Q\}_{Q\in\mathcal{Q}}\big)>D^{-k}.
\]
Since the positive constants $D^k$ and $D^{-k}$ do not depend on the choice of a family 
$\mathcal{Q}\in\mathfrak{S}$ and a sequence
$\{t_Q\}_{Q\in\mathcal{Q}}\in l^{\varphi(Q)}(\mathcal{Q})=l^{\psi(Q)}(\mathcal{Q})$,
the claim of the lemma follows from the above inequalities.
\end{proof}

Finally, a pointwise characterization of domination is available when the 
dominated function satisfies the $\Delta_2$-condition. The next theorem, 
up to a change of notation, reproduces~\cite[Theorem~5.6.10]{DHHR11}. 
However, we remove some unnecessary assumptions from the original statement and 
provide a slightly more detailed exposition of the proof for the reader's convenience.
\begin{theorem}[Pointwise characterization of domination]\label{th:pointwise-domin}
Let $\varphi,\psi\in\Phi(\mathscr{D})$ and $\psi$ satisfy the $\Delta_2$-condition.
Then $\psi\preceq\varphi$ if and only if there exist $A_2>1$ and 
a function $b:\mathscr{D}\to[0,\infty)$ with
\[
\|b\|_{\mathfrak{S},1}:=\sup_{\mathcal{Q}\in\mathfrak{S}}
\sum_{Q\in\mathcal{Q}}\mu(Q)\,b(Q)\le A_2
\]
such that, for all $Q\in\mathscr{D}$ and $t\ge0$,
\begin{equation}\label{eq:point-0}
\mu(Q)\,\varphi(Q,t)\le1
\implies
\psi(Q,t)\le A_2\,\varphi(Q,t)+b(Q).
\end{equation}
\end{theorem}
\begin{proof}
Suppose first that $\psi\preceq\varphi$. Fix a nonnegative sequence 
$\{t_Q\}_{Q\in\mathscr{D}}$ on the set of dyadic cubes 
$\mathscr{D}$. Let a family $\mathcal{Q}\in\mathfrak{S}$ be such that 
\begin{equation}\label{eq:point-1}
\sum_{Q\in\mathcal{Q}}\mu(Q)\,\varphi(Q,t_Q)
=\rho_{\varphi(Q)}\big(\{t_Q\}_{Q\in\mathcal{Q}}\big)\le4.
\end{equation}
Then $\rho_{\varphi(Q)}(\frac14\{t_Q\}_{Q\in\mathcal{Q}})\le1$ by the convexity of
$\rho_{\varphi(Q)}$. Since $\psi$ satisfies the 
$\Delta_2$-condition, it follows from Lemma~\ref{le:domination} that 
$\rho_{\psi(Q)}(\frac14\{t_Q\}_{Q\in\mathcal{Q}})\le A_1$, where $A_1>1$ does not 
depend on $\mathcal{Q}$. 
If $D\ge2$ is the $\Delta_2$-constant of $\psi$, we get
\begin{align}\label{eq:point-2}
\sum_{Q\in\mathcal{Q}}\mu(Q)\,\psi(Q,t_Q)&=
\rho_{\psi(Q)}\big(\{t_Q\}_{Q\in\mathcal{Q}}\big)
\nonumber \\
&\le D^2\,\rho_{\psi(Q)}\left(\frac14\,\{t_Q\}_{Q\in\mathcal{Q}}\right)
\le D^2 A_1=:A_2.
\end{align}
Implication ``\eqref{eq:point-1}$\implies$\eqref{eq:point-2}'' allows us to 
apply~\cite[Lemma~5.6.1]{DHHR11} with $X:=\mathscr{D}$ and 
$Y:=\mathfrak{S}$, and conclude that for every $\{t_Q\}_{Q\in\mathscr{D}}$, 
there exists $a_{\{t_Q\}}:\mathscr{D}\to[0,\infty)$ satisfying
\begin{equation}\label{eq:point-3}
\sup_{\mathcal{Q}\in\mathfrak{S}}\sum_{Q\in\mathcal{Q}}a_{\{t_Q\}}(Q)\le A_2
\end{equation}
such that for all $Q\in\mathscr{D}$ there holds
\[
\mu(Q)\,\varphi(Q,t_Q)\le1
\implies
\mu(Q)\,\psi(Q,t_Q)\le A_2\,\mu(Q)\,\varphi(Q,t_Q)+a_{\{t_Q\}}(Q).
\]

Define $b:\mathscr{D}\times[0,\infty)\to[0,\infty)$ by
\[
b(Q,t):=\left\{
\begin{array}{ll}
     \frac1{\mu(Q)}\inf\limits_{\substack{\{t_Q\}_{Q\in\mathscr{D}}\\ \text{with }t_Q=t}}
     a_{\{t_Q\}}(Q) & \text{if}\quad\mu(Q)\,\varphi(Q,t_Q)\le1,  \\
     0 & \text{otherwise}. 
\end{array}
\right.
\]
Then for all $Q\in\mathscr{D}$ and $t\ge0$,
\begin{equation}\label{eq:point-4}
\mu(Q)\,\varphi(Q,t)\le1
\implies
\psi(Q,t)\le A_2\,\varphi(Q,t)+b(Q,t);
\end{equation}
moreover, for all nonnegative sequences $\{t_Q\}_{Q\in\mathscr{D}}$ and all families 
$\mathcal{Q}\in\mathfrak{S}$, there holds
\begin{equation}\label{eq:point-5}
\sum_{Q\in\mathcal{Q}}\mu(Q)\,b(Q,t_Q)\le\sum_{Q\in\mathcal{Q}}a_{\{t_Q\}}(Q)\le A_2
\end{equation}
in view of~\eqref{eq:point-3}. Finally, define $b:\mathscr{D}\to[0,\infty)$ by 
\[
b(Q):=\sup_{t\ge0}b(Q,t).
\]
Note that this function is indeed finite-valued, since choosing $\mathcal{Q}:=\{Q\}$ and
$t_Q:=t$ in estimate~\eqref{eq:point-5} implies $b(Q,t)\le A_2/\mu(Q)$ for all $t\ge0$. 
Then our claim~\eqref{eq:point-0} is an immediate consequence 
of~\eqref{eq:point-4}; it only remains to check that $\|b\|_{\mathfrak{S},1}\le A_2$.

Take an arbitrary family $\mathcal{Q}\in\mathfrak{S}$ and enumerate its cubes as
$\mathcal{Q}=\{Q_k\}_{k=1}^\infty$. Fix an $\varepsilon>0$. By the definition of $b$, 
for every $k\in\mathbb{N}$ there exists $t_{Q_k}\ge0$ such that 
\[
b(Q_k,t_{Q_k})>b(Q_k)-\frac{\varepsilon}{2^k\,\mu(Q_k)}.
\]
This and~\eqref{eq:point-5} together give 
\[
\sum_{Q\in\mathcal{Q}}\mu(Q)\,b(Q)
\le\sum_{k=1}^\infty\mu(Q_k)\,b(Q_k,t_{Q_k})+\sum_{k=1}^\infty\frac{\varepsilon}{2^k}
\le A_2+\varepsilon,
\]
where $A_2$ does not depend on $\mathcal{Q}$.
By letting $\varepsilon\to0^+$ and taking the supremum over all $\mathcal{Q}\in\mathfrak{S}$,
we arrive at $\|b\|_{\mathfrak{S},1}\le A_2$.

Conversely, if there exist $b:\mathscr{D}\to[0,\infty)$ and $A_2>0$ such that 
$\|b\|_{\mathfrak{S},1}\le A_2$ and~\eqref{eq:point-0} hold, then
for all $\mathcal{Q}\in\mathfrak{S}$ and 
$\{t_Q\}_{Q\in\mathcal{Q}}\in l^{\varphi(Q)}(\mathcal{Q})$, 
\[
\sum_{Q\in\mathcal{Q}}\mu(Q)\,\varphi(Q,|t_Q|)\le1
\implies
\sum_{Q\in\mathcal{Q}}\mu(Q)\,\psi(Q,|t_Q|)\le2A_2.
\]
The $\Delta_2$-condition for $\psi$ and Lemma~\ref{le:domination} imply 
$\psi\preceq\varphi$.
\end{proof}

\subsection{Strong domination between generalized $\Phi$-functions on cubes}

Following the approach initiated by Diening in~\cite{LD05}, we complement the notion 
of domination for generalized $\Phi$-functions on cubes with a stronger relation,
known as strong domination. This concept will play a key role in Theorem~\ref{th:sufficient},
where a sufficient condition for the boundedness of $M$ is formulated in terms
of strong domination. The following definition is a dyadic counterpart 
of~\cite[Definition~5.5.1]{DHHR11}. 
\begin{definition}\label{def:strong-domin}
Let $\varphi,\psi\in\Phi(\mathscr{D})$. We say that $\psi$ is {\it strongly dominated} by
$\varphi$, and write $\psi\ll\varphi$, if for every $B_1>0$ there exist $B_2>0$ such that 
\begin{equation}\label{eq:def-strong-domin}
\sum_{k=-\infty}^\infty\sum_{Q\in\mathcal{Q}_k}\mu(Q)\,\varphi(Q,2^k)\le B_1
\implies
\sum_{k=-\infty}^\infty\sum_{Q\in\mathcal{Q}_k}\mu(Q)\,\psi(Q,2^k)\le B_2
\end{equation}
for any sequence of families $\mathcal{Q}_k\in\mathfrak{S}$, $k\in\mathbb{Z}$. We 
allow some families $\mathcal{Q}_k$ to be empty; in such cases, the corresponding terms
in the above sums are understood to be zero. 
\end{definition}
Strong domination $\ll$ is called so because it is, at least, not weaker than 
domination~$\preceq$, as shown in~\cite[Lemma~5.5.5]{DHHR11} in the context of 
Euclidean cubes. For completeness, we reproduce the proof of this lemma in our context, 
while also correcting some misprints present in the original source.
\begin{lemma}\label{le:strong-weak}
Let $\varphi,\psi\in\Phi(\mathscr{D})$ be such that $\psi\ll\varphi$. Then 
$\psi\preceq\varphi$.
\end{lemma}
\begin{proof}
Let $\mathcal{Q}\in\mathfrak{S}$ and 
$\{t_Q\}_{Q\in\mathcal{Q}}\in l^{\varphi(Q)}(\mathcal{Q})$ satisfy
$\|\{t_Q\}_{Q\in\mathcal{Q}}\|_{l^{\varphi(Q)}(\mathcal{Q})}\le1$, which is equivalent
by the unit ball property in Lemma~\ref{le:spaces}(e) to
\[
\rho_{\varphi(Q)}\big(\{t_Q\}_{Q\in\mathcal{Q}}\big)
=\sum_{Q\in\mathcal{Q}}\mu(Q)\,\varphi(Q,|t_Q|)\le1.
\]
Assume without loss of generality that $t_Q\ne0$ for all 
$Q\in\mathcal{Q}$ (otherwise exclude all cubes $Q$ for which $t_Q=0$ from $\mathcal{Q}$; the
above semimodular sum will not change). Then for every $Q\in\mathcal{Q}$ there exists
$k_Q\in\mathbb{Z}$ such that $2^{k_Q}\le|t_Q|<2^{k_Q+1}$. Define the sequence 
of families $\mathcal{Q}_k\in\mathfrak{S}$, for all $k\in\mathbb{Z}$, by
$\mathcal{Q}_k:=\{Q\in\mathcal{Q}:k_Q=k\}$. Then 
\begin{equation}\label{eq:str-domin}
\sum_{k=-\infty}^\infty\sum_{Q\in\mathcal{Q}_k}\mu(Q)\,\varphi(Q,2^k)
\le\sum_{Q\in\mathcal{Q}}\mu(Q)\,\varphi(Q,|t_Q|)\le1.
\end{equation}

Since $\psi\ll\varphi$, by Definition~\ref{def:strong-domin} one may find $B_2>1$ for
$B_1:=1$ such that~\eqref{eq:def-strong-domin} holds. Choose an $m\in\mathbb{N}$ 
satisfying $2^{m-1}<B_2\le2^m$. Then inequality~\eqref{eq:str-domin} implies,
by the estimate $|t_Q|<2^{k_Q+1}$ for all $Q\in\mathcal{Q}$, the convexity of $\psi$
in the second variable, and assumption~\eqref{eq:def-strong-domin}, that 
\begin{align*}
\sum_{Q\in\mathcal{Q}}\mu(Q)\,\psi(Q,2^{-(m+1)}|t_Q|)
&\le\sum_{k=-\infty}^\infty\sum_{Q\in\mathcal{Q}_k}\mu(Q)\,\psi(Q,2^{k-m})
\\
&\le2^{-m}\sum_{k=-\infty}^\infty\sum_{Q\in\mathcal{Q}_k}\mu(Q)\,\psi(Q,2^k)
\le2^{-m}B_2\le1,
\end{align*}
which gives $\|\{t_Q\}_{Q\in\mathcal{Q}}\|_{l^{\psi(Q)}(\mathcal{Q})}\le2^{m+1}<4B_2$
by the unit ball property. A scaling argument proves that
\[
\big\|\{t_Q\}_{Q\in\mathcal{Q}}\big\|_{l^{\psi(Q)}(\mathcal{Q})}\le4B_2\,
\big\|\{t_Q\}_{Q\in\mathcal{Q}}\big\|_{l^{\varphi(Q)}(\mathcal{Q})}
\]
for all $\mathcal{Q}\in\mathfrak{S}$ and 
$\{t_Q\}_{Q\in\mathcal{Q}}\in l^{\varphi(Q)}(\mathcal{Q})$, therefore
$\psi\preceq\varphi$ by Definition~\ref{def:domination}.
\end{proof}

It has been observed in~\cite[Remark~5.5.4]{DHHR11} that under the $\Delta_2$-condition on
the (strongly) dominant function, it is possible to replace $2^k$ in the definition of
strong domination~\cite[Definition~5.5.1]{DHHR11} by $\alpha^k$ for any $\alpha>1$. 
However, this fact was not proved there. After reconstructing the proof, we found it to be 
nontrivial; so we prove here the analogous fact in the context of
Definition~\ref{def:strong-domin}.
\begin{lemma}\label{le:alpha^k}
Let $\varphi,\psi\in\Phi(\mathscr{D})$ and $\varphi$ satisfy the $\Delta_2$-condition. 
Suppose that $\alpha>1$. Then $\psi\ll\varphi$ if and only if for every $C_1>0$ there
exists $C_2>0$ such that 
\begin{equation}\label{eq:def-equiv}
\sum_{m=-\infty}^\infty\sum_{Q\in\mathcal{Q}_m}\mu(Q)\,\varphi(Q,\alpha^m)\le C_1
\implies
\sum_{m=-\infty}^\infty\sum_{Q\in\mathcal{Q}_m}\mu(Q)\,\psi(Q,\alpha^m)\le C_2
\end{equation}
for any sequence of families $\mathcal{Q}_m\in\mathfrak{S}$, $m\in\mathbb{Z}$,
allowing some $\mathcal{Q}_m$ to be empty.
\end{lemma}
\begin{proof}
We prove the lemma for $\alpha>2$, and discuss the case $1<\alpha<2$ at the conclusion 
of the proof.

\medskip\noindent
\emph{Part 1 (the necessity part): from $2$ to $\alpha$.}
Assume that $\psi\ll\varphi$ and take a sequence of families $\mathcal{Q}_m\in\mathfrak{S}$, 
$m\in\mathbb{Z}$, satisfying
\begin{equation}\label{ineq:assum-1}
\sum_{m=-\infty}^\infty\sum_{Q\in\mathcal{Q}_m}\mu(Q)\,\varphi(Q,\alpha^m)\le C_1.
\end{equation}
For every $m\in\mathbb{Z}$, define 
$k(m):=\min\{k\in\mathbb{Z}:\alpha^m\le2^k\}.$
Intuitively, $2^{k(m)}$ is the closest to $\alpha^m$ power of two from above. Clearly, 
\begin{equation}\label{eq:est}
2^{k(m)-1}<\alpha^m\le 2^{k(m)}
\quad\text{ for all }m\in\mathbb{Z},
\end{equation}
and $k(\cdot)$ is non-decreasing. Moreover, since $\alpha>2$, the function $k(\cdot)$ is
strictly increasing. Indeed, if $k(m_1)=k(m_2)$ for some $m_1\le m_2$, then 
estimate~\eqref{eq:est} gives 
\[
\alpha^{m_2}\le 2^{k(m_2)}=2\times 2^{k(m_1)-1}<2\alpha^{m_1},
\]
whence $m_2-m_1<\log_\alpha2<1$, so $m_1=m_2$. Thus, the sequence $\{k(m)\}_{m\in\mathbb{Z}}$
does not contain equal numbers.

In view of the $\Delta_2$-condition for $\varphi$ with constant $D$, the monotonicity
of $\varphi$ in the second variable, and
estimate~\eqref{eq:est}, the assumption in~\eqref{ineq:assum-1} implies 
\begin{align*}
\sum_{m=-\infty}^\infty\sum_{Q\in\mathcal{Q}_m}\mu(Q)\,\varphi(Q,2^{k(m)})   
&\le D\sum_{m=-\infty}^\infty\sum_{Q\in\mathcal{Q}_m}\mu(Q)\,\varphi(Q,2^{k(m)-1})
\\
&\le D\sum_{m=-\infty}^\infty\sum_{Q\in\mathcal{Q}_m}\mu(Q)\,\varphi(Q,\alpha^m)\le 
DC_1=:B_1.
\end{align*}
Define $\widetilde{\mathcal{Q}}_{k}:=\mathcal{Q}_m$ if $k=k(m)$ for some
$m\in\mathbb{Z}$, and $\widetilde{\mathcal{Q}}_{k}:=\emptyset$ for all 
$k\in\mathbb{Z}\setminus\{k(m)\}_{m\in\mathbb{Z}}$. Then the above estimate can be rewritten
as
\[
\sum_{k=-\infty}^\infty\sum_{Q\in\widetilde{\mathcal{Q}}_k}\mu(Q)\,\varphi(Q,2^k)\le B_1, 
\]
so applying Definition~\ref{def:strong-domin} together with estimate~\eqref{eq:est},
we conclude that
\begin{align*}
\sum_{m=-\infty}^\infty\sum_{Q\in\mathcal{Q}_m}\mu(Q)\,\psi(Q,\alpha^m)
&\le\sum_{m=-\infty}^\infty\sum_{Q\in\mathcal{Q}_m}\mu(Q)\,\psi(Q,2^{k(m)})
\\
&=\sum_{k=-\infty}^\infty\sum_{Q\in\widetilde{\mathcal{Q}}_k}\mu(Q)\,\psi(Q,2^k)
\le B_2,
\end{align*}
where $B_2$ does not depend on the sequence 
$\{\widetilde{\mathcal{Q}}_k\}_{k\in\mathbb{Z}}$ and hence on 
$\{\mathcal{Q}_m\}_{m\in\mathbb{Z}}$. Thus we proved~\eqref{eq:def-equiv} with $C_2:=B_2$.

\medskip\noindent
\emph{Part 2 (the sufficiency part): from $\alpha$ to $2$.}
Conversely, assume~\eqref{eq:def-equiv} and take a sequence of families
$\mathcal{Q}_k\in\mathfrak{S}$, $k\in\mathbb{Z}$, satisfying
\begin{equation}\label{ineq:assum-2}
\sum_{k=-\infty}^\infty\sum_{Q\in\mathcal{Q}_k}\mu(Q)\,\varphi(Q,2^k)\le B_1.
\end{equation}
Define the number $m(k):=\min\{m\in\mathbb{Z}:2^k\le\alpha^m\}$ for each $k\in\mathbb{Z}$,
so that $\alpha^{m(k)}$ is the closest to $2^k$ power of $\alpha$ from above. By this
definition,
\begin{equation}\label{eq:est-2}
\alpha^{m(k)-1}<2^k\le\alpha^{m(k)}
\quad\text{for all }k\in\mathbb{Z},
\end{equation}
and $m(\cdot)$ is non-decreasing. Since $2<\alpha$, several consecutive values of $k$ may
have the same $m(k)$, but the number of such $k$ never exceeds $N:=[\log_2\alpha]+1$. Indeed,
if $m(k_1)=m(k_2)$ for some $k_1\le k_2$, then it follows from estimate~\eqref{eq:est-2} that
\[
2^{k_2}\le \alpha^{m(k_2)}=\alpha\times\alpha^{m(k_1)-1}<\alpha\,2^{k_1},
\]
whence $k_2-k_1<\log_2\alpha$. Thus the maximum possible distance between such points $k_1$
and $k_2$ equals $[\log_2\alpha]$, hence the maximum possible number of integers in the 
interval $[k_1,k_2]$ is $[\log_2\alpha]+1$. 

Another consequence of $2<\alpha$ is that for every $m\in\mathbb{Z}$, there is 
$k\in\mathbb{Z}$ such that $m(k)=m$. This follows from the observation that 
$m=m(k)$ if and only if 
\[
(m-1)\log_2\alpha<k\le m\log_2\alpha
\]
by estimate~\eqref{eq:est-2}, and since $m\log_2\alpha-(m-1)\log_2\alpha=\log_2\alpha>1$,
there is always at least one integer $k$ in the above interval. 

Therefore, the collection
\[
\mathfrak{C}_m:=\{\mathcal{Q}_k\ :\ m(k)=m,\,k\in\mathbb{Z}\},
\]
for any $m\in\mathbb{Z}$, is nonempty and contains at most $N$ elements.
Re-denote them by $\mathfrak{C}_m:=\{\widetilde{\mathcal{Q}}_m^i\}_{i=1}^N$,
assuming that we may add some empty $\widetilde{\mathcal{Q}}_m^i$ to make the number of
families in $\mathfrak{C}_m$ exactly $N$. 
Consider the finite number of sequences 
\[
\{\widetilde{\mathcal{Q}}_m^{i}\}_{m\in\mathbb{Z}},
\quad
i=1,\ldots,N.
\]
By construction, 
\begin{enumerate}
\item[1)] each family $\mathcal{Q}_k$, $k\in\mathbb{Z}$, belongs to a unique 
sequence $\{\widetilde{\mathcal{Q}}_m^{i}\}_{m\in\mathbb{Z}}$; 
\item[2)] every sequence $\{\widetilde{\mathcal{Q}}_m^{i}\}_{m\in\mathbb{Z}}$
contains not more than one family $\mathcal{Q}_k$ with $m(k)=m$ for each $m\in\mathbb{Z}$;
\item[3)] each $\widetilde{\mathcal{Q}}_m^i$ is either
some $\mathcal{Q}_k$ with $m(k)=m$ or the empty set.
\end{enumerate}
Now, noting that for all $m\in\mathbb{Z}$ and $Q\in\mathscr{D}$ there holds
\begin{equation}\label{eq:estimate}
\varphi(Q,\alpha^m)=\varphi(Q,2^{\log_2\alpha}\alpha^{m-1})
\le\varphi(Q,2^N\alpha^{m-1})\le D^N\,\varphi(Q,\alpha^{m-1})
\end{equation}
due to the $\Delta_2$-condition for $\varphi$, we apply
estimate~\eqref{eq:est-2} and assumption~\eqref{ineq:assum-2} to
obtain, for each $i=1,\ldots,N$, that
\begin{align*}
\sum_{m=-\infty}^\infty\sum_{Q\in\widetilde{\mathcal{Q}}_m^i}\mu(Q)\,\varphi(Q,\alpha^m)
&\le D^N
\sum_{m=-\infty}^\infty\sum_{Q\in\widetilde{\mathcal{Q}}_m^i}\mu(Q)\,\varphi(Q,\alpha^{m-1})
\\
&\le D^N
\sum_{k=-\infty}^\infty\sum_{Q\in\mathcal{Q}_k}\mu(Q)\,\varphi(Q,\alpha^{m(k)-1})
\\
&\le D^N
\sum_{k=-\infty}^\infty\sum_{Q\in\mathcal{Q}_k}\mu(Q)\,\varphi(Q,2^k)
\le D^N B_1=:C_1,
\end{align*}
and hence
\[
\sum_{m=-\infty}^\infty\sum_{Q\in\widetilde{\mathcal{Q}}_m^i}\mu(Q)\,\psi(Q,\alpha^m)
\le C_2 
\]
holds by our assumption~\eqref{eq:def-equiv} with a constant $C_2>0$ independent of $i$
and families $\{\widetilde{\mathcal{Q}}_m^{i}\}_{m\in\mathbb{Z}}$.
Summing the above inequalities over $i=1,\ldots,N$, we conclude that
\begin{align*}
\sum_{k=-\infty}^\infty\sum_{Q\in\mathcal{Q}_k}\mu(Q)\,\psi(Q,2^k)
&\le\sum_{k=-\infty}^\infty\sum_{Q\in\mathcal{Q}_k}\mu(Q)\,\psi(Q,\alpha^{m(k)})
\\
&=\sum_{i=1}^N
\sum_{m=-\infty}^\infty\sum_{Q\in\widetilde{\mathcal{Q}}_m^i}\mu(Q)\,\psi(Q,\alpha^m)
\le N C_2=:B_2,
\end{align*}
where $B_2>0$ does not depend on the family 
$\{\mathcal{Q}_k\}_{k\in\mathbb{Z}}$. Thus, we proved~\eqref{eq:def-strong-domin} and 
therefore $\psi\ll\varphi$.

\medskip\noindent
\emph{The situation $1<\alpha<2$ is considered analogously}.
In this case, the necessity part with passing from $2$ to $\alpha$ will be similar to 
{\it Part~2} of the above proof, since we pass to a smaller number $\alpha<2$.
One just has to replace $\alpha$ with $2$ and $m$ with $k$ in \emph{Part~2} and repeat 
the argument almost verbatim, with the simpler estimate
$\varphi(Q,2^k)\le D\,\varphi(Q,2^{k-1})$ in place of~\eqref{eq:estimate}.
The sufficiency part will follow by the same argument as in {\it Part~1}.
\end{proof}

Under the $\Delta_2$-condition on the (strongly) dominated function, it suffices to
verify Definition~\ref{def:strong-domin} for $B_1:=1$ only. This is another
observation related to the definition of strong domination, left without proof 
in~\cite[Remark~5.5.4]{DHHR11}. We check this statement for the reader's convenience.
\begin{lemma}\label{le:strongwith1}
Let $\varphi,\psi\in\Phi(\mathscr{D})$ and $\psi$ satisfy the $\Delta_2$-condition. 
Suppose that there exists $C>0$ such that 
\begin{equation}\label{eq:claim-1}
\sum_{m=-\infty}^\infty\sum_{Q\in\mathcal{Q}_m}\mu(Q)\,\varphi(Q,2^m)\le1
\implies
\sum_{m=-\infty}^\infty\sum_{Q\in\mathcal{Q}_m}\mu(Q)\,\psi(Q,2^m)\le C
\end{equation}
for any sequence of families $\mathcal{Q}_m\in\mathfrak{S}$, $m\in\mathbb{Z}$ 
(some $\mathcal{Q}_m$ may be empty). Then $\psi\ll\varphi$.
\end{lemma}
\begin{proof}
To check~\eqref{eq:def-strong-domin}, take a sequence of families
$\mathcal{Q}_k\in\mathfrak{S}$, $k\in\mathbb{Z}$, such that
\[
\sum_{k=-\infty}^\infty\sum_{Q\in\mathcal{Q}_k}\mu(Q)\,\varphi(Q,2^k)\le B_1
\]
for some $B_1>0$. Choose a number $l\in\mathbb{N}$ satisfying
$B_1\le2^l$. Then it follows by the convexity of $\varphi$ in the second variable that 
\[
\sum_{k=-\infty}^\infty\sum_{Q\in\mathcal{Q}_k}\mu(Q)\,\varphi(Q,2^{k-l})
\le\frac1{2^l}
\sum_{k=-\infty}^\infty\sum_{Q\in\mathcal{Q}_k}\mu(Q)\,\varphi(Q,2^k)
\le\frac{B_1}{2^l}\le1.
\]
Denote $m:=k-l$ and $\widetilde{\mathcal{Q}}_m:=\mathcal{Q}_k$ for all $k\in\mathbb{Z}$. 
In the new notation, the above estimate can be rewritten as
\[
\sum_{m=-\infty}^\infty\sum_{Q\in\widetilde{\mathcal{Q}}_m}\mu(Q)\,\varphi(Q,2^m)\le1,
\]
and this implies by~\eqref{eq:claim-1} and the $\Delta_2$-condition for 
$\psi$ with constant $D$ that
\begin{align*}
\sum_{k=-\infty}^\infty\sum_{Q\in\mathcal{Q}_k}\mu(Q)\,\psi(Q,2^k)
&\le D^l \sum_{k=-\infty}^\infty\sum_{Q\in\mathcal{Q}_k}\mu(Q)\,\psi(Q,2^{k-l})
\\
&=D^l \sum_{m=-\infty}^\infty\sum_{Q\in\widetilde{\mathcal{Q}}_m}\mu(Q)\,\psi(Q,2^m)
\le D^l\,C=:B_2.
\end{align*}
Note that $B_2>0$ does not depend on the exact sequence 
$\{\mathcal{Q}_k\}_{k\in\mathbb{Z}}$, since
the constant~$C$ is independent of $\{\widetilde{\mathcal{Q}}_m\}_{m\in\mathbb{Z}}$. Thus 
$\psi\ll\varphi$ by Definition~\ref{def:strong-domin}.
\end{proof}

\section{Self-Improving Class $\mathscr{A^D}$}\label{sec:A}
In this section, we introduce a special class $\mathscr{A^D}$ of generalized 
$\Phi$-functions $\varphi$ on a space of homogeneous type $X$, characterized by the uniform
boundedness of the averaging operators $T_\mathcal{Q}$, $\mathcal{Q}\in\mathfrak{S}$, 
on the corresponding spaces $L^\varphi(X)$. Our class $\mathscr{A^D}$ is a dyadic variant
of the class $\mathcal{A}$ introduced by Diening in~\cite[Section~3]{LD05}, whose idea was 
to generalize the concept of the Muckenhoupt classes $A_p$, $1<p<\infty$, from classical
weighted Lebesgue spaces to Musielak--Orlicz spaces. We show that the class $\mathscr{A^D}$
is ``self-improving'' (Theorem~\ref{th:self-improvment}), give its characterization 
in terms of domination between generalized $\Phi$-functions over cubes
(Theorem~\ref{th:1char}), and obtain a refined domination property of $\mathscr{A^D}$
(Theorem~\ref{th:refinement}).

\subsection{Definition and the duality property of the class $\mathscr{A^D}$}

To define the class $\mathscr{A^D}$, we use the averaging operators over families of 
dyadic cubes introduced in Definition~\ref{def:averaging-operators}.
\begin{definition}
We say that $\varphi\in\Phi(X)$ is in the {\it class $\mathscr{A^D}$} if the averaging
operators $T_\mathcal{Q}$ are bounded on $L^\varphi(X)$ uniformly over all families 
$\mathcal{Q}\in\mathfrak{S}$. The smallest constant $C>0$ for which 
\[
\|T_\mathcal{Q}f\|_\varphi\le C\,\|f\|_\varphi
\]
holds for all $\mathcal{Q}\in\mathfrak{S}$ and $f\in L^\varphi(X)$ will be called
the \emph{$\mathscr{A^D}$-constant} of $\varphi$.
\end{definition}
Proper functions $\varphi\in\Phi(X)$ belong to $\mathscr{A^D}$ together with their
conjugates $\varphi^*$, as the next lemma shows (cf.~\cite[Lemma~5.2.2]{DHHR11}). 
In other words, the class $\mathscr{A^D}$ behaves well with respect to duality.
\begin{lemma}\label{le:A-duality}
Let $\varphi\in\Phi(X)$ be proper. Then $\varphi\in\mathscr{A^D}$ if and only if
$\varphi^*\in\mathscr{A^D}$.
\end{lemma}
\begin{proof}
Let $\varphi\in\mathscr{A^D}$. Since $\varphi$ is proper, its conjugate function 
$\varphi^*$ is also proper. Hence, we can apply 
Lemma~\ref{le:norm-conj} to $\varphi^*$, observing that $(\varphi^*)^*=\varphi$. Then
this, self-duality of averaging operators~\eqref{eq:self-dual}, H\"older's inequality 
in Lemma~\ref{le:Holder}, and the assumption $\varphi\in\mathscr{A^D}$ imply, for any 
$\mathcal{Q}\in\mathfrak{S}$ and $f\in L^{\varphi^*}(X)$, that
\begin{align*}
\|T_\mathcal{Q}f\|_{\varphi^*}
\le\sup_{\substack{g\in L^\varphi(X):\\\|g\|_\varphi\le1}}
\langle T_\mathcal{Q}f,|g|\rangle 
&=\sup_{\substack{g\in L^\varphi(X):\\\|g\|_\varphi\le1}}
\langle|f|,T_\mathcal{Q}g\rangle
\\
&\le2\,\|f\|_{\varphi^*}
\sup_{\substack{g\in L^\varphi(X):\\\|g\|_\varphi\le1}}
\|T_\mathcal{Q}g\|_\varphi 
\le 2C\,\|f\|_{\varphi^*},
\end{align*}
where $C$ is the $\mathscr{A^D}$-constant of $\varphi$. Thus $\varphi^*\in\mathscr{A^D}$. 
For the converse claim, it suffices to replace $\varphi$ with $\varphi^*$ 
and use $(\varphi^*)^*=\varphi$ in the above estimate.
\end{proof}

\subsection{Lemma about perturbations of a cube-constant function}

To obtain finer properties of the class $\mathscr{A^D}$---such as the ``self-improvement''
property---we study a wider class $\mathscr{A_\infty^D}$, which is a dyadic variant of 
Diening's class $\mathcal{A}_\infty$ introduced in~\cite[Section~5]{LD05} as an analogue of
the Muckenhoupt class $A_\infty$.
\begin{definition}\label{def:A-inf}
By $\mathscr{A^D_\infty}$ we denote the set of all $\varphi\in\Phi(X)$ with the following
property: For every $\alpha\in(0,1)$ there exists $\beta\in(0,1)$ such that if $N\subseteq X$
(measurable) and $\mathcal{Q}\in\mathfrak{S}$ satisfy 
\begin{equation}\label{eq:A-inf-if}
\mu(N\cap Q)\le\alpha\,\mu(Q)\quad
\text{for all }Q\in\mathcal{Q},
\end{equation}
then
\begin{equation}\label{eq:A-inf-then}
\bigg\|\sum_{Q\in\mathcal{Q}}t_Q\,\chi_{N\cap Q}\bigg\|_\varphi\le\beta\,
\bigg\|\sum_{Q\in\mathcal{Q}}t_Q\,\chi_{Q}\bigg\|_\varphi
\end{equation}
for any sequence $\{t_Q\}_{Q\in\mathcal{Q}}\in\mathbb{C}^\mathcal{Q}$. If 
$\varphi\in\mathscr{A^D_\infty}$, 
the smallest constant $\beta$ for $\alpha=\frac12$ is called the 
{\it $\mathscr{A^D_\infty}$-constant} of $\varphi$.
\end{definition}
Condition~\eqref{eq:A-inf-then} admits the following interpretation.
A family $\mathcal{Q}\in\mathfrak{S}$ and a sequence
$\{t_Q\}_{Q\in\mathcal{Q}}\in\mathbb{C}^\mathcal{Q}$ generate the \emph{cube-constant} 
function $\sum_{Q\in\mathcal{Q}}t_Q\,\chi_{Q}$ on $\mathcal{Q}$, that is, the function with
the constant values $t_Q$ on the corresponding cubes $Q\in\mathcal{Q}$ and zero otherwise.
We fix a family $\mathcal{Q}\in\mathfrak{S}$ and a set $N$ that is ``non-expansive'' with
respect to $\mathcal{Q}$ in the sense of inequalities~\eqref{eq:A-inf-if}. 
Condition~\eqref{eq:A-inf-then} says that the $L^\varphi$-norms of all cube-constant
functions on $\mathcal{Q}$ restricted to the non-expansive set $N$ are uniformly controlled
by the $L^\varphi$-norms of their non-restricted originals---if $\varphi$ is in the 
class $\mathscr{A^D_\infty}$.

We next show that the class $\mathscr{A^D_\infty}$ is wider than some nice part of
$\mathscr{A^D}$. Our Lemma~\ref{le:A-A-inf} is an analogue 
of~\cite[Corollary~5.4.10]{DHHR11} obtained by 
compressing~\cite[Lemma~5.4.8(b) and Lemma~5.4.9]{DHHR11} into one proof.
\begin{lemma}\label{le:A-A-inf}
If $\varphi\in\mathscr{A^D}$ satisfies the $\Delta_2$-condition, 
then $\varphi\in\mathscr{A^D_\infty}$.
\end{lemma}
\begin{proof}
Let $\alpha\in(0,1)$ and let $N$, $\mathcal{Q}$ satisfy~\eqref{eq:A-inf-if}. Take an 
arbitrary sequence $\{t_Q\}_{Q\in\mathcal{Q}}$. Define 
$P:=\bigcup_{Q\in\mathcal{Q}}(Q\setminus N)$ and 
$f:=\sum_{Q\in\mathcal{Q}}t_Q\,\chi_{P\cap Q}$. 
Let $\|f\|_\varphi<\infty$ (otherwise $\|\sum_{Q\in\mathcal{Q}}t_Q\,\chi_Q\|_\varphi=\infty$
and~\eqref{eq:A-inf-then} holds trivially). Then for any $Q\in\mathcal{Q}$,
relation~\eqref{eq:A-inf-if} implies
\[
\mu(P\cap Q)=\mu(Q)-\mu(N\cap Q)\ge(1-\alpha)\,\mu(Q),
\]
and thus
\[
M_Qf=t_Q\,\frac{\mu(P\cap Q)}{\mu(Q)}\ge t_Q\,(1-\alpha).
\]
Therefore, since $\varphi\in\mathscr{A^D}$,
\begin{equation}\label{eq:ineq1}
(1-\alpha)\,\bigg\|\sum_{Q\in\mathcal{Q}}t_Q\,\chi_Q\bigg\|_\varphi
\le\bigg\|\sum_{Q\in\mathcal{Q}}(M_Qf)\,\chi_Q\bigg\|_\varphi=\|T_\mathcal{Q}f\|_\varphi
\le C\,\|f\|_\varphi,
\end{equation}
where $C$ is the $\mathscr{A^D}$-constant of $\varphi$. 

Assume next, without loss of 
generality, that $\|\sum_{Q\in\mathcal{Q}}t_Q\,\chi_Q\|_\varphi=1$. Then~\eqref{eq:ineq1}
implies $\|f\|_\varphi\ge(1-\alpha)/C$. It follows from Lemma~\ref{le:delta2}(d) that there 
exists a number $\beta_0\in(0,1)$, depending only on $\alpha$, $C$, and the 
$\Delta_2$-constant of $\varphi$, such that $\rho_\varphi(f)\ge\beta_0$. 
Since $N\cap P=\emptyset$, and 
$\rho_\varphi(\sum_{Q\in\mathcal{Q}}t_Q\,\chi_Q)=1$ by the unit sphere property 
in~Lemma~\ref{le:delta2}(c), we have
\[
\rho_\varphi\bigg(\sum_{Q\in\mathcal{Q}}t_Q\,\chi_{N\cap Q}\bigg)=
\rho_\varphi\bigg(\sum_{Q\in\mathcal{Q}}t_Q\,\chi_Q\bigg)-
\rho_\varphi\bigg(\sum_{Q\in\mathcal{Q}}t_Q\,\chi_{P\cap Q}\bigg)\le1-\beta_0.
\]
Finally, Lemma~\ref{le:delta2}(e) ensures that there exists $\beta\in(0,1)$, which
depends only on $\beta_0$ and the $\Delta_2$-constant of $\varphi$, such that 
\[
\bigg\|\sum_{Q\in\mathcal{Q}}t_Q\,\chi_{N\cap Q}\bigg\|_\varphi\le\beta.
\]
A scaling argument then yields~\eqref{eq:A-inf-then}.
\end{proof}

The following lemma about perturbations of cube-constant functions in the spaces
$L^\varphi(X)$ with $\varphi\in\mathscr{A^D_\infty}$ is the deepest result of this section.
Building on the Calder\'on--Zygmund-type decomposition in Lemma~\ref{le:CZ},
it extends \cite[Lemma~5.4.12]{DHHR11} to spaces of homogeneous type.
\begin{lemma}\label{le:perturbation}
Let $\varphi\in\Phi(X)$ be in $\mathscr{A^D_\infty}$. Then there exist $\delta\in(0,1)$ and
$G>2$ such that 
\[
\bigg\| \sum_{Q\in\mathcal{Q}}t_Q\bigg|\frac{f}{M_Q f}\bigg|^\delta\chi_Q\bigg\|_\varphi
\le
G\,\bigg\|\sum_{Q\in\mathcal{Q}}t_Q\chi_Q\bigg\|_\varphi
\]
for all $\mathcal{Q}\in\mathfrak{S}$, all $\{t_Q\}_{Q\in\mathcal{Q}}$ with $t_Q\ge0$, and all
$f\in L^1_\textnormal{loc}(X)$ with $M_Qf\ne0$, $Q\in\mathcal{Q}$.
\end{lemma}
\begin{proof}
Take a family $\mathcal{Q}\in\mathfrak{S}$, a non-negative sequence 
$\{t_Q\}_{Q\in\mathcal{Q}}$, and an $f\in L^1_\text{loc}(X)$ with $M_Qf\ne0$ for all
$Q\in\mathcal{Q}$. Fix a dyadic cube $Q\in\mathcal{Q}$ and find the number 
$t\in\{1,\ldots,K\}$ of the dyadic grid for which $Q\in\mathscr{D}^t$. 
Define $f_Q\in L^1(X)$ by $f_Q:=f\chi_Q$.

\medskip\noindent
\emph{Step 1: Smart decomposition of $Q$ into nested sets.}
We claim that $Q$ coincides with a certain superlevel set of the $t$-dyadic maximal
function of $f_Q$; namely,
\begin{equation}\label{eq:cube-Q}
\left\{x\in X\ :\ M^{\mathscr{D}^t}f_Q(x)>(1-\varepsilon)\, M_Qf\right\}=Q,
\end{equation}
where $\varepsilon\in(0,1)$ is the parameter of the dyadic system $\mathscr{D}$
associated with $X$ that is specified in Theorem~\ref{th:HK}(d). Indeed, if $x\in Q$, 
then by definition 
$M^{\mathscr{D}^t}f_Q(x)\ge M_Qf_Q=M_Qf>(1-\varepsilon)\, M_Qf$ since $M_Qf>0$.
Conversely, let $ M^{\mathscr{D}^t}f_Q(x)>(1-\varepsilon)\, M_Qf$. Then there exists
$Q'\in\mathscr{D}^t$ such that $x\in Q'$ and 
\[
\frac{1}{\mu(Q')}\int_{Q'\cap Q}|f(y)|\,d\mu(y)=M_{Q'}f_Q>(1-\varepsilon)\,M_Qf>0.
\]
This implies that $Q'\cap Q\ne\emptyset$, hence either $Q'\subseteq Q$ or
$Q\varsubsetneq Q'$. The former case yields $x\in Q$. The latter case is impossible, 
since $Q\varsubsetneq Q'$ would result in
\[
\frac{\mu(Q)}{\mu(Q')}=\frac{M_{Q'}f_Q}{M_Qf}>1-\varepsilon,
\]
which contradicts Lemma~\ref{le:cones}(b). We thus proved~\eqref{eq:cube-Q}.

Next, we introduce the sequence of superlevel sets
\[
E_Q^k:=\left\{x\in X\ :\ M^{\mathscr{D}^t}f_Q(x)>
(1-\varepsilon)\left(\frac2\varepsilon\right)^k
M_Qf\right\},
\]
where $k\in\mathbb{N}_0:=\mathbb{N}\cup\{0\}$. They are nested:
\[
E_Q^{k+1}\subseteq E_Q^k\subseteq\ldots\subseteq E_Q^0=Q.
\]
Each set $E_Q^k$, being a superlevel set of $M^{\mathscr{D}^t}f_Q$, decomposes into 
the disjoint union of all its Calder\'on--Zygmund cubes provided by Lemma~\ref{le:CZ}
(and thus is measurable).

\medskip\noindent
\emph{Step 2: Auxiliary claim:} Every Calder\'on--Zygmund cube $V$ of 
the superlevel set $E_Q^k$, $k\in\mathbb{N}_0$, satisfies
\begin{equation}\label{eq:non-expansiveness}
\mu(E_Q^{k+1}\cap V)\le\frac12\,\mu(V).
\end{equation}
\emph{Proof of auxiliary claim.}
If $E_Q^{k+1}\cap V=\emptyset$, the inequality trivially holds. Suppose that 
$E_Q^{k+1}\cap V\ne\emptyset$. Let $W$ be a Calder\'on--Zygmund cube of $E_Q^{k+1}$ which
intersects $V$. By Lemma~\ref{le:CZ}(ii), there holds
$M_Wf_Q>(1-\varepsilon)(\frac2\varepsilon)^{k+1}\, M_Qf$, or
\begin{equation}\label{in:W}
\mu(W)\,M_Qf<
\frac1{1-\varepsilon}\left(\frac\varepsilon2\right)^{k+1}\int_W|f_Q(y)|\,d\mu(y).
\end{equation}
Since two cubes $W$ and $V$ from $\mathscr{D}^t$ have nonempty intersection, one of them is
contained in the other one; $E_Q^{k+1}\subseteq E_Q^k$ implies $W\subseteq V$ (otherwise
$V\varsubsetneq W\subseteq E_Q^{k+1}\subseteq E_Q^k$ would contradict the maximality of $V$
in $E_Q^k$, see Lemma~\ref{le:CZ}(iii)).
Therefore, $E_Q^{k+1}\cap V$ is the union of all Calder\'on--Zygmund cubes $W$ of $E_Q^{k+1}$
which intersect $V$. Summing inequalities~\eqref{in:W} over all such cubes $W$ implies
\begin{equation}\label{in:V}
\mu(E_Q^{k+1}\cap V)\,M_Qf\le
\frac1{1-\varepsilon}\left(\frac\varepsilon2\right)^{k+1}\int_V|f_Q(y)|\,d\mu(y).
\end{equation}
Next, we apply the upper estimate in Theorem~\ref{le:CZ}(ii) to the Calder\'on--Zygmund cube
$V$ of $E_Q^k$ to obtain
\[
\int_V|f_Q(y)|\,d\mu(y)=\mu(V)\,M_Vf_Q\le\frac{1-\varepsilon}\varepsilon
\left(\frac2\varepsilon\right)^k\mu(V)\,M_Qf.
\]
Combined with~\eqref{in:V} and $M_Qf>0$, this gives the desired 
estimate~\eqref{eq:non-expansiveness}.
The auxiliary claim is proved, and we continue the main proof now.$\quad\square$

\medskip\noindent
\emph{Step 3: The class $\mathscr{A^D_\infty}$ comes into play.}
Let $\{V_{Q,j}^k\}_j$ be the collection of all Calder\'on--Zygmund cubes of 
$E_Q^k$, $k\in\mathbb{N}_0$. These cubes are disjoint, and by the auxiliary claim,
\[
\mu(E_Q^{k+1}\cap V^k_{Q,j})\le\frac12\,\mu(V^k_{Q,j})
\]
for all $j$. Since $E_Q^k\subseteq Q$ and the cubes $Q\in\mathcal{Q}$ are pairwise disjoint, 
the collection $\{V^k_{Q,j}\}_{Q,j}$ consists of pairwise disjoint cubes for each fixed $k$.
The set $F^k:=\bigcup_{Q\in\mathcal{Q}}E_Q^k$ is measurable and satisfies 
\[
\mu(F^{k+1}\cap V^k_{Q,j})=\mu(E_Q^{k+1}\cap V^k_{Q,j})\le\frac12\,\mu(V^k_{Q,j})
\]
for all $Q\in\mathcal{Q}$ and $j$. Hence for every $k\in\mathbb{N}_0$, the set $F^{k+1}$ and
the family $\{V^k_{Q,j}\}_{Q,j}\in\mathfrak{S}$ are exactly as in Definition~\ref{def:A-inf}
with $\alpha=\frac12$. Thus it follows from 
$\varphi\in\mathscr{A^D_\infty}$ that 
\begin{equation}\label{eq:ineq-beta}
\bigg\|\sum_{Q\in\mathcal{Q}}\sum_{j}t_Q\,\chi_{F^{k+1}\cap V^k_{Q,j}}\bigg\|_\varphi
\le\beta\,\bigg\|\sum_{Q\in\mathcal{Q}}\sum_{j}t_Q\,\chi_{V^k_{Q,j}}\bigg\|_\varphi,
\end{equation}
where $\beta\in(0,1)$ is the $\mathscr{A^D_\infty}$-constant of $\varphi$, holds for 
all $k\in\mathbb{N}_0$. Note that we applied Definition~\ref{def:A-inf} to the sequence
whose terms corresponding to the cubes of the collection $\{V^k_{Q,j}\}_{j}$ are all equal
to $t_Q$, $Q\in\mathcal{Q}$.

Since $E^k_Q=\bigcup_{j}V^k_{Q,j}$, inequality~\eqref{eq:ineq-beta} can be rewritten
as
\[
\bigg\|\sum_{Q\in\mathcal{Q}}t_Q\,\chi_{F^{k+1}\cap E_Q^k}\bigg\|_\varphi
\le\beta\,\bigg\|\sum_{Q\in\mathcal{Q}}t_Q\,\chi_{E_Q^k}\bigg\|_\varphi.
\]
The monotonicity of $E_Q^k$ implies $F^{k+1}\cap E_Q^k=E_Q^{k+1}\cap E_Q^k=E_Q^{k+1}$
and hence
\[
\bigg\|\sum_{Q\in\mathcal{Q}}t_Q\,\chi_{E_Q^{k+1}}\bigg\|_\varphi
\le\beta\,\bigg\|\sum_{Q\in\mathcal{Q}}t_Q\,\chi_{E_Q^k}\bigg\|_\varphi,
\]
which gives by induction, for any $k\in\mathbb{N}_0$,
\[
\bigg\|\sum_{Q\in\mathcal{Q}}t_Q\,\chi_{E_Q^k}\bigg\|_\varphi
\le\beta^k\,\bigg\|\sum_{Q\in\mathcal{Q}}t_Q\,\chi_{E_Q^{0}}\bigg\|_\varphi
=\beta^k\bigg\|\sum_{Q\in\mathcal{Q}}t_Q\,\chi_{Q}\bigg\|_\varphi.
\]

\medskip\noindent
\emph{Step 4: Perturbation.} Since $f_Q\in L^1(X)$, it follows from 
Lemma~\ref{le:maxfun} that $|f_Q|\le M^{\mathscr{D}^t}f_Q$ holds $\mu$-almost everywhere.
This, the definition of $E_Q^{k+1}$, and the above inductive inequality imply, 
for all $k\in\mathbb{N}_0$ and all $\delta>0$, that
\begin{align*}
\bigg\|\sum_{Q\in\mathcal{Q}}t_Q\bigg|\frac{f_Q}{M_Q f}\bigg|^\delta
\chi_{E_Q^k\setminus E_Q^{k+1}}\bigg\|_\varphi
&\le\bigg\|\sum_{Q\in\mathcal{Q}}t_Q
\bigg(\frac{M^{\mathscr{D}^t}f_Q}{M_Q f}\bigg)^\delta
\chi_{E_Q^k\setminus E_Q^{k+1}}\bigg\|_\varphi \\
&\le\bigg\|\sum_{Q\in\mathcal{Q}}t_Q
\bigg((1-\varepsilon)\bigg(\frac2\varepsilon\bigg)^{k+1}\bigg)^\delta
\chi_{E_Q^k\setminus E_Q^{k+1}}\bigg\|_\varphi \\
&\le\bigg(\frac2\varepsilon\bigg)^{\delta(k+1)}
\,\bigg\|\sum_{Q\in\mathcal{Q}}t_Q\,\chi_{E_Q^k}\bigg\|_\varphi \\
&\le\beta^k\left(\frac2\varepsilon\right)^{\delta(k+1)}
\,\bigg\|\sum_{Q\in\mathcal{Q}}t_Q\,\chi_{Q}\bigg\|_\varphi.
\end{align*}
We fix $\delta\in(0,1)$ such that $\gamma:=\beta\,(\frac2\varepsilon)^\delta<1$. Note that
$\delta$ depends only on $\beta$, which is the $\mathscr{A^D_\infty}$-constant of $\varphi$,
and on the parameter $\varepsilon$ of the system $\mathscr{D}$. Then for all 
$k\in\mathbb{N}_0$, 
\[
\bigg\|\sum_{Q\in\mathcal{Q}}t_Q\bigg|\frac{f_Q}{M_Q f}\bigg|^\delta
\chi_{E_Q^k\setminus E_Q^{k+1}}\bigg\|_\varphi
\le\frac{2\gamma^k}{\varepsilon}
\,\bigg\|\sum_{Q\in\mathcal{Q}}t_Q\,\chi_{Q}\bigg\|_\varphi.
\]
Using the monotonicity of $E_Q^k$, Lemma~\ref{le:spaces}(d), and the above inequality,
we finally get
\begin{align*}
\bigg\|\sum_{Q\in\mathcal{Q}}t_Q\bigg|\frac{f}{M_Qf}\bigg|^\delta\chi_Q\bigg\|_\varphi
&=\bigg\|\sum_{k=0}^\infty\sum_{Q\in\mathcal{Q}}t_Q
\bigg|\frac{f_Q}{M_Qf}\bigg|^\delta\chi_{E_Q^k\setminus E_Q^{k+1}}
\bigg\|_\varphi 
\\
&\le\sum_{k=0}^\infty\,\bigg\|\sum_{Q\in\mathcal{Q}}t_Q
\bigg|\frac{f_Q}{M_Qf}\bigg|^\delta\chi_{E_Q^k\setminus E_Q^{k+1}}
\bigg\|_\varphi 
\\
&\le\sum_{k=0}^\infty\frac{2\gamma^k}{\varepsilon}
\,\bigg\|\sum_{Q\in\mathcal{Q}}t_Q\,\chi_Q\bigg\|_\varphi \\
&=\frac{2}{\varepsilon(1-\gamma)}
\,\bigg\|\sum_{Q\in\mathcal{Q}}t_Q\,\chi_Q\bigg\|_\varphi.
\end{align*}
This is the claim with $G:=\frac2{\varepsilon(1-\gamma)}>2$.
\end{proof}

\subsection{Self-improvement property of the class $\mathscr{A^D}$}

With the help of Lemma~\ref{le:perturbation} about perturbations of a cube-constant 
function, we now prove the ``self-improvement'' property of the class $\mathscr{A^D}$.
It says, in essence, that if we take a ``nice enough'' function $\varphi\in\mathscr{A^D}$ 
and slow down the growth of each $\Phi$-function $t\mapsto\varphi(x,t)$ just a little 
bit---by replacing $t$ with $t^s$ for some $s$ slightly below one---the resulting function 
$(x,t)\mapsto\varphi(x,t^s)$ still belongs to the class $\mathscr{A^D}$.
The following theorem generalizes the self-improvement property of 
Diening's class $\mathcal{A}$, established in~\cite[Theorems~5.4.13 and 5.4.15]{DHHR11} 
as an analogue of the left-openness of the classical Muckenhoupt classes.
\begin{theorem}\label{th:self-improvment}
Suppose that $\varphi\in\mathscr{A^D}$ is proper and $\varphi^*$ satisfies the
$\Delta_2$-condition. Then there exists $s_0\in(0,1)$ such that for any $s\in[s_0,1]$,
all functions $\psi\in\Phi(X)$ satisfying $\psi(x,t)\approx\varphi(x,t^s)$
uniformly in $x\in X$ and $t\ge0$ belong to the class $\mathscr{A^D}$.
\end{theorem}
\begin{proof}
Due to Lemma~\ref{le:A-duality}, it follows from $\varphi\in\mathscr{A^D}$ that
$\varphi^*\in\mathscr{A^D}$. Since $\varphi^*$ also satisfies the $\Delta_2$-condition, 
we have $\varphi^*\in\mathscr{A^D_\infty}$ by Lemma~\ref{le:A-A-inf}. Therefore, we can apply
Lemma~\ref{le:perturbation} to the function $\varphi^*$, find $\delta\in(0,1)$ and $G>2$ as
in the lemma, and define $r_0:=1+\delta>1$.

We first show that the operators $T_{r_0,\mathcal{Q}}$ are bounded on $L^\varphi(X)$
uniformly over all families $\mathcal{Q}\in\mathfrak{S}$. Fix such a $\mathcal{Q}$, take 
$g\in L^\varphi(X)$ and assume, without loss of generality, that $M_{r_0,Q}g\ne0$ 
for all $Q\in\mathcal{Q}$. Note that this implies $M_Qg\ne0$ for all $Q\in\mathcal{Q}$.
Then for any $h\in L^{\varphi^*}(X)$, self-duality of averaging 
operators~\eqref{eq:self-dual} yields
\begin{align*}
\left\langle T_{r_0,\mathcal{Q}}g,|h|\right\rangle
&=\left\langle\sum_{Q\in\mathcal{Q}}\chi_Q\,M_{r_0,Q}g,|h|\right\rangle
=\left\langle\sum_{Q\in\mathcal{Q}}\chi_Q\,
\frac{(M_{r_0,Q}g)^{r_0}}{(M_{r_0,Q}g)^\delta},|h|\right\rangle \\
&\le\left\langle\sum_{Q\in\mathcal{Q}}\chi_Q\,
\frac{M_Q(|g|^{r_0})}{(M_Qg)^\delta},|h|\right\rangle
=\left\langle\sum_{Q\in\mathcal{Q}}\chi_Q\,M_Q
\left(\frac{|g|^{r_0}}{(M_Qg)^\delta}\right),|h|\right\rangle \\
&=\left\langle T_\mathcal{Q}
\left(\sum_{Q\in\mathcal{Q}}\chi_Q\,\frac{|g|^{r_0}}{(M_Qg)^\delta}\right),|h|\right\rangle
=\left\langle\sum_{Q\in\mathcal{Q}}\chi_Q\,\frac{|g|^{r_0}}{(M_Qg)^\delta},
T_\mathcal{Q}h\right\rangle \\
&=\left\langle\sum_{Q\in\mathcal{Q}}\chi_Q\,\frac{|g|^{1+\delta}}{(M_Qg)^\delta},
\sum_{Q\in\mathcal{Q}}\chi_Q\,M_Qh\right\rangle \\
&=\left\langle\sum_{Q\in\mathcal{Q}}\chi_Q\,|g|,\sum_{Q\in\mathcal{Q}}
\chi_Q\left|\frac{g}{M_Qg}\right|^\delta M_Qh\right\rangle.
\end{align*}
By successive application of Lemma~\ref{le:norm-conj}, 
the above estimate and H\"older's inequality in Lemma~\ref{le:Holder}, 
Lemma~\ref{le:perturbation} for $\varphi^*$, and $\varphi^*\in\mathscr{A^D}$, we obtain
\begin{align*}
\|T_{r_0,\mathcal{Q}}g\|_\varphi
&\le\sup_{\substack{h\in L^{\varphi^*}(X):\\\|h\|_{\varphi^*}\le1}}
\langle T_{r_0,\mathcal{Q}}g,|h|\rangle \\
&\le2
\sup_{\substack{h\in L^{\varphi^*}(X):\\\|h\|_{\varphi^*}\le1}}
\left\|\sum_{Q\in\mathcal{Q}}\chi_Q\,|g|\right\|_\varphi
\left\|\sum_{Q\in\mathcal{Q}}
\chi_Q\left|\frac{g}{M_Qg}\right|^\delta M_Qh \right\|_{\varphi^*} \\
&\le2G
\sup_{\substack{h\in L^{\varphi^*}(X):\\\|h\|_{\varphi^*}\le1}}
\|g\|_\varphi\,\|T_\mathcal{Q}h\|_{\varphi^*} \\
&\le2GC^*\|g\|_\varphi,
\end{align*}
where $C^*$ is the $\mathscr{A^D}$-constant of $\varphi^*$. Since $G$ is 
independent of $\mathcal{Q}$ and $g$, we conclude that the operators $T_{r_0,\mathcal{Q}}$
are bounded on $L^\varphi(X)$ uniformly over all $\mathcal{Q}\in\mathfrak{S}$.

Now define $s_0:=1/r_0$. Then for any $s\in[s_0,1]$, it follows from 
Lemma~\ref{le:comparison}(b) and the uniform boundedness of $T_{r_0,\mathcal{Q}}$, 
$\mathcal{Q}\in\mathfrak{S}$, on $L^\varphi(X)$ that
\begin{align*}
\|T_\mathcal{Q}f\|_\psi
&=\|(T_{1/s,\mathcal{Q}}(|f|^s))^{1/s}\|_\psi
\approx\|T_{1/s,\mathcal{Q}}(|f|^s)\|_\varphi^{1/s} \\
&\le\|T_{r_0,\mathcal{Q}}(|f|^s)\|_\varphi^{1/s}\le c\,\|\,|f|^s\|_\varphi^{1/s}
\approx \|f\|_\psi
\end{align*}
uniformly in $f\in L^\psi(X)$ and $\mathcal{Q}\in\mathfrak{S}$. Hence $\psi\in\mathscr{A^D}$.
\end{proof}
\begin{remark}[Track of constants]\label{re:track}
The \emph{threshold of admissible slowdown $s_0$}---defined for each ``nice enough'' function
$\varphi\in\mathscr{A^D}$ in Theorem~\ref{th:self-improvment}---depends exclusively on:
\begin{enumerate}
    \item[(i)] the $\mathscr{A^D}$-constant of $\varphi$,
    \item[(ii)] the $\Delta_2$-constant of $\varphi^*$,
    \item[(iii)] the parameter $\varepsilon$ of the dyadic system
    $\mathscr{D}$ specified in Theorem~\ref{th:HK}(d).
\end{enumerate}
Indeed, this threshold is set to be $s_0:=1/(1+\delta)$, where $\delta$
originates from Lemma~\ref{le:perturbation} applied to the function
$\varphi^*\in\mathscr{A^D_\infty}$.
The proof of Lemma~\ref{le:perturbation} reveals that $\delta$ depends only on the
$\mathscr{A^D_\infty}$-constant of $\varphi^*$ and the parameter $\varepsilon$ of the 
system $\mathscr{D}$. The $\mathscr{A^D_\infty}$-constant of $\varphi^*$, in turn,
depends only on its $\mathscr{A^D}$-constant and $\Delta_2$-constant, as can be seen by
carefully tracking the constants in Lemma~\ref{le:A-A-inf}. Finally, 
the proof of Lemma~\ref{le:A-duality} shows that the $\mathscr{A^D}$-constant of $\varphi^*$ 
is at most twice the $\mathscr{A^D}$-constant of $\varphi$.
\end{remark}

\subsection{The class $\mathscr{A^D}$ and the domination of mirror means
over standard means}

A ``regular'' subclass of $\mathscr{A^D}$---consisting of proper functions
$\varphi\in\Phi(X)$ that satisfy the $\Delta_2$-condition---can be characterized via the
domination relation between the standard means $M_Q\varphi$ and the mirror means 
$(M_Q\varphi^*)^*$ viewed as generalized $\Phi$-functions on $\mathscr{D}$. The following
theorem is a counterpart of~\cite[Theorems~5.2.15 and 5.2.18]{DHHR11}, and we simplify the
original proof by working with the semimodulars rather than the norms.
\begin{theorem}[Partial characterization of the class $\mathscr{A^D}$]\label{th:1char}
Let $\varphi\in\Phi(X)$ be proper and satisfy the $\Delta_2$-condition. Then 
$\varphi\in\mathscr{A^D}$ if and only if $M_Q\varphi\preceq(M_Q\varphi^*)^*.$
\end{theorem}
\begin{proof}
Notice that $M_Q\varphi,(M_Q\varphi^*)^*\in\Phi(\mathscr{D})$ by Lemma~\ref{le:mean-Phi}
since $\varphi$ is proper. Assume first that $M_Q\varphi\preceq(M_Q\varphi^*)^*$. Take a
family $\mathcal{Q}\in\mathfrak{S}$ and a function $f\in L^\varphi(X)$ with 
$\|f\|_\varphi\le1$, or equivalently $\rho_\varphi(f)\le1$ by the unit ball property
in Lemma~\ref{le:spaces}(e). It follows from Lemma~\ref{le:mirror-modular} that 
\[
\rho_{(M_Q\varphi^*)^*}\big(\{M_Qf\}_{Q\in\mathcal{Q}}\big)
\le \rho_\varphi(f)\le1.
\]
Since $M_Q\varphi$ satisfies the $\Delta_2$-condition together with $\varphi$, 
by Lemma~\ref{le:domination} there exists a constant $A_1>1$, 
independent of $\mathcal{Q}$ and $f$, such that 
\begin{align*}
\rho_\varphi(T_\mathcal{Q}f)
&=\int_X \varphi\bigg(x,\sum_{Q\in\mathcal{Q}}\chi_Q(x)\,M_Qf \bigg)d\mu(x)
=\sum_{Q\in\mathcal{Q}}\int_Q\varphi(x,M_Qf)\,d\mu(x)
\\
&=\sum_{Q\in\mathcal{Q}}\mu(Q)\,(M_Q\varphi)(M_Qf)
=\rho_{M_Q\varphi}\big(\{M_Qf\}_{Q\in\mathcal{Q}}\big)\le A_1.
\end{align*}
The convexity of $\rho_\varphi$ implies $\rho_\varphi(T_\mathcal{Q}f/A_1)\le1$, hence
$\|T_\mathcal{Q}f\|_\varphi\le A_1$ by the unit ball property. A scaling argument gives 
$\|T_\mathcal{Q}f\|_\varphi\le A_1\|f\|_\varphi$ for all  $\mathcal{Q}\in\mathfrak{S}$ and
$f\in L^\varphi(X)$, therefore $\varphi\in\mathscr{A^D}$.

Conversely, let $\varphi\in\mathscr{A^D}$. Fix a family $Q\in\mathfrak{S}$ and
a sequence $\{t_Q\}_{Q\in\mathcal{Q}}\in l^{(M_Q\varphi^*)^*}(\mathcal{Q})$ 
satisfying $\rho_{(M_Q\varphi^*)^*}(\{t_Q\}_{Q\in\mathcal{Q}})\le1$. 
By Lemma~\ref{le:mirror-modular}, there exists a function $f\in L^\varphi(X)$, with 
$M_Qf=|t_Q|$ for all $Q\in\mathcal{Q}$, such that $\rho_\varphi(f)\le2$. The convexity of
$\rho_\varphi$ implies $\rho_\varphi(f/2)\le1$, and hence $\|f\|_\varphi\le2$. If 
$C$ is the $\mathscr{A^D}$-constant of $\varphi$, then there holds
$\|T_\mathcal{Q}f\|_\varphi\le 2C\le2^k$ for some $k\in\mathbb{N}$. Hence the
$\Delta_2$-condition for $\varphi$ with constant $D\ge2$ and the unit ball property yield
\[
\rho_\varphi(T_\mathcal{Q}f)\le D^k\,\rho_\varphi(2^{-k}\,T_\mathcal{Q}f)\le D^k=:A_2.
\]
Note that $A_2>1$ depends only on the $\mathscr{A^D}$-constant and the $\Delta_2$-constant 
of $\varphi$. Thus,
\[
\rho_{M_Q\varphi}\big(\{t_Q\}_{Q\in\mathcal{Q}}\big)
=\rho_{M_Q\varphi}\big(\{M_Qf\}_{Q\in\mathcal{Q}}\big)
=\rho_\varphi(T_\mathcal{Q}f) 
\le A_2
\]
with the constant $A_2$ independent of $\mathcal{Q}\in\mathfrak{S}$ and 
$\{t_Q\}_{Q\in\mathcal{Q}}\in l^{(M_Q\varphi^*)^*}(\mathcal{Q})$. We conclude by
Lemma~\ref{le:domination} that $M_Q\varphi\preceq(M_Q\varphi^*)^*$.
\end{proof}

To derive a refined domination property of $\mathscr{A^D}$ from the 
partial characterization in Theorem~\ref{th:1char}, we need the following auxiliary 
lemma. It shows that the trivial domination $M_Q\varphi\preceq M_{s,Q}\varphi$, 
$s\ge1$, can be reversed in the case of a ``regular'' function $\varphi\in\mathscr{A^D}$, 
yielding $M_{s,Q}\varphi\preceq M_Q\varphi$ for some $s>1$. Our proof of this result
simplifies the argument for the reverse H\"older-type estimate in~\cite[Lemma~5.7.4]{DHHR11} 
and offers a more detailed treatment of~\cite[Remark~5.7.8]{DHHR11}.
\begin{lemma}\label{le:RHI}
Let $\varphi\in\mathscr{A^D}$ be proper and satisfy the $\Delta_2$-condition. Then there 
exists $s>1$ such that $M_{s,Q}\varphi\in\Phi(\mathscr{D})$ and 
$M_{s,Q}\varphi\preceq M_Q\varphi$.
\end{lemma}
\begin{proof}
Notice first that $M_Q\varphi\in\Phi(\mathscr{D})$ by Lemma~\ref{le:mean-Phi} since $\varphi$
is proper. By Lemma~\ref{le:A-A-inf}, we have $\varphi\in\mathscr{A^D_\infty}$ due to the 
$\Delta_2$-condition for $\varphi$. Then we apply Lemma~\ref{le:perturbation} to the
function $\varphi$, find $\delta\in(0,1)$ and $G>2$ as in the lemma, and define
$s:=1+\delta>1$.

Take a family $\mathcal{Q}\in\mathfrak{S}$ and a sequence 
$\{t_Q\}_{Q\in\mathcal{Q}}\in l^{M_Q\varphi}(\mathcal{Q})$, with $t_Q>0$ for all 
$Q\in\mathcal{Q}$, satisfying
\begin{equation}\label{eq:modular1}
\rho_\varphi\bigg(\sum_{Q\in\mathcal{Q}}t_Q\chi_Q\bigg)
=\sum_{Q\in\mathcal{Q}}\mu(Q)\,(M_Q\varphi)(t_Q)
=\rho_{M_Q\varphi}\big(\{t_Q\}_{Q\in\mathcal{Q}}\big)
\le1. 
\end{equation}
Note that $t\chi_Q\in L^\varphi(X)$ for all $Q\in\mathscr{D}$ and $t\ge0$ since $\varphi$ is 
proper; then, due to the $\Delta_2$-condition for $\varphi$, it follows from 
Lemma~\ref{le:delta2}(a)--(b) that $\rho_\varphi(t\chi_Q)\in(0,\infty)$ for all 
$Q\in\mathcal{Q}$ and $t>0$.
This implies that the function $f\in L^1_\text{loc}(X)$ defined by 
\[
f(x):=\sum_{Q\in\mathcal{Q}}\chi_Q(x)\,\varphi(x,t_Q),
\quad x\in X,
\]
is admissible in Lemma~\ref{le:perturbation}, because its means 
\[
M_Qf=(M_Q\varphi)(t_Q)=\frac{\rho_\varphi(t_Q\chi_Q)}{\mu(Q)}
\]
are positive and finite for all $Q\in\mathcal{Q}$. Applying Lemma~\ref{le:perturbation}
with this choice of $f$ and the unit ball property in view of 
inequality~\eqref{eq:modular1}, we obtain that
\[
\bigg\| \sum_{Q\in\mathcal{Q}}t_Q\bigg|\frac{f}{M_Q f}\bigg|^\delta\chi_Q\bigg\|_\varphi
\le
G\,\bigg\|\sum_{Q\in\mathcal{Q}}t_Q\chi_Q\bigg\|_\varphi
\le G\le2^k
\]
for some $k\in\mathbb{N}$, independent of $\mathcal{Q}$ and $\{t_Q\}_{Q\in\mathcal{Q}}$. 
By the $\Delta_2$-condition for $\varphi$ with constant $D\ge2$
and the unit ball property, this leads to
\begin{equation}\label{eq:modular2}
\rho_\varphi\bigg(\sum_{Q\in\mathcal{Q}}t_Q\bigg|\frac{f}{M_Q f}\bigg|^\delta\chi_Q\bigg)
\le D^k\,\rho_\varphi
\bigg(\sum_{Q\in\mathcal{Q}}\frac{t_Q}{2^k}\bigg|\frac{f}{M_Q f}\bigg|^\delta\chi_Q\bigg)
\le D^k.
\end{equation}
Due to the convexity of $\varphi$ in the second variable, we have the simple inequality
\[
c\,\varphi(x,t)\le\left\{\begin{array}{ll}
\varphi(x,t)&\textnormal{if }c\in[0,1]\\
\varphi(x,ct)&\textnormal{if }c>1
\end{array}\right.
\le\varphi(x,t)+\varphi(x,ct)
\]
for all $c\ge0$, $x\in X$ and $t\ge0$. Therefore, using the identity 
$M_Qf=(M_Q\varphi)(t_Q)$ for all $Q\in\mathcal{Q}$, the above inequality, and estimates for 
semimodulars~\eqref{eq:modular1} and \eqref{eq:modular2}, we obtain that
\begin{align}\label{eq:modular3}
\sum_{Q\in\mathcal{Q}}\mu(Q)(M_{s,Q}\varphi)(t_Q)
&\le\sum_{Q\in\mathcal{Q}}\mu(Q)\,\frac{\big((M_{s,Q}\varphi)(t_Q)\big)^{1+\delta}}
{\big((M_Q\varphi)(t_Q)\big)^\delta}
\nonumber \\
&=\sum_{Q\in\mathcal{Q}}\int_Q\varphi(x,t_Q)\,\bigg|\frac{f(x)}{M_Qf}\bigg|^\delta d\mu(x)
\nonumber \\
&\le\sum_{Q\in\mathcal{Q}}\int_Q \left(\varphi(x,t_Q)+
\varphi\bigg(x,t_Q\,\bigg|\frac{f(x)}{M_Qf}\bigg|^\delta\,\bigg)\right)\,d\mu(x)
\nonumber \\
&=\rho_\varphi\bigg(\sum_{Q\in\mathcal{Q}}t_Q\chi_Q\bigg)
+\rho_\varphi\bigg(\sum_{Q\in\mathcal{Q}}t_Q\bigg|\frac{f}{M_Qf}\bigg|^\delta\chi_Q\bigg)
\nonumber \\
&\le1+D^k=:A_1,
\end{align}
where the constant $A_1>1$ is independent of $\mathcal{Q}$ and $\{t_Q\}_{Q\in\mathcal{Q}}$.

It follows from the part already proved that 
$M_{s,Q}\varphi\in\Phi(\mathscr{D})$. Indeed, for any 
$Q\in\mathscr{D}$, the convexity and left-continuity of $t\mapsto(M_{s,Q}\varphi)(t)$
as well as $(M_{s,Q}\varphi)(0)=0$ and the limit $(M_{s,Q}\varphi)(t)\to\infty$ as
$t\to\infty$ follow as in Lemma~\ref{le:mean-Phi} solely due to $\varphi\in\Phi(X)$. For the
limit $(M_{s,Q}\varphi)(t)\to0$ as $t\to0^+$, we use the fact that 
$M_Q\varphi\in\Phi(\mathscr{D})$ and the obtained
implication ``\eqref{eq:modular1}~$\implies$~\eqref{eq:modular3}'' with 
$\mathcal{Q}=\{Q\}$. Since $(M_Q\varphi)(t)\to0$ as $t\to0^+$, we can choose 
$t_Q>0$ satisfying $\mu(Q)\,(M_Q\varphi)(t_Q)\le 1$. Then for all $t\le t_Q$,
we conclude by the convexity of $M_{s,Q}\varphi$ and inequality~\eqref{eq:modular3} that
\[
(M_{s,Q}\varphi)(t)\le \frac{t}{t_Q}\,(M_{s,Q}\varphi)(t_Q)
\le t\,\frac{A_1}{t_Q\,\mu(Q)}\to0
\quad\text{as }t\to0^+.
\]

We can now rewrite implication ``\eqref{eq:modular1}~$\implies$~\eqref{eq:modular3}'' 
in the form
\[
\rho_{M_Q\varphi}\big(\{t_Q\}_{Q\in\mathcal{Q}}\big)\le1
\implies
\rho_{M_{s,Q}\varphi}\big(\{t_Q\}_{Q\in\mathcal{Q}}\big)\le A_1
\]
and note that it holds for all families $\mathcal{Q}\in\mathfrak{S}$ and all sequences
$\{t_Q\}_{Q\in\mathcal{Q}}\in l^{M_Q\varphi}(\mathcal{Q})$, not only with positive
terms, because: (i)~the semimodulars for $\{t_Q\}_{Q\in\mathcal{Q}}$ and 
$\{|t_Q|\}_{Q\in\mathcal{Q}}$ coincide; (ii)~the semimodular sums do not change if we 
exclude the cubes $Q$ for which $t_Q=0$ from the family $\mathcal{Q}$.
Since $M_{s,Q}\varphi$ satisfies the $\Delta_2$-condition together with $\varphi$, we 
conclude by Lemma~\ref{le:domination} that $M_{s,Q}\varphi\preceq M_Q\varphi$.
\end{proof}

We now establish a refined domination result for the class $\mathscr{A^D}$, showing that 
not only does $\varphi\in\mathscr{A^D}$ imply $M_Q\varphi\preceq(M_Q\varphi^*)^*$ 
when $\varphi$ has suitable properties, but also that the stronger relation
$M_{s,Q}\varphi\preceq(M_{s,Q}\varphi^*)^*$ holds for some $s>1$. The next theorem
is based on~\cite[Lemma~5.7.9]{DHHR11}, where we carefully correct a flaw contained 
in the original proof. Up to a change of notation, it was implicitly assumed in
\cite[Lemma~5.7.9]{DHHR11} that the function $M_{s,Q}\varphi^*$ is proper as an element of
$\Phi(\mathcal{Q})$ for every $\mathcal{Q}\in\mathfrak{S}$, since otherwise one could not 
apply~\cite[Remark~5.2.19]{DHHR11}; however, under the lemma's hypotheses one could 
only guarantee that $M_Q\varphi^*\in\Phi(\mathcal{Q})$ is proper for every 
$\mathcal{Q}\in\mathfrak{S}$ in view of~\cite[Lemma~5.2.11]{DHHR11}. We avoid this 
unjustified assumption by using Lemma~\ref{le:domin-conj}, which is a finer version of 
\cite[Remark~5.2.19]{DHHR11}.
\begin{theorem}[Refined domination property of $\mathscr{A^D}$]\label{th:refinement}
Let $\varphi\in\mathscr{A^D}$ be proper, and both $\varphi$ and $\varphi^*$ satisfy the 
$\Delta_2$-condition. Then $M_{s,Q}\varphi\preceq(M_{s,Q}\varphi^*)^*$ for some $s>1$.
\end{theorem}
\begin{proof}
By the partial characterization of the class~$\mathscr{A^D}$ in Theorem~\ref{th:1char}, 
we have that $M_Q\varphi\preceq(M_Q\varphi^*)^*$. By Lemma~\ref{le:A-duality}, the proper
function $\varphi^*$ belongs to $\mathscr{A^D}$, so we can apply Lemma~\ref{le:RHI} to 
$\varphi$ and $\varphi^*$ and find $s>1$ such that 
$M_{s,Q}\varphi,M_{s,Q}\varphi^*\in\Phi(\mathscr{D})$ and there hold dominations
$M_{s,Q}\varphi\preceq M_Q\varphi$ and $M_{s,Q}\varphi^*\preceq M_Q\varphi^*$ 
(we may take $s$ as the minimum between those
numbers $s$ which exist for $\varphi$ and $\varphi^*$ by Lemma~\ref{le:RHI}). 
Notice that $M_Q\varphi^*$ is proper as an element of $\Phi(\mathcal{Q})$ for every
$\mathcal{Q}\in\mathfrak{S}$ due to Lemma~\ref{le:means-proper}. Then 
Lemma~\ref{le:domin-conj} allows us to establish domination
$(M_Q\varphi^*)^*\preceq(M_{s,Q}\varphi^*)^*$ between the conjugate functions,
and we obtain overall that
\[
M_{s,Q}\varphi\preceq M_Q\varphi\preceq(M_Q\varphi^*)^*
\preceq(M_{s,Q}\varphi^*)^*.
\]
This implies $M_{s,Q}\varphi\preceq(M_{s,Q}\varphi^*)^*$ by transitivity of
domination.
\end{proof}
\begin{remark}[Track of constants]
The refined parameter $s>1$ in Theorem~\ref{th:refinement} is determined by 
the constants $\delta$ existing for $\varphi$ and $\varphi^*$ by 
Lemma~\ref{le:perturbation}. Hence $s$ is a function of the 
$\mathscr{A^D}$-constant of $\varphi$, the $\Delta_2$-constants of $\varphi$ and
$\varphi^*$, and the parameter $\varepsilon$ of the dyadic system $\mathscr{D}$ specified 
in Theorem~\ref{th:HK}(d). We refer the reader to Remark~\ref{re:track} for a more detailed
explanation.
\end{remark}

\section{A Sufficient Condition for the Boundedness of $M$}\label{sec:sufficient}
In this section, we define a more restrictive class 
of generalized $\Phi$-functions on $X$ than the class $\mathscr{A^D}$. This 
new class, called $\mathscr{A}^{\mathscr{D}}_\textnormal{strong}$, 
is a dyadic variant of the class
$\mathcal{A}_\textnormal{strong}$ introduced in~\cite[Definition~5.5.6]{DHHR11}. 
Then we prove a sufficient condition for the boundedness of the maximal operator~$M$,
formulated in terms of the class $\mathscr{A}^{\mathscr{D}}_\textnormal{strong}$ 
(Theorem~\ref{th:sufficient}).

\subsection{The class $\mathscr{A}^{\mathscr{D}}_\textnormal{strong}$ as a more
restrictive class than $\mathscr{A^D}$}

In accordance with~\cite[Definition~5.5.6]{DHHR11}, the class 
$\mathscr{A}^{\mathscr{D}}_\textnormal{strong}$ is defined as follows.
\begin{definition}\label{def:A-strong}
We say that $\varphi\in\Phi(X)$ is in the 
\emph{class~$\mathscr{A}^{\mathscr{D}}_\textnormal{strong}$} if 
$M_Q\varphi,(M_Q\varphi^*)^*\in\Phi(\mathscr{D})$ and $M_Q\varphi\ll(M_Q\varphi^*)^*$.
\end{definition}
Let us clarify the idea behind this definition. We know from Theorem~\ref{th:1char} 
that a proper function $\varphi\in\Phi(X)$ satisfying
the $\Delta_2$-condition belongs to the class $\mathscr{A^D}$ if and only if 
$M_Q\varphi\preceq(M_Q\varphi^*)^*$. By Lemma~\ref{le:strong-weak}, 
the assumption $M_Q\varphi\ll(M_Q\varphi^*)^*$ in Definition~\ref{def:A-strong} is 
stronger then $M_Q\varphi\preceq(M_Q\varphi^*)^*$, which leads to the following 
remark:
\begin{remark}\label{re:strong-usual}
If $\varphi\in\Phi(X)$ is proper and satisfies the $\Delta_2$-condition, then
$\varphi\in\mathscr{A}^{\mathscr{D}}_\textnormal{strong}$ implies 
$\varphi\in\mathscr{A^D}$.
\end{remark}
\noindent
Thus, the subset of functions in $\mathscr{A}^{\mathscr{D}}_\textnormal{strong}$ that 
are proper and satisfy the $\Delta_2$-condition forms a narrower class than 
$\mathscr{A^D}$. This is the motivation for the name 
``$\mathscr{A}^{\mathscr{D}}_\textnormal{strong}$.'' It remains an 
open question, however, whether $\varphi\in\mathscr{A^D}$ implies 
$\varphi\in\mathscr{A}^{\mathscr{D}}_\textnormal{strong}$ for a proper function 
$\varphi\in\Phi(X)$ satisfying the $\Delta_2$-condition; in other words, whether
the converse of Remark~\ref{re:strong-usual} holds. We will show in
Theorem~\ref{th:coincidence} that this is true at least for the functions 
$\varphi=\varphi_{p(\cdot)}$ that generate variable Lebesgue spaces $L^{p(\cdot)}(X)$ 
with the exponents satisfying $1<p_-\le p_+<\infty$.

The next lemma shows, in essence, that if we take a ``regular'' function 
$\varphi\in\mathscr{A}^{\mathscr{D}}_\textnormal{strong}$ (proper and satisfying the 
$\Delta_2$-condition) and ``speed it up'' in $t$---by replacing $t$ with a larger 
$t^s$---the resulting function $(x,t)\mapsto\varphi(x,t^s)$ still belongs to the class
$\mathscr{A}^{\mathscr{D}}_\textnormal{strong}$. This is a variant 
of~\cite[Lemma~5.5.9]{DHHR11} with the added requirement that $\varphi$ satisfy 
the $\Delta_2$-condition, which allows us to apply Lemma~\ref{le:mirror-mean}; see
the discussion preceding Lemma~\ref{le:mirror-mean}. We also take this opportunity to 
correct a misprint in the formulation of~\cite[Lemma~5.5.9]{DHHR11}: in its conclusion,
the roles of $\varphi$ and $\psi$ are mistakenly swapped.
\begin{lemma}\label{le:strong-strong}
Let $\varphi,\psi\in\Phi(X)$ be proper with
$\psi(x,t)\approx\varphi(x,t^s)$, uniformly in $x\in X$ and $t\ge0$, for some $s\ge1$. 
Suppose that $\varphi$ satisfies the $\Delta_2$-condition. If 
$\varphi\in\mathscr{A}^{\mathscr{D}}_\textnormal{strong}$, then
$\psi\in\mathscr{A}^{\mathscr{D}}_\textnormal{strong}$.
\end{lemma}
\begin{proof}
By assumption, there exists $c>1$ such that, for all $x\in X$ and $t\ge0$,
\begin{equation}\label{eq:connection}
\frac1c\,\varphi(x,t^s)\le\psi(x,t)\le c\,\varphi(x,t^s).
\end{equation}
Due to this relation, the $\Delta_2$-condition for $\varphi$ implies the 
$\Delta_2$-condition for $\psi$. The second inequality in~\eqref{eq:connection} 
immediately gives
\begin{equation}\label{eq:mean-1}
(M_Q\psi)(t)\le c\,(M_Q\varphi)(t^s)    
\end{equation}
for all $Q\in\mathscr{D}$ and $t\ge0$. The first inequality in~\eqref{eq:connection}, 
Lemma~\ref{le:mirror-mean} applied to $\varphi$ and $\psi$, and the monotonicity of 
means due to Jensen's inequality yield
\begin{align}
(M_Q\varphi^*)^*(t^s)&=\inf_{\substack{f\in L^0(X):\\M_Qf\ge t^s}}M_Q(\varphi(f))
\le c\inf_{\substack{f\in L^0(X):\\M_Qf\ge t^s}}M_Q(\psi(|f|^{1/s}))
\nonumber \\
&=c\inf_{\substack{g\in L^0(X):\\M_{s,Q}g\ge t}}M_Q(\psi(g))
\le c\inf_{\substack{g\in L^0(X):\\M_Qg\ge t}}M_Q(\psi(g))
=c\,(M_Q\psi^*)^*(t)
\label{eq:mean-2}
\end{align}
for all $Q\in\mathscr{D}$ and $t\ge0$. To prove that 
$\psi\in\mathscr{A}^{\mathscr{D}}_\textnormal{strong}$
if $\varphi\in\mathscr{A}^{\mathscr{D}}_\textnormal{strong}$, fix a number $B_1>0$ and 
a sequence of families $\mathcal{Q}_k\in\mathfrak{S}$, $k\in\mathbb{Z}$, such that 
\[
\sum_{k=-\infty}^\infty\sum_{Q\in\mathcal{Q}_k}\mu(Q)\,(M_Q\psi^*)^*(2^k)\le B_1.
\]
In view of~\eqref{eq:mean-2}, this implies 
\[
\sum_{k=-\infty}^\infty\sum_{Q\in\mathcal{Q}_k}\mu(Q)\,(M_Q\varphi^*)^*((2^s)^k)\le c\,B_1
=:C_1.
\]
Then $\varphi\in\mathscr{A}^{\mathscr{D}}_\textnormal{strong}$ and
Lemma~\ref{le:alpha^k} with $\alpha=2^s$ ensure that for this $C_1>0$, there is a
constant $C_2>0$ independent of $\{\mathcal{Q}_k\}_{k\in\mathbb{Z}}$ such that 
inequality~\eqref{eq:mean-1} leads to
\[
\sum_{k=-\infty}^\infty\sum_{Q\in\mathcal{Q}_k}\mu(Q)\,(M_Q\psi)(2^k)
\le 
c\sum_{k=-\infty}^\infty\sum_{Q\in\mathcal{Q}_k}\mu(Q)\,(M_Q\varphi)((2^s)^k)
\le c\,C_2=:B_2.
\]
By Definition~\ref{def:strong-domin} we conclude that $M_Q\psi\ll(M_Q\psi^*)^*$ and thus
$\psi\in\mathscr{A}^{\mathscr{D}}_\textnormal{strong}$.
\end{proof}

\subsection{Boundedness of the maximal operator $M$ on $L^\varphi(X)$}

To prove the main theorem of this section---Theorem~\ref{th:sufficient} with a 
sufficient condition for the boundedness of $M$---we need the following auxiliary lemma.
Our Lemma~\ref{le:aux-to-suff} serves as an analogue of~\cite[Lemma~5.5.7]{DHHR11} for the 
dyadic maximal function on spaces of homogeneous type, but it is obtained via a
\emph{completely different proof}. The proof of~\cite[Lemma~5.5.7]{DHHR11} in $\mathbb{R}^n$
relies fundamentally on the Besicovitch covering theorem, a tool specific to the Euclidean
setting. We circumvent the need for the Besicovitch theorem by using the finite collection
of adjacent dyadic grids $\mathscr{D}^t$ and the Calder\'on--Zygmund-type decomposition of
the superlevel sets of $t$-dyadic maximal functions from Lemma~\ref{le:CZ}. 
\begin{lemma}\label{le:aux-to-suff}
Let $\varphi\in\mathscr{A^D}$ be such that $\chi_F\notin L^\varphi(X)$ for all 
$F\subseteq X$ with $\mu(F)=\infty$. 
For every $f\in L^\varphi(X)$ and $\lambda>0$, there
exist $K$ families $\mathcal{Q}_1,\ldots,\mathcal{Q}_K\in\mathfrak{S}$, where $K$ 
is the number of adjacent dyadic grids in $\mathscr{D}=\bigcup_{t=1}^K\mathscr{D}^t$, such that 
\[
M_Q f>\lambda\quad\text{for all } Q\in\mathcal{Q}_t,\; t\in\{1,\ldots, K\}, 
\]
and
\[
\int_X\varphi(x,\lambda)\,\chi_{\{M^\mathscr{D}f>\lambda\}}(x)\,d\mu(x)
\le\sum_{t=1}^K \sum_{Q\in\mathcal{Q}_t}\mu(Q)\,(M_Q\varphi)(\lambda).
\]
\end{lemma}
\begin{proof}
Let $f\in L^\varphi(X)$ and $\lambda>0$. Note first that 
\begin{equation}\label{eq:union}
\left\{x\in X\,:\,M^\mathscr{D}f(x)>\lambda\right\}=
\bigcup_{t=1}^K \left\{x\in X\,:\,M^{\mathscr{D}^t}f(x)>\lambda\right\}.
\end{equation}
Indeed, if $M^\mathscr{D}f(x)>\lambda$, then there exists a cube $Q\in\mathscr{D}$ 
such that $Q\ni x$ and $M_Qf>\lambda$. At the same time, $Q\in\mathscr{D}$ implies that 
$Q\in\mathscr{D}^t$ for some $t\in\{1,\ldots,K\}$. Therefore
$M^{\mathscr{D}^t}f(x)\ge M_Qf>\lambda$, and thus
$x$ belongs to the union on the right-hand side of~\eqref{eq:union}. Conversely, if 
$M^{\mathscr{D}^t}f(x)>\lambda$ for some $t\in\{1,\ldots,K\}$, there is a cube 
$Q'\in\mathscr{D}^t\subseteq\mathscr{D}$ such that
$Q'\ni x$ and $M_{Q'}f>\lambda$. Thus $M^\mathscr{D}f(x)\ge M_{Q'}f>\lambda$.

As in Lemma~\ref{le:cones}, denote by $Q_k^t(x)\in\mathscr{D}^t_k$ the unique cube from 
the $k$-th generation of cubes in $\mathscr{D}^t$ that contains a point $x\in X$.
Assumption $\varphi\in\mathscr{A^D}$ implies, in particular, 
that $T_{\{Q^t_k(x)\}}f=M_{Q^t_k(x)}f\,\chi_{Q^t_k(x)}\in L^\varphi(X)$ 
and 
\begin{equation}\label{eq:means-to-0}
M_{Q^t_k(x)}f
=\frac{\|T_{\{Q^t_k(x)\}}f\|_\varphi}{\|\chi_{Q^t_k(x)}\|_\varphi}
\le C\,\frac{\|f\|_\varphi}{\|\chi_{Q^t_k(x)}\|_\varphi}<\infty
\end{equation}
for all $x\in X$, $t\in\{1,\ldots,K\}$ and $k\in\mathbb{Z}$, where $C$ is the 
$\mathscr{A^D}$-constant of $\varphi$. This has two consequences:
\begin{itemize}
    \item[(i)] $f\in L^1_\text{loc}(X)$. Indeed, take a ball $B:=B(x,r)\subseteq X$.
    By Theorem~\ref{th:HK}, there exists a containing cube $Q^t_k(x)\supseteq B$, with some 
    $t\in\{1,\ldots,K\}$ and $k\in\mathbb{Z}$. Then 
    \[
    \int_B|f(y)|\,d\mu(y)\le \mu(Q^t_k(x))\,M_{Q^t_k(x)}f<\infty.
    \]
    \item[(ii)] For all $x\in X$ and $t\in\{1,\ldots,K\}$, the means $M_{Q^t_k(x)}f\to0$ as 
    $k\to-\infty$. To check this, observe that by Lemma~\ref{le:cones}(a), the set  
    $F:=\bigcup_{k\in\mathbb{Z}}Q_k^t(x)$ has infinite measure. Consequently, 
    $\chi_{Q^t_k(x)}\nearrow\chi_F\notin L^\varphi(X)$
    as $k\to-\infty$. Therefore $\|\chi_{Q^t_k(x)}\|_\varphi\to\infty$ as $k\to-\infty$, 
    since otherwise, if $\sup_{k\in\mathbb{Z}}\|\chi_{Q^t_k(x)}\|_\varphi<\infty$, 
    the Fatou property in Lemma~\ref{le:spaces}(c) would imply that $\chi_F\in L^\varphi(X)$.
    Hence, we conclude from estimate~\eqref{eq:means-to-0} that 
    $M_{Q^t_k(x)}f\to0$ as $k\to-\infty$.
\end{itemize}
In view of Remark~\ref{re:CZ}, this allows us to apply the Calder\'on--Zygmund-type
Lemma~\ref{le:CZ} to the function $f$ and a fixed number $\lambda>0$. By this lemma, 
for each $t\in\{1,\ldots,K\}$ there exists a countable family 
$\mathcal{Q}_t\in\mathfrak{S}$ of cubes from $\mathscr{D}^t$ such that $M_Qf>\lambda$
for all $Q\in\mathcal{Q}_t$ and 
\[
\left\{x\in X\,:\,M^{\mathscr{D}^t}f(x)>\lambda\right\}=\bigcup_{Q\in\mathcal{Q}_t}Q.
\]
This and~\eqref{eq:union} imply, for any $x\in X$, that
\begin{equation}\label{eq:chi-sum}
\chi_{\{M^\mathscr{D}f>\lambda\}}(x)\le\sum_{t=1}^K \chi_{\{M^{\mathscr{D}^t}f>\lambda\}}(x)
=\sum_{t=1}^K\sum_{Q\in\mathcal{Q}_t}\chi_Q(x).
\end{equation}
Thus
\begin{align*}
\int_X\varphi(x,\lambda)\,\chi_{\{M^\mathscr{D}f>\lambda\}}(x)\,d\mu(x)
&\le\int_X\varphi(x,\lambda)\sum_{t=1}^K\sum_{Q\in\mathcal{Q}_t}\chi_Q(x)\,d\mu(x) \\
&=\sum_{t=1}^K \sum_{Q\in\mathcal{Q}_t}\mu(Q)\,(M_Q\varphi)(\lambda),
\end{align*}
which concludes the proof of the lemma.
\end{proof}
The construction in the proof of Lemma~\ref{le:aux-to-suff} yields  
the following natural result.
\begin{lemma}\label{le:aefinite}
Let $\varphi\in\mathscr{A^D}$ be such that $\chi_F\notin L^\varphi(X)$ for all 
$F\subseteq X$ with $\mu(F)=\infty$. Then for any $f\in L^\varphi(X)$, the 
function $M^\mathscr{D}f$ is finite $\mu$-almost everywhere on $X$.
\end{lemma}
\begin{proof}
Just as in the proof of Lemma~\ref{le:aux-to-suff}, for every $n\in\mathbb{N}$, 
there exist $K$ families $\mathcal{Q}_1,\ldots,\mathcal{Q}_K\in\mathfrak{S}$ such that 
$M_Qf>n$ for all $Q\in\mathcal{Q}_t$, $t\in\{1,\ldots,K\}$, and 
inequality~\eqref{eq:chi-sum} holds. Since $\varphi\in\mathscr{A^D}$, we consequently have
\begin{align*}
\|\chi_{\{M^\mathscr{D}f>n\}}\|_\varphi
&\le\frac1n\sum_{t=1}^K \bigg\|\sum_{Q\in\mathcal{Q}_t}n\chi_Q\bigg\|_\varphi
\le\frac1n\sum_{t=1}^K \bigg\|\sum_{Q\in\mathcal{Q}_t}\chi_Q\,M_Qf\bigg\|_\varphi
\\
&=\frac1n\sum_{t=1}^K \|T_{\mathcal{Q}_t}f\|_\varphi\le \frac{CK}{n}\|f\|_\varphi,
\end{align*}
where $C$ is the $\mathscr{A^D}$-constant of $\varphi$. This implies that 
\[
\sup_{n\in\mathbb{N}}\|\chi_{\{M^\mathscr{D}f>n\}}\|_\varphi\le
CK\,\|f\|_\varphi<\infty,
\]
from which it follows by the Fatou property in Lemma~\ref{le:spaces}(c) that 
\[
\|\chi_{\{M^\mathscr{D}f=\infty\}}\|_\varphi
=\lim_{n\to\infty}\|\chi_{\{M^\mathscr{D}f>n\}}\|_\varphi
\le \lim_{n\to\infty}\frac{CK}{n}\|f\|_\varphi=0. 
\]
Thus $M^\mathscr{D}f(x)<\infty$ for $\mu$-almost every $x\in X$.
\end{proof}
We are now ready to present and prove a sufficient condition for the boundedness of 
the maximal operator $M$ on the space $L^\varphi(X)$. The following theorem, which is a 
counterpart of~\cite[Theorem~5.5.12]{DHHR11}, asserts that if a ``nice enough'' function 
$\varphi\in\Phi(X)$ is such that its slowed-down version $(x,t)\mapsto\varphi(x,t^{s_0})$ 
for some $s_0\in(0,1)$ belongs to the class $\mathscr{A}^{\mathscr{D}}_\textnormal{strong}$,
then $M$ is bounded on $L^\varphi(X)$. 

Compared to the original result in~\cite{DHHR11},
we simplify the assumptions by formulating our theorem for $\Phi$-functions instead of 
$N$-functions and by imposing the $\Delta_2$-condition only on the function $\varphi$.
Alongside this, we streamline the original proof while retaining its overall structure. The 
key difference of our proof, however, lies in the use of the newly obtained 
Lemma~\ref{le:aux-to-suff} in place of~\cite[Lemma~5.5.7]{DHHR11}. 

\begin{theorem}\label{th:sufficient}
Let $\varphi,\psi_0\in \Phi(X)$ be proper with $\psi_0(x,t)\approx\varphi(x,t^{s_0})$, 
uniformly in $x\in X$ and $t\ge0$,
for some $s_0\in(0,1)$. Suppose that $\varphi$ satisfies the $\Delta_2$-condition and 
$\chi_F\notin L^\varphi(X)$ for all $F\subseteq X$ with $\mu(F)=\infty$. If
$\psi_0\in\mathscr{A}^{\mathscr{D}}_\textnormal{strong}$, then
$M$ is bounded on $L^\varphi(X)$.
\end{theorem}
\begin{proof}
Due to the pointwise equivalence between the classical and the dyadic maximal
functions in Lemma~\ref{le:M-and-dyadicM}, the maximal operator $M$ is bounded on
$L^\varphi(X)$ if and only if the dyadic maximal operator $M^\mathscr{D}$ is 
bounded on $L^\varphi(X)$. Hence it is enough to prove the claim
for $M^\mathscr{D}$ instead of $M$. In fact, it suffices to show that there exists 
$C\ge1$ such that for all $f\in L^\varphi(X)$,
\begin{equation}\label{eq:suff-goal}
\int_X\varphi(x,|f(x)|)\,d\mu(x)\le1
\implies
\int_X\varphi(x,M^\mathscr{D}f(x))\,d\mu(x)\le C.
\end{equation}
Indeed, let $\|f\|_\varphi\le1$; then $\rho_\varphi(f)\le1$ by the unit ball
property in Lemma~\ref{le:spaces}(e). When~\eqref{eq:suff-goal} is proved, this will imply 
$\rho_\varphi(\frac{M^\mathscr{D}f}{C})\le\frac1C\,\rho_\varphi(M^\mathscr{D}f)\le1$ 
due to the convexity of~$\rho_\varphi$. Then it follows by the unit ball property that
$\|M^\mathscr{D}f\|_\varphi\le C$,
and a scaling argument yields $\|M^\mathscr{D}f\|_\varphi\le C\,\|f\|_\varphi$ for all 
$f\in L^\varphi(X)$.

\medskip\noindent
\emph{Step 1: Estimating the integral $\int_X\varphi(\cdot,M^\mathscr{D}f)\,d\mu$ by a sum.}
Observe first that the pointwise relation between $\varphi$ and $\psi_0$ guarantees that 
$\psi_0$ satisfies the $\Delta_2$-condition together with $\varphi$ and 
there holds $\varphi(x,t)\approx\psi_0(x,t^{1/s_0})$, 
uniformly in $x\in X$ and $t\ge0$, with $1/s_0>1$. Since also 
$\psi_0\in\mathscr{A}^{\mathscr{D}}_\textnormal{strong}$, it follows by 
Lemma~\ref{le:strong-strong} that $\varphi\in\mathscr{A}^{\mathscr{D}}_\textnormal{strong}$.
The function $\varphi$ is proper and satisfies the $\Delta_2$-condition, 
so it follows by Remark~\ref{re:strong-usual} that $\varphi\in\mathscr{A^D}$. 
This and the assumption that $\chi_F\notin L^\varphi(X)$ for all $F\subseteq X$
with $\mu(F)=\infty$ yield, by Lemma~\ref{le:aefinite}, that $M^\mathscr{D}f$ is 
finite $\mu$-almost everywhere on $X$.

Take an $f\in L^\varphi(X)$ such that $\int_X\varphi(x,|f(x)|)\,d\mu(x)\le1$.
We start with the estimate
\begin{equation}\label{eq:suff-1}
 \varphi(x,M^\mathscr{D}f(x))\le c\sum_{k=-\infty}^\infty
 \varphi(x,2^{k})\,\chi_{\{M^\mathscr{D}f>2^{k+1}\}}(x)
\end{equation}
for all $x\in X$, where $c=c(D)$ depends only on the $\Delta_2$-constant $D\ge2$
of~$\varphi$. 
To check~\eqref{eq:suff-1}, fix an $x\in X$ and find the number $m\in\mathbb{Z}$
satisfying $2^{m-1}<M^\mathscr{D}f(x)\le2^{m}$; otherwise, if $M^\mathscr{D}f(x)=0$, 
inequality~\eqref{eq:suff-1} obviously holds. Then the $\Delta_2$-condition 
for $\varphi$ implies
\begin{align*}
\sum_{k=-\infty}^\infty\varphi(x,2^{k})\,\chi_{\{M^\mathscr{D}f>2^{k+1}\}}(x)
&=\sum_{k=-\infty}^{m-2}\varphi(x,2^{k})
=\sum_{k=-\infty}^{m-2}\varphi\left(x,\frac{2^{m}}{2^{m-k}}\right)
\\
&\ge\varphi(x,2^{m})\sum_{k=-\infty}^{m-2}\frac{1}{D^{m-k}}
\ge\frac{\varphi(x,M^\mathscr{D}f(x))}{D(D-1)},
\end{align*}
which is exactly~\eqref{eq:suff-1} with $c=D(D-1)$. As a consequence
of~\eqref{eq:suff-1}, 
\[
\int_X\varphi(x,M^\mathscr{D}f(x))\,d\mu(x)\le
c\sum_{k=-\infty}^\infty\int_X\varphi(x,2^k)\,
\chi_{\{M^\mathscr{D}f>2^{k+1}\}}(x)\,d\mu(x).
\]

For each $k\in\mathbb{Z}$, let us decompose the function $f$ as
$f=f_k+\tilde{f}_{k}$, where
\[
f_k:=f\chi_{\{|f|>2^{k}\}}
\quad\text{and}\quad
\tilde{f_k}:=f\chi_{\{|f|\le2^{k}\}}.
\]
Consequently, for any $k\in\mathbb{Z}$, we have
\begin{equation}\label{eq:suff-2}
\{x\in X\,:\,M^\mathscr{D}f(x)>2^{k+1}\}\subseteq
\{x\in X\,:\,M^\mathscr{D}f_k(x)>2^k\}.
\end{equation}
Indeed, if $M^\mathscr{D}f(x)>2^{k+1}$, we can find 
a cube $Q\in\mathscr{D}$ such that $Q\ni x$ and $M_Qf>2^{k+1}$. Then 
$2^{k+1}<M_Qf\le M_Qf_k+M_Q\tilde{f}_k$, hence at least one of the inequalities 
$M_Qf_k>2^k$ or $M_Q\tilde{f}_k>2^k$ is true. However, by definition $|\tilde{f}_k|\le2^k$
on $X$, and thus $M_Q\tilde{f}_k\le2^k$. Therefore $M_Qf_k>2^k$, which implies
$M^\mathscr{D}f_k(x)>2^k$.
Further, inclusion~\eqref{eq:suff-2} gives
\[
\int_X\varphi(x,M^\mathscr{D}f(x))\,d\mu(x)
\le c\sum_{k=-\infty}^\infty\int_X\varphi(x,2^k)\,
\chi_{\{M^\mathscr{D}f_k>2^k\}}(x)\,d\mu(x).
\]

The assumptions of the theorem, together with the fact that $\varphi\in\mathscr{A^D}$ 
established above, allow us to apply 
Lemma~\ref{le:aux-to-suff} to the functions $f_k\in L^\varphi(X)$ and numbers $2^k>0$ 
for every $k\in\mathbb{Z}$, yielding families 
$\mathcal{Q}_{k,1},\ldots,\mathcal{Q}_{k,K}\in\mathfrak{S}$ such that
\[
M_Qf_k>2^k
\quad\text{for all }Q\in\mathcal{Q}_{k,t},\; t=1,\ldots,K,
\]
and
\[
\int_X\varphi(x,2^k)\,\chi_{\{M^\mathscr{D}f_k>2^k\}}(x)\,d\mu(x)
\le\sum_{t=1}^K\sum_{Q\in\mathcal{Q}_{k,t}}\mu(Q)\,(M_Q\varphi)(2^k).
\]
This and $\psi_0(x,t)\approx\varphi(x,t^{s_0})$, uniformly in $x\in X$ and 
$t\ge0$, together lead to
\begin{align}\label{eq:suff-3}
\int_X\varphi(x,M^\mathscr{D}f(x))\,d\mu(x)
&\le c\sum_{t=1}^K
\sum_{k=-\infty}^\infty\sum_{Q\in\mathcal{Q}_{k,t}}
\mu(Q)\,(M_Q\varphi)(2^k)
\nonumber \\
&\le c\sum_{t=1}^K
\left(\sum_{k=-\infty}^\infty\sum_{Q\in\mathcal{Q}_{k,t}}
\mu(Q)\,(M_Q\psi_0)(2^{k/s_0})\right)
\end{align}
with some constant $c>1$, renewed in the second line, independent of $f$.

\medskip\noindent
\emph{Step 2: The use of strong domination.}
Since $\psi_0\in\mathscr{A}^{\mathscr{D}}_\textnormal{strong}$, we have 
$M_Q\psi_0\ll(M_Q\psi_0^*)^*$. Note that $(M_Q\psi_0^*)^*$ satisfies the 
$\Delta_2$-condition together with $\psi_0$ in view of Lemma~\ref{le:delta2-means}.
We will show that there is a constant $B_1>0$ such that
for each $t\in\{1,\ldots,K\}$,
\begin{equation}\label{eq:suff-4}
S_t:=\sum_{k=-\infty}^\infty\sum_{Q\in\mathcal{Q}_{k,t}}
\mu(Q)\,(M_Q\psi^*_0)^*(2^{k/s_0})\le B_1.
\end{equation}
Once~\eqref{eq:suff-4} is proved, the strong domination $M_Q\psi_0\ll(M_Q\psi_0^*)^*$
and Lemma~\ref{le:alpha^k} with $\alpha=2^{1/s_0}$ imply the existence of a constant 
$B_2>0$, independent of $t$, which bounds the sum in parentheses 
in~\eqref{eq:suff-3}. Consequently, inequality~\eqref{eq:suff-3} gives
\[
\int_X\varphi(x,M^\mathscr{D}f(x))\,d\mu(x)\le cKB_2=:C
\]
for any $f\in L^\varphi(X)$ satisfying $\int_X\varphi(x,|f(x)|)\,d\mu(x)\le1$. Thus
claim~\eqref{eq:suff-goal}, as well as the entire 
theorem, is proved. It only remains to show~\eqref{eq:suff-4}.

\medskip\noindent
\emph{Step 3: Proof of the remaining upper bound for $S_t$.} 
For every $k\in\mathbb{Z}$, we trivially have
\[
0\le|f|^{1-1/s_0}\,2^{k(1/s_0-1)}\,\chi_{\{|f|>2^k\}}\le1
\]
on the set $\{x\in X:f(x)\ne0\}$ since $s_0\in(0,1)$. This estimate, 
the convexity of $\psi_0$ in the second variable, the relation 
$\psi_0(x,t)\approx\varphi(x,t^{s_0})$ uniformly in $x\in X$ and $t\ge0$,
and the definition of $f_k$ imply, for all
$k\in\mathbb{Z}$, that
\begin{align*}
\psi_0(|f_k|\,2^{k(1/s_0-1)})&=\psi_0(|f|\,2^{k(1/s_0-1)}\,\chi_{\{|f|>2^k\}})
\\
&\le\psi_0(|f|^{1/s_0})\,
|f|^{1-1/s_0}\,2^{k(1/s_0-1)}\,\chi_{\{|f|>2^k\}}
\\
&\le c\,\varphi(|f|)\,|f|^{1-1/s_0}\,2^{k(1/s_0-1)}\,\chi_{\{|f|>2^k\}}
\end{align*}
on $\{x\in X:f(x)\ne0\}$; here, as before, $\varphi(f)$ denotes 
$\varphi(\cdot,|f(\cdot)|)$. Using the fact that $M_Qf_k>2^k$ for all $k\in\mathbb{Z}$ and 
$Q\in\mathcal{Q}_{k,t}$, inequality~\eqref{ineq:for-fav-lemma} in 
Lemma~\ref{le:mirror-inequality} and 
the above estimate, we obtain
\begin{align*}
S_t&\le\sum_{k=-\infty}^\infty\sum_{Q\in\mathcal{Q}_{k,t}}
\mu(Q)\,(M_Q\psi^*_0)^*(2^{k(1/s_0-1)}\,M_Qf_k)
\\
&\le\sum_{k=-\infty}^\infty\sum_{Q\in\mathcal{Q}_{k,t}}
\mu(Q)\,M_Q(\psi_0(|f_k|\,2^{k(1/s_0-1)}))
\\
&\le\sum_{k=-\infty}^\infty\int_{\{f\ne0\}}\psi_0(x,|f_k(x)|\,2^{k(1/s_0-1)})\,d\mu(x)
\\
&\le c\int_{\{f\ne0\}}\varphi(x,|f(x)|)\,|f(x)|^{1-1/s_0}
\bigg(
\sum_{k=-\infty}^\infty2^{k(1/s_0-1)}\,\chi_{\{|f|>2^k\}}(x)
\bigg)d\mu(x).
\end{align*}
For every $x\in X$ such that $f(x)\ne0$, we can find the number $m\in\mathbb{N}$ satisfying
$2^m<|f(x)|\le2^{m+1}$ and bound the sum in the parentheses by 
\begin{align*}
\sum_{k=-\infty}^\infty2^{k(1/s_0-1)}\,\chi_{\{|f|>2^k\}}(x)=
\sum_{k=-\infty}^m 2^{k(1/s_0-1)}=\frac{2^{m(1/s_0-1)}}{1-2^{1-1/s_0}}
\le\frac{|f(x)|^{1/s_0-1}}{1-2^{1-1/s_0}}.
\end{align*}
In combination with the previous 
estimate this implies, due to our choice of $f\in L^\varphi(X)$ with 
$\int_X\varphi(x,|f(x)|)\,d\mu(x)\le1$, that
\[
S_t\le \frac{c}{1-2^{1-1/s_0}}\int_{\{f\ne0\}}\varphi(x,|f(x)|)\,d\mu(x)
\le\frac{c}{1-2^{1-1/s_0}}=:B_1,
\]
where $B_1>0$ is independent of $t\in\{1,\ldots,K\}$. The theorem is thus proved. 
\end{proof}

\section{``Variable-Lebesgue'' Functions $\varphi_{p(\cdot)}$ in the Class $\mathscr{A^D}$}
\label{sec:V-L}
In this section, we move from the general framework of Musielak--Orlicz spaces to 
the special case of variable Lebesgue spaces $L^{p(\cdot)}(X):=L^{p(\cdot)}(X,d,\mu)$,
with exponents $p(\cdot)\in\mathcal{P}(X):=\mathcal{P}(X,d,\mu)$
satisfying $1<p_-\le p_+<\infty$, over an unbounded space of homogeneous type 
$(X,d,\mu)$ equipped with a Borel-semiregular measure $\mu$. As before, a dyadic system
$\mathscr{D}$ provided by Theorem~\ref{th:HK} is associated with $X$. 

Recall from Subsection~\ref{subseq:VLS} that $L^{p(\cdot)}(X)$ is a 
Musielak--Orlicz space generated by either of the generalized $\Phi$-functions 
$\varphi_{p(\cdot)}:=\bar{\varphi}_{p(\cdot)}$ or 
$\varphi_{p(\cdot)}:=\widetilde{\varphi}_{p(\cdot)}$ on $X$. 
A remarkable feature of the ``variable-Lebesgue'' functions $\varphi_{p(\cdot)}$
is that the condition $\varphi_{p(\cdot)}\in\mathscr{A^D}$ implies 
$\varphi_{p(\cdot)}\in\mathscr{A}^{\mathscr{D}}_\textnormal{strong}$ whenever 
$1<p_-\le p_+<\infty$---even though the class $\mathscr{A}^{\mathscr{D}}_\textnormal{strong}$
is expected to be more restrictive than the class $\mathscr{A^D}$. This fact is a dyadic 
analogue of~\cite[Theorem~5.7.1]{DHHR11} and the main result 
of this section (Theorem~\ref{th:coincidence}).
The abbreviated notations $p(\cdot)\in\mathscr{A^D}$ for
$\varphi_{p(\cdot)}\in\mathscr{A^D}$ and 
$p(\cdot)\in\mathscr{A}^{\mathscr{D}}_\textnormal{strong}$ for 
$\varphi_{p(\cdot)}\in\mathscr{A}^{\mathscr{D}}_\textnormal{strong}$ 
are adopted in what follows.

\subsection{Five auxiliary lemmas}

The proof of Theorem~\ref{th:coincidence} is intricate and relies on several auxiliary 
lemmas, which we present in this subsection. With the exception of the subsidiary
Lemma~\ref{le:gammas}, the following lemmas establish necessary conditions for 
the domination $M_{s,Q}\varphi_{p(\cdot)}\preceq(M_{s,Q}\varphi^*_{p(\cdot)})^*$ 
with $s\ge1$. Note, incidentally, that the standard means $M_{s,Q}\varphi_{p(\cdot)}$ and
the mirror means $(M_{s,Q}\varphi^*_{p(\cdot)})^*$ are generalized $N$-functions 
when $1<p_-\le p_+<\infty$.
\begin{remark}\label{re:N-means}
For any exponent $p(\cdot)\in\mathcal{P}(X)$ with $1<p_-\le p_+<\infty$ and any parameter $s\ge1$, 
we have $M_{s,Q}\varphi_{p(\cdot)},(M_{s,Q}\varphi^*_{p(\cdot)})^*\in N(\mathscr{D})$. 
Indeed,
\begin{enumerate}
    \item[(i)] $\varphi_{p(\cdot)},\varphi^*_{p(\cdot)}\in N(X)$ in view of 
Lemma~\ref{le:phi-bar-tilde};
    \item[(ii)] again by Lemma~\ref{le:phi-bar-tilde}, uniformly in $x\in X$ and $t\ge0$ there 
hold
\begin{align*}
(\varphi_{p(\cdot)}(x,t))^s&\approx(\bar{\varphi}_{p(\cdot)}(x,t))^s
=\bar{\varphi}_{sp(\cdot)}(x,t),
\\
(\varphi^*_{p(\cdot)}(x,t))^s&\approx(\bar{\varphi}_{p'(\cdot)}(x,t))^s
=\bar{\varphi}_{sp'(\cdot)}(x,t);
\end{align*}
    \item[(iii)] then, since $\bar{\varphi}_{sp(\cdot)}$ and
$\bar{\varphi}_{sp'(\cdot)}$ are proper, so are $(\varphi_{p(\cdot)})^s$ and 
$(\varphi^*_{p(\cdot)})^s$ by Lemma~\ref{le:approper};
    \item[(iv)] both $(\varphi_{p(\cdot)})^s$ and 
$(\varphi^*_{p(\cdot)})^s$ satisfy the $\Delta_2$-condition, since 
$\varphi_{p(\cdot)}$ and $\varphi^*_{p(\cdot)}$ do so by Lemma~\ref{le:phi-bar-tilde};
    \item[(v)] the fact that
$M_{s,Q}\varphi_{p(\cdot)},(M_{s,Q}\varphi^*_{p(\cdot)})^*\in N(\mathscr{D})$
follows by Lemma~\ref{le:mean-N}.
\end{enumerate}
\end{remark}

Our first Lemma~\ref{le:means-of-char} is not specific to the 
``variable-Lebesgue'' setting and is stated for arbitrary generalized $\Phi$-functions.
It is a refined version of~\cite[Lemma~5.7.14]{DHHR11}, where we remove the 
unnecessary assumptions that $\varphi$ is a generalized $N$-function on $X$ and $\varphi^*$ 
satisfies the $\Delta_2$-condition.
\begin{lemma}\label{le:means-of-char}
Let $\varphi\in\Phi(X)$ be proper and satisfy the $\Delta_2$-condition. For some 
$s\ge1$, let $M_{s,Q}\varphi,(M_{s,Q}\varphi^*)^*\in\Phi(\mathscr{D})$ be such that 
$M_{s,Q}\varphi\preceq(M_{s,Q}\varphi^*)^*$. Then uniformly in $Q\in\mathscr{D}$,
\begin{equation}\label{eq:two-equiv}
\mu(Q)\,(M_{s,Q}\varphi)\left(\frac1{\|\chi_Q\|_\varphi}\right)\approx1
\quad\text{and}\quad
\mu(Q)\,(M_{s,Q}\varphi^*)^*\left(\frac1{\|\chi_Q\|_\varphi}\right)\approx1.
\end{equation}
\end{lemma}
\begin{proof}
Notice that since $\varphi$ is proper, we have 
$M_Q\varphi,(M_Q\varphi^*)^*\in\Phi(\mathscr{D})$
by Lemma~\ref{le:mean-Phi}. Jensen's inequality implies
$(M_Q\varphi)(t)\le(M_{s,Q}\varphi)(t)$ and $(M_Q\varphi^*)(t)\le(M_{s,Q}\varphi^*)(t)$
for all $Q\in\mathscr{D}$ and $t\ge0$, the latter being equivalent to 
$(M_{s,Q}\varphi^*)^*(t)\le(M_Q\varphi^*)^*(t)$ by Lemma~\ref{le:conjugation}(a).
As we remarked after Definition~\ref{def:domination}, these pointwise inequalities 
trivially yield dominations $M_Q\varphi\preceq M_{s,Q}\varphi$ and 
$(M_{s,Q}\varphi^*)^*\preceq(M_Q\varphi^*)^*$. Since we also have 
$(M_Q\varphi^*)^*(t)\le(M_Q\varphi)(t)$ for all $Q\in\mathscr{D}$ and $t\ge0$ 
by Lemma~\ref{le:mirror-inequality}, it follows that
\[
M_{s,Q}\varphi\preceq(M_{s,Q}\varphi^*)^*\preceq(M_Q\varphi^*)^*\preceq
M_Q\varphi\preceq M_{s,Q}\varphi,
\]
and hence 
\begin{equation}\label{eq:equiv-domin}
M_Q\varphi\cong M_{s,Q}\varphi\cong(M_{s,Q}\varphi^*)^*.
\end{equation}
Consider the family $\mathcal{Q}=\{Q\}$ consisting of a single cube $Q\in\mathscr{D}$. 
Clearly, $M_Q\varphi\in\Phi(\mathcal{Q})$.
For the corresponding semimodular $\rho_{M_Q\varphi}$ on 
$\mathbb{C}^\mathcal{Q}=\mathbb{C}$, we obtain
\[
\rho_{M_Q\varphi}\left(\frac1{\|\chi_Q\|_\varphi}\right)
=\mu(Q)\,(M_Q\varphi)\left(\frac1{\|\chi_Q\|_\varphi}\right)
=\rho_\varphi\left(\frac{\chi_Q}{\|\chi_Q\|_\varphi}\right)=1
\]
by the unit sphere property in Lemma~\ref{le:delta2}(c) as a consequence of the 
$\Delta_2$-condition 
for $\varphi$. This, equivalences in terms of domination~\eqref{eq:equiv-domin},
and the fact that $M_Q\varphi$, $M_{s,Q}\varphi$ and $(M_{s,Q}\varphi^*)^*$ satisfy the 
$\Delta_2$-condition together with $\varphi$ (see Lemma~\ref{le:delta2-means})
give the desired 
equivalences~\eqref{eq:two-equiv}, uniformly over $Q\in\mathscr{D}$, 
by Lemma~\ref{le:mod-equiv-1}. 
\end{proof}

The next lemma already depends on the specific form of the ``variable-Lebesgue'' functions
$\varphi_{p(\cdot)}$. Being a counterpart of~\cite[Lemma~5.7.16]{DHHR11}, our 
Lemma~\ref{le:alpha}, however, features a much more elaborate proof than the original.
We draw the reader's special attention to \emph{Step~3} in the proof below, as it clarifies
a rather non-trivial transition that was left unexplained in the original argument.
\begin{lemma}\label{le:alpha}
Let $p(\cdot)\in\mathcal{P}(X)$ with $1<p_-\le p_+<\infty$. Let 
$M_{s,Q}\varphi_{p(\cdot)}\preceq(M_{s,Q}\varphi^*_{p(\cdot)})^*$ for some $s\ge1$.
\emph{Define the ratio function} 
\begin{align}\label{eq:def-alpha}
\alpha_{s,Q}\ :\ &\mathscr{D}\times(0,\infty)\to(0,\infty),
\nonumber \\
&(Q,t)\mapsto\frac{(M_{s,Q}\varphi_{p(\cdot)})(t)}{(M_{s,Q}\varphi^*_{p(\cdot)})^*(t)}
=:\alpha_{s,Q}(t)
\end{align}
\emph{with the convention to write $\alpha_Q$ instead of $\alpha_{1,Q}$.}
Then for any $C_1,C_2>0$ there exists $C_3\ge1$ such that for all $Q\in\mathscr{D}$,
\[
t\in\left[C_1\min\left\{1,\frac1{\|\chi_Q\|_{p(\cdot)}}\right\},\,
C_2\max\left\{1,\frac1{\|\chi_Q\|_{p(\cdot)}}\right\}\right]
\implies\alpha_{s,Q}(t)\le C_3.
\]
\end{lemma}
\begin{proof}
Note first that the condition $p_+<\infty$ implies that $\varphi_{p(\cdot)}$ satisfies
the $\Delta_2$-condition with constant $D:=2^{p_+}$, and so do $M_{s,Q}\varphi_{p(\cdot)}$
and $(M_{s,Q}\varphi^*_{p(\cdot)})^*$ by Lemma~\ref{le:delta2-means}. Then for any
functions $\gamma_1,\gamma_2\ :\ \mathscr{D}\times[0,\infty)\to[0,\infty)$ satisfying 
\[
\gamma_1(Q,t)\approx\gamma_2(Q,t)
\]
uniformly in $Q\in\mathscr{D}$ and $t\ge0$, we have
\begin{align*}
(M_{s,Q}\varphi_{p(\cdot)})(\gamma_1(Q,t))&\approx 
(M_{s,Q}\varphi_{p(\cdot)})(\gamma_2(Q,t)),
\\
(M_{s,Q}\varphi^*_{p(\cdot)})^*(\gamma_1(Q,t))&\approx
(M_{s,Q}\varphi^*_{p(\cdot)})^*(\gamma_2(Q,t))
\end{align*}
uniformly in $Q\in\mathscr{D}$ and $t\ge0$. Indeed, since there is a number 
$m\in\mathbb{N}$ such that $2^{-m}\,\gamma_1(Q,t)\le\gamma_2(Q,t)\le 2^m\,\gamma_1(Q,t)$
holds for all $Q\in\mathscr{D}$ and all $t\ge0$, the above equivalences follow from 
the inequalities
\begin{align*}
(M_{s,Q}\varphi_{p(\cdot)})(\gamma_2(Q,t))&\le
D^m\,(M_{s,Q}\varphi_{p(\cdot)})(\gamma_1(Q,t)), 
\\
(M_{s,Q}\varphi_{p(\cdot)})(\gamma_2(Q,t))&\ge
D^{-m}(M_{s,Q}\varphi_{p(\cdot)})(\gamma_1(Q,t))
\end{align*}
and the analogous inequalities for $(M_{s,Q}\varphi^*_{p(\cdot)})^*$. We will implicitly 
use these equivalences while passing to equivalent arguments inside
$M_{s,Q}\varphi_{p(\cdot)}$ and $(M_{s,Q}\varphi^*_{p(\cdot)})^*$.
With this in mind, we begin the proof of the lemma.
We prove it for $\varphi_{p(\cdot)}=\bar{\varphi}_{p(\cdot)}$.

\medskip\noindent
{\it Step 1: The proof that $\alpha_{s,Q}(1)\approx1$.}
Consider first the auxiliary function 
\begin{align*}
\beta_{s,Q}\ :\ &\mathscr{D}\times[0,\infty)\to[0,\infty),
\\
&(Q,t)\mapsto(M_{s,Q}(\bar{\varphi}^*_{p(\cdot)})')(t)=:\beta_{s,Q}(t).
\end{align*}
It is going to play the key role in the proof. Recalling 
from the discussion after Remark~\ref{re:N-function} that 
$(\bar{\varphi}^*_{p(\cdot)})'=(\bar{\varphi}'_{p(\cdot)})^{-1}$, and using  
$1<p_-\le p_+<\infty$, we obtain from a direct computation that
\[
\beta_{s,Q}(t)=
\left(
\fint_Q \left(\frac{t}{p(x)}\right)^{\frac{s}{p(x)-1}}d\mu(x)
\right)^{1/s}.
\]
Let us list a few properties of this function:
\begin{enumerate}
    \item[(a)] $\beta_{s,Q}(0)=0$ for all $Q\in\mathscr{D}$;
    \item[(b)] $t\mapsto\beta_{s,Q}(t)$ is increasing on $[0,\infty)$ and, by the
    monotone convergence theorem applied as in the proof of Lemma~\ref{le:mean-Phi},
    $\beta_{s,Q}(t)\to\infty$ as $t\to\infty$
    for all $Q\in\mathscr{D}$;
    \item[(c)] $\beta_{s,Q}(t)\approx(\fint_Q t^\frac{s}{p(x)-1}\,d\mu(x))^{1/s}$
    uniformly in $Q\in\mathscr{D}$ and $t\ge0$, 
    since due to the condition $1<p_-\le p_+<\infty$ we have, for all $x\in X$, that
    \[
    1<p_-^\frac{s}{p_+-1}\le p(x)^{\frac{s}{p(x)-1}}\le p_+^\frac{s}{p_--1}<\infty;
    \]
    \item[(d)] $\beta_{s,Q}(1)\approx1$ uniformly in $Q\in\mathscr{D}$ 
    as a consequence of (c);
    \item[(e)] $\beta_{s,Q}$ is finite-valued in view of (c) and 
    $1<p_-\le p_+<\infty$;
    \item[(f)] $t\mapsto\beta_{s,Q}(t)$ is continuous for all $Q\in\mathscr{D}$, because 
    for each sequence $t_k\to t$, there is a number $C>0$ such that all $t_k\le C$, and hence
    \[
    \left(\frac{t_k}{p(\cdot)}\right)^{\frac{s}{p(\cdot)-1}}\le
    \left(\frac{C}{p(\cdot)}\right)^{\frac{s}{p(\cdot)-1}}\in L^1(Q)
    \]
    due to (e); consequently, the dominated convergence theorem~\cite[Theorem~2.24]{GF99}
    yields
    \[
    \lim_{k\to\infty}\beta_{s,Q}(t_k)=
    \left(\fint_Q\lim_{k\to\infty}\left(\frac{t_k}{p(x)}\right)^
    {\frac{s}{p(x)-1}}d\mu(x)\right)^{1/s}=\beta_{s,Q}(t).
    \]
\end{enumerate}
Next, since $\bar{\varphi}^*_{p(\cdot)}$ belongs to the class $N(X)$ and satisfies the 
$\Delta_2$-condition in view of Lemma~\ref{le:phi-bar-tilde}, it follows by 
Lemma~\ref{le:N-ineq}(c) that
$\bar{\varphi}^*_{p(\cdot)}(x,t)\approx t\,(\bar{\varphi}^*_{p(\cdot)})'(x,t)$ uniformly 
in $x\in X$ and $t\ge0$. Hence  
\[
\beta_{s,Q}(t)\approx
\frac{(M_{s,Q}\bar{\varphi}^*_{p(\cdot)})(t)}{t}
\]
uniformly in $Q\in\mathscr{D}$ and $t>0$. This and Lemma~\ref{le:N-ineq}(d),
applied to the function $(M_{s,Q}\bar{\varphi}^*_{p(\cdot)})^*\in N(\mathscr{D})$ 
satisfying the $\Delta_2$-condition, together give
\begin{align}
(M_{s,Q}\bar{\varphi}^*_{p(\cdot)})^*(\beta_{s,Q}(t))
&\approx(M_{s,Q}\bar{\varphi}^*_{p(\cdot)})^*
\left(\frac{(M_{s,Q}\bar{\varphi}^*_{p(\cdot)})(t)}{t}\right)
\nonumber \\
&\approx(M_{s,Q}\bar{\varphi}^*_{p(\cdot)})(t)
\approx t\,\beta_{s,Q}(t)
\label{eq:gig-1}
\end{align}
uniformly in $Q\in\mathscr{D}$ and $t>0$. Since $\beta_{s,Q}(1)\approx1$, the above 
equivalence yields
\[
(M_{s,Q}\bar{\varphi}^*_{p(\cdot)})^*(1)\approx
(M_{s,Q}\bar{\varphi}^*_{p(\cdot)})^*(\beta_{s,Q}(1))\approx
1\times\beta_{s,Q}(1)\approx1
\]
uniformly in $Q\in\mathscr{D}$. Combining this with the simple calculation 
$(M_{s,Q}\bar{\varphi}_{p(\cdot)})(1)=1$, we obtain $\alpha_{s,Q}(1)\approx1$ 
uniformly in $Q\in\mathscr{D}$.

\medskip\noindent
{\it Step 2: Useful inequalities for $\alpha_{s,Q}(\beta_{s,Q}(t))$.}
Applying equivalence~\eqref{eq:gig-1} together with 
$\beta_{s,Q}(t)\approx(\fint_Q t^\frac{s}{p(x)-1}\,d\mu(x))^{1/s}$
and the identity
\[
\frac{(M_{s,Q}\bar{\varphi}_{p(\cdot)})(t)}{t}=
\left(\fint_Q t^{s(p(y)-1)}\,d\mu(y)\right)^{1/s},
\]
we get, uniformly in $Q\in\mathscr{D}$ and $t>0$, that  

\begin{align*}\label{eq:gig-2}
\alpha_{s,Q}(\beta_{s,Q}(t))
&=\frac{(M_{s,Q}\bar{\varphi}_{p(\cdot)})
(\beta_{s,Q}(t))}
{(M_{s,Q}\bar{\varphi}^*_{p(\cdot)})^*
(\beta_{s,Q}(t))}
\approx\frac{(M_{s,Q}\bar{\varphi}_{p(\cdot)})
(\beta_{s,Q}(t))}
{t\,\beta_{s,Q}(t)}
\\
&\approx\frac{(M_{s,Q}\bar{\varphi}_{p(\cdot)})
((\fint_Q t^\frac{s}{p(x)-1}\,d\mu(x))^{1/s})}
{t\,(\fint_Q t^\frac{s}{p(x)-1}\,d\mu(x))^{1/s}}
\\
&=\frac1t\left(
\fint_Q\left(\fint_Q t^\frac{s}{p(x)-1}\,d\mu(x)\right)^{p(y)-1}d\mu(y)
\right)^{1/s}
\\
&=\left(
\fint_Q\left(\fint_Q t^\frac{s(p(y)-p(x))}{(p(x)-1)(p(y)-1)}\,d\mu(x)\right)^{p(y)-1}d\mu(y)
\right)^{1/s}.
\end{align*}
Define
\begin{align*}
\alpha_{s,Q}^{>}(t)&:=\left(
\fint_Q\left(\fint_Q t^\frac{s(p(y)-p(x))}{(p(x)-1)(p(y)-1)}
\chi_{\{x\in Q\,:\,p(y)>p(x)\}}(x)\,d\mu(x)\right)^{p(y)-1}d\mu(y)
\right)^{1/s},
\\
\alpha_{s,Q}^{\le}(t)&:=\left(
\fint_Q\left(\fint_Q t^\frac{s(p(y)-p(x))}{(p(x)-1)(p(y)-1)}
\chi_{\{x\in Q\,:\,p(y)\le p(x)\}}(x)\,d\mu(x)\right)^{p(y)-1}d\mu(y)
\right)^{1/s}.
\end{align*}
Then, uniformly in $Q\in\mathscr{D}$ and $t>0$, we have the equivalence
\[
\alpha_{s,Q}(\beta_{s,Q}(t))\approx\alpha_{s,Q}^{>}(t)+\alpha_{s,Q}^{\le}(t),
\]
which follows from the simple observations that:
\begin{itemize}
    \item[(i)] $(f+g)^{p(y)-1}\approx f^{p(y)-1}+g^{p(y)-1}$ uniformly in $y\in X$ and 
    $f,g\ge0$, due to
    \begin{align*}
    \frac{f^{p(y)-1}+g^{p(y)-1}}2\le(f+g)^{p(y)-1}&\le\max{\{1,2^{p(y)-2}\}}\,
    (f^{p(y)-1}+g^{p(y)-1})
    \\
    &\le\max\{1,2^{p_+-2}\}\,(f^{p(y)-1}+g^{p(y)-1});
    \end{align*}
    \item[(ii)] $(a+b)^{1/s}\approx a^{1/s}+b^{1/s}$ uniformly in $a,b\ge0$, since in view 
    of $1/s\in(0,1]$,
    \[
    \frac{a^{1/s}+b^{1/s}}2\le(a+b)^{1/s}\le a^{1/s}+b^{1/s}.
    \]
\end{itemize}
Therefore, there exists a constant $C\ge1$ such that for all $Q\in\mathscr{D}$ and $t>0$,
\begin{equation}\label{eq:gig-2}
\frac1C\,\alpha_{s,Q}(\beta_{s,Q}(t))\,\le\,\alpha_{s,Q}^{>}(t)+\alpha_{s,Q}^{\le}(t)
\,\le\, C\,\alpha_{s,Q}(\beta_{s,Q}(t)).
\end{equation}

By definition, for all $Q\in\mathscr{D}$,
\begin{enumerate}
    \item[] $t\mapsto\alpha_{s,Q}^{>}(t)$ is non-decreasing on $(0,\infty)$;
    \item[] $t\mapsto\alpha_{s,Q}^{\le}(t)$ is non-increasing on $(0,\infty)$;
    \item[] $0\le\alpha_{s,Q}^{>}(t)\le1$ for $0<t\le1$;
    \item[] $0\le\alpha_{s,Q}^{\le}(t)\le1$ for $t\ge1$.
\end{enumerate}
In view of these properties, equivalence \eqref{eq:gig-2} implies that 
for any $0<t_1\le t_2\le1$,
\begin{align}
\alpha_{s,Q}(\beta_{s,Q}(t_2))&\le C\,(\alpha_{s,Q}^{>}(t_2)+\alpha_{s,Q}^{\le}(t_2))
\le C\,(1+\alpha_{s,Q}^{\le}(t_1))
\nonumber \\
&\le C\,(C\,\alpha_{s,Q}(\beta_{s,Q}(t_1))+1)
\le C^2\,(\alpha_{s,Q}(\beta_{s,Q}(t_1))+1),
\label{eq:t<1}
\end{align}
while for any $1\le t_3\le t_4$, we have
\begin{align}
\alpha_{s,Q}(\beta_{s,Q}(t_3))&\le C\,(\alpha_{s,Q}^{>}(t_3)+\alpha_{s,Q}^{\le}(t_3))
\le C\,(\alpha_{s,Q}^{>}(t_4)+1)
\nonumber \\
&\le C\,(C\,\alpha_{s,Q}(\beta_{s,Q}(t_4))+1)
\le C^2\,(\alpha_{s,Q}(\beta_{s,Q}(t_4))+1)
\label{eq:t>1}
\end{align}
for all cubes $Q\in\mathscr{D}$.

\medskip\noindent
{\it Step 3: Corresponding inequalities for $\alpha_{s,Q}(t)$.}
Note that for any $K>0$, there holds
\begin{equation}\label{eq:gig-3}
\alpha_{s,Q}(t)\approx\alpha_{s,Q}(Kt)
\end{equation}
uniformly in
$Q\in\mathscr{D}$ and $t>0$. This follows from the corresponding equivalences
$(M_{s,Q}\bar{\varphi}_{p(\cdot)})(Kt)\approx(M_{s,Q}\bar{\varphi}_{p(\cdot)})(t)$ 
and $(M_{s,Q}\bar{\varphi}^*_{p(\cdot)})^*(Kt)\approx
(M_{s,Q}\bar{\varphi}^*_{p(\cdot)})^*(t)$, uniformly in $Q\in\mathscr{D}$ and $t>0$, 
which are true since $Kt\approx t$ uniformly in $t>0$.
By property~(d) of the function $\beta_{s,Q}$, there exists $K_0\ge1$ such that 
$\frac1{K_0}\le\beta_{s,Q}(1)\le K_0$ for all $Q\in\mathscr{D}$. Denote 
\[
\widetilde{\beta}_{s,Q}(t):=K_0\,\beta_{s,Q}(t).
\]
Then for any $Q\in\mathscr{D}$, we have $\widetilde{\beta}_{s,Q}(0)=0$, 
$\widetilde{\beta}_{s,Q}(1)\ge1$, and the function
$t\mapsto\widetilde{\beta}_{s,Q}(t)$ is continuously increasing on $[0,\infty)$. 
This ensures that for any $0<\tau_1\le\tau_2\le1$, the values 
\[
t_1:=\widetilde{\beta}_{s,Q}^{-1}(\tau_1)
\quad\text{and}\quad
t_2:=\widetilde{\beta}_{s,Q}^{-1}(\tau_2),
\]
where $\widetilde{\beta}_{s,Q}^{-1}$ denotes the inverse of $\widetilde{\beta}_{s,Q}$ 
with respect to the second (real) variable,
satisfy 
\[
0<t_1\le t_2\le1
\]
for all $Q\in\mathscr{D}$. Therefore, by equivalence~\eqref{eq:gig-3} and
inequality~\eqref{eq:t<1},
\begin{align*}
\alpha_{s,Q}(\tau_2)&=\alpha_{s,Q}(\widetilde{\beta}_{s,Q}(t_2))
=\alpha_{s,Q}(K_0\,\beta_{s,Q}(t_2))
\\
&\approx\alpha_{s,Q}(\beta_{s,Q}(t_2))\le
C^2\,(\alpha_{s,Q}(\beta_{s,Q}(t_1))+1)
\\
&\approx\alpha_{s,Q}(K_0\,\beta_{s,Q}(t_1))+1=\alpha_{s,Q}(\tau_1)+1
\end{align*}
uniformly in $Q\in\mathscr{D}$ and $0<\tau_1\le\tau_2\le1$. Hence there exists a 
constant $\widetilde{C}\ge1$ such that for all $Q\in\mathscr{D}$, we have
\begin{equation}\label{eq:gig-4}
\alpha_{s,Q}(\tau_2)\le \widetilde{C}\,(\alpha_{s,Q}(\tau_1)+1)
\quad\text{if }0<\tau_1\le\tau_2\le1.
\end{equation}

Similarly, denote
\[
\widehat{\beta}_{s,Q}(t):=\frac1{K_0}\,\beta_{s,Q}(t).
\]
Then for any $Q\in\mathscr{D}$, we have that $\widehat{\beta}_{s,Q}(1)\le1$, the function
$t\mapsto\widehat{\beta}_{s,Q}(t)$ is continuously increasing on $[0,\infty)$, and 
$\widehat{\beta}_{s,Q}(t)\to\infty$ as $t\to\infty$. This ensures that for any 
$1\le\tau_3\le\tau_4$, the values
\[
t_3:=\widehat{\beta}_{s,Q}^{-1}(\tau_3)
\quad\text{and}\quad
t_4:=\widehat{\beta}_{s,Q}^{-1}(\tau_4),
\]
where $\widehat{\beta}_{s,Q}^{-1}$ denotes the inverse of $\widehat{\beta}_{s,Q}$ with
respect to the second (real) variable, satisfy
\[
1\le t_3\le t_4
\]
for all $Q\in\mathscr{D}$. Therefore, by equivalence~\eqref{eq:gig-3} and 
inequality~\eqref{eq:t>1}, 
\begin{align*}
\alpha_{s,Q}(\tau_3)&=\alpha_{s,Q}(\widehat{\beta}_{s,Q}(t_3))
=\alpha_{s,Q}\left(\frac{1}{K_0}\,\beta_{s,Q}(t_3)\right)
\\
&\approx\alpha_{s,Q}(\beta_{s,Q}(t_3))\le
C^2\,(\alpha_{s,Q}(\beta_{s,Q}(t_4))+1)
\\
&\approx\alpha_{s,Q}\left(\frac1{K_0}\,\beta_{s,Q}(t_4)\right)+1=\alpha_{s,Q}(\tau_4)+1
\end{align*}
uniformly in $Q\in\mathscr{D}$ and $1\le\tau_3\le\tau_4$.
Hence there exists a constant $\widehat{C}\ge1$ such that for all $Q\in\mathscr{D}$,
we have
\begin{equation}\label{eq:gig-5}
\alpha_{s,Q}(\tau_3)\le \widehat{C}\,(\alpha_{s,Q}(\tau_4)+1)
\quad\text{if }1\le\tau_3\le\tau_4.
\end{equation}

\medskip\noindent
{\it Step 4: Proof of the main claim.} In \emph{Step~1}, we showed that 
$\alpha_{s,Q}(1)\approx1$ uniformly in $Q\in\mathscr{D}$. It easily follows by
Lemma~\ref{le:means-of-char} that we also have
$\alpha_{s,Q}(1/\|\chi_Q\|_{p(\cdot)})\approx1$ uniformly in $Q\in\mathscr{D}$. 
These two facts combined give
\[
\alpha_{s,Q}\left(\min\left\{1,\frac1{\|\chi_Q\|_{p(\cdot)}}\right\}\right)\approx1
\quad\text{and}\quad
\alpha_{s,Q}\left(\max\left\{1,\frac1{\|\chi_Q\|_{p(\cdot)}}\right\}\right)\approx1
\]
uniformly in $Q\in\mathscr{D}$. For an arbitrary $Q\in\mathscr{D}$ and some
$C_1,C_2>0$, let us take a value
\begin{equation}\label{eq:gig-t}
t\in\left[C_1\min\left\{1,\frac1{\|\chi_Q\|_{p(\cdot)}}\right\},\,
C_2\max\left\{1,\frac1{\|\chi_Q\|_{p(\cdot)}}\right\}\right].
\end{equation}
For $t<1$, inequality~\eqref{eq:gig-4} yields
\begin{align*}
\alpha_{s,Q}(t)&\le\widetilde{C}\,\left(\alpha_{s,Q}\left(C_1
\min\left\{1,\frac1{\|\chi_Q\|_{p(\cdot)}}\right\}\right)+1\right)
\\
&\approx\alpha_{s,Q}\left(
\min\left\{1,\frac1{\|\chi_Q\|_{p(\cdot)}}\right\}\right)+1
\approx1
\end{align*}
uniformly in $Q\in\mathscr{D}$ and $t<1$; otherwise, if $t\ge1$, 
inequality~\eqref{eq:gig-5} gives
\begin{align*}
\alpha_{s,Q}(t)&\le\widehat{C}\,\left(\alpha_{s,Q}\left(C_2
\max\left\{1,\frac1{\|\chi_Q\|_{p(\cdot)}}\right\}\right)+1\right)
\\
&\approx\alpha_{s,Q}\left(
\max\left\{1,\frac1{\|\chi_Q\|_{p(\cdot)}}\right\}\right)+1
\approx1
\end{align*}
uniformly in $Q\in\mathscr{D}$ and $t\ge1$.
The above two estimates guarantee that there exists $C_3\ge1$ such that 
$\alpha_{s,Q}(t)\le C_3$ for all $Q\in\mathscr{D}$ and any 
$t$ from the interval~\eqref{eq:gig-t}. This proves the lemma for 
$\varphi_{p(\cdot)}=\bar{\varphi}_{p(\cdot)}$.

For $\varphi_{p(\cdot)}=\widetilde{\varphi}_{p(\cdot)}$, the same proof applies almost 
verbatim, with natural amendments in the steps involving direct computations. Besides,
\emph{Step~3} becomes significantly simpler in the case 
$\varphi_{p(\cdot)}=\widetilde{\varphi}_{p(\cdot)}$, since one has $K_0=1$ and thus 
can take $\widetilde{\beta}_{s,Q}=\widehat{\beta}_{s,Q}=\beta_{s,Q}$.
\end{proof}

The next lemma is a refinement of Theorem~\ref{th:pointwise-domin} in the case of specific
domination $M_{s,Q}\varphi_{p(\cdot)}\preceq(M_{s,Q}\varphi^*_{p(\cdot)})^*$. 
Building on the above Lemma~\ref{le:alpha}, it establishes
pointwise inequalities implied by this domination in analogy 
with~\cite[Lemma~5.7.20]{DHHR11}.
\begin{lemma}\label{le:p-pointwise}
Let $p(\cdot)\in\mathcal{P}(X)$ with $1<p_-\le p_+<\infty$. Let 
$M_{s,Q}\varphi_{p(\cdot)}\preceq(M_{s,Q}\varphi^*_{p(\cdot)})^*$ for some $s\ge1$.
Then there exist $K_0\ge1$ and $b:\mathscr{D}\to[0,\infty)$ with
$\|b\|_{\mathfrak{S},1}<\infty$ such that for all $Q\in\mathscr{D}$ and $t\ge0$
satisfying
\[
\mu(Q)\,(M_{s,Q}\varphi^*_{p(\cdot)})^*(t)\le1,
\]
there holds
\[
(M_{s,Q}\varphi_{p(\cdot)})(t)\le 
\left\{\begin{array}{ll}
     K_0\,(M_{s,Q}\varphi^*_{p(\cdot)})^*(t)+b(Q) &\text{if }0\le t<1,\\
     K_0\,(M_{s,Q}\varphi^*_{p(\cdot)})^*(t) &\text{if }t\ge1. 
\end{array}\right.
\]
\end{lemma}
\begin{proof}
Note that $\varphi_{p(\cdot)}$ is proper and satisfies the $\Delta_2$-condition due to 
Lemma~\ref{le:phi-bar-tilde}.
Let $\mu(Q)\,(M_{s,Q}\varphi^*_{p(\cdot)})^*(t)\le1$. Assume first that $t\ge1$. By
Lemma ~\ref{le:means-of-char}, there exists $C\ge1$ such that for all $Q\in\mathscr{D}$
we have
\[
\frac1C\le
\mu(Q)\,(M_{s,Q}\varphi^*_{p(\cdot)})^*\left(\frac1{\|\chi_Q\|_{p(\cdot)}}\right)
\le C,
\]
so it follows by assumption and the convexity of $(M_{s,Q}\varphi^*_{p(\cdot)})^*$
in the second variable that
\begin{align*}
\mu(Q)\,(M_{s,Q}\varphi^*_{p(\cdot)})^*(t)\le1&\le 
C\,\mu(Q)\,(M_{s,Q}\varphi^*_{p(\cdot)})^*\left(\frac1{\|\chi_Q\|_{p(\cdot)}}\right)
\\
&\le\mu(Q)\,(M_{s,Q}\varphi^*_{p(\cdot)})^*\left(\frac{C}{\|\chi_Q\|_{p(\cdot)}}\right).
\end{align*}
Since $t\mapsto(M_{s,Q}\varphi^*_{p(\cdot)})^*(t)$ is increasing due to 
$(M_{s,Q}\varphi^*_{p(\cdot)})^*\in N(\mathscr{D})$, see Remark~\ref{re:N-means}, we derive 
from the above inequality that $1\le t\le C/\|\chi_Q\|_{p(\cdot)}$ with the constant
$C$ independent of a cube $Q\in\mathscr{D}$. Then by Lemma~\ref{le:alpha} applied with 
$C_1:=1$ and $C_2:=C$, there is $C_3\ge1$ such that for all $Q\in\mathscr{D}$,
\begin{equation}\label{eq:sp-point-1}
(M_{s,Q}\varphi_{p(\cdot)})(t)\le C_3\,(M_{s,Q}\varphi^*_{p(\cdot)})^*(t).
\end{equation}

Otherwise, if $0\le t<1$, we use the fact that $M_{s,Q}\varphi_{p(\cdot)}$ satisfies the 
$\Delta_2$-condition and 
apply Theorem~\ref{th:pointwise-domin} to the specific
domination $M_{s,Q}\varphi_{p(\cdot)}\preceq(M_{s,Q}\varphi^*_{p(\cdot)})^*$ to
find $A_2>1$ and a function $b:\mathscr{D}\to[0,\infty)$ with 
$\|b\|_{\mathfrak{S},1}\le A_2<\infty$
such that for all $Q\in\mathscr{D}$,
\begin{equation}\label{eq:sp-point-2}
(M_{s,Q}\varphi_{p(\cdot)})(t)\le A_2\,(M_{s,Q}\varphi^*_{p(\cdot)})^*(t)+b(Q).
\end{equation}
Inequalities~\eqref{eq:sp-point-1} and \eqref{eq:sp-point-2} yield the claim with 
$K_0:=\max\{C_3,A_2\}.$
\end{proof}

The following Lemma~\ref{le:gammas} is preparatory for the proof of 
Lemma~\ref{le:alpha-s-1}. It is a simplified version of~\cite[Lemma~5.7.23]{DHHR11},
with a refined proof removing a vicious cycle hidden in the original. 
In the original proof, to apply~\cite[Remark~2.6.7]{DHHR11}, one had to know that the
function $\gamma$ (the same as in the proof of Lemma~\ref{le:gammas} below) satisfies 
the $\Delta_2$-condition;
however, the $\Delta_2$-condition for $\gamma$ was mistakenly inferred from an
argument involving~\cite[Remark~2.6.7]{DHHR11} itself. We break this vicious cycle by 
establishing the $\Delta_2$-condition for $\gamma$ before we apply 
Lemma~\ref{le:N-ineq}(c), which corresponds to \cite[Remark~2.6.7]{DHHR11} in our text.
\begin{lemma}\label{le:gammas}
Let $\varphi\in N(\mathscr{D})$ be such that $\varphi$ and $\varphi^*$ satisfy the 
$\Delta_2$-condition. Suppose that $\psi\in N(\mathscr{D})$ satisfies 
$\psi(Q,t)\approx(\varphi(Q,t^{1/s}))^s$, uniformly in $Q\in\mathscr{D}$ and $t\ge0$, 
for some $s>0$. Then $\psi^*(Q,t)\approx(\varphi^*(Q,t^{1/s}))^s$ uniformly in 
$Q\in\mathscr{D}$ and $t\ge0$.
\end{lemma}
\begin{proof}
Consider the auxiliary function $\gamma:\mathscr{D}\times[0,\infty)\to[0,\infty)$ 
defined by
\[
\gamma(Q,t):=\int_0^t \left(\varphi'(Q,\tau^{1/s})\right)^s d\tau,
\]
where $\varphi'$ is the right derivative of $\varphi\in N(\mathscr{D})$ with respect to the 
second variable.
Let us check that $\gamma\in N(\mathscr{D})$. Since $\varphi\in N(\mathscr{D})$, it 
follows by Remark~\ref{re:N-function} that for any $Q\in\mathscr{D}$, the right derivative
$t\mapsto\varphi'(Q,t)$ is non-decreasing and right-continuous, 
$\varphi'(Q,0)=0$, $\varphi'(Q,t)>0$ for $t>0$, and $\varphi'(Q,t)\to\infty$ as 
$t\to\infty$. Due to the equality $\gamma'(Q,t)=(\varphi'(Q,t^{1/s}))^s$, the right
derivative $t\mapsto\gamma'(Q,t)$ satisfies exactly the same properties for every
$Q\in\mathscr{D}$. Again by Remark~\ref{re:N-function}, this is equivalent to 
$\gamma\in N(\mathscr{D})$.

Lemma~\ref{le:N-ineq}(a) applied to $\gamma$, 
equality $\gamma'(Q,t)=(\varphi'(Q,t^{1/s}))^s$,
Lemma~\ref{le:N-ineq}(c) applied to $\varphi$,
and the $\Delta_2$-condition for $\varphi$ with constant $D\ge2$
together give
\begin{align*}
\gamma(Q,2t)&\le2t\,\gamma'(Q,2t)=2t\,\big(\varphi'(Q,(2t)^{1/s})\big)^s
\approx2t\left(\frac{\varphi(Q,(2t)^{1/s})}{(2t)^{1/s}}\right)^s
\\
&=\frac{t}{2}\left(\frac{\varphi(Q,2^{2/s}(\frac{t}{2})^{1/s})}{(\frac{t}{2})^{1/s}}\right)^s
\le\frac{t}{2}\left(\frac{D^{[2/s]+1}\,\varphi(Q,(\frac{t}{2})^{1/s})}
{(\frac{t}{2})^{1/s}}\right)^s
\\
&\approx\frac{t}{2}\left(\varphi'\left(Q,\Big(\dfrac{t}{2}\Big)^{1/s}\right)\right)^s
=\frac{t}{2}\,\gamma'\left(Q,\frac{t}{2}\right)\le\gamma(Q,t)
\end{align*}
uniformly in $Q\in\mathscr{D}$ and $t>0$, hence $\gamma$ satisfies the $\Delta_2$-condition.
We also deduce from $\gamma'(Q,t)=(\varphi'(Q,t^{1/s}))^s$ that
\begin{align*}
(\gamma^*)'(Q,t)=(\gamma')^{-1}(Q,t)
&=\inf\{u\ge0\ :\ \gamma'(Q,u)>t\}
\\
&=\inf\left\{u\ge0\ :\ \big(\varphi'(Q,u^{1/s})\big)^s>t\right\}
\\
&=\left(\inf\{u\ge0\ :\ \varphi'(Q,u)>t^{1/s}\}\right)^s
\\
&=\left((\varphi')^{-1}(Q,t^{1/s})\right)^s
=\left((\varphi^*)'(Q,t^{1/s})\right)^s,
\end{align*}
and since $\varphi^*$ satisfies the $\Delta_2$-condition, we repeat the above chain of
inequalities with $\gamma^*$ and $\varphi^*$ instead of $\gamma$ and $\varphi$, 
respectively, and conclude that $\gamma^*$ also satisfies the $\Delta_2$-condition. 

Now Lemma~\ref{le:N-ineq}(c) applied to $\gamma$ and $\varphi$, both of which satisfy
the $\Delta_2$-condition, gives
\[
\gamma(Q,t)\approx t\,\gamma'(Q,t)=\big(t^{1/s}\,\varphi'(Q,t^{1/s})\big)^s
\approx\big(\varphi(Q,t^{1/s})\big)^s\approx\psi(Q,t)
\]
uniformly in $Q\in\mathscr{D}$ and $t\ge0$. Therefore, there is a number
$m\in\mathbb{N}$ such that
\[
2^{-m}\,\gamma(Q,t)\le\psi(Q,t)\le2^m\,\gamma(Q,t)
\]
for all $Q\in\mathscr{D}$ and $t\ge0$. By Lemma~\ref{le:conjugation} and
the $\Delta_2$-condition for $\gamma^*$ with some constant $D^*\ge2$, this implies
\begin{align*}
\psi^*(Q,t)&\ge2^m\,\gamma^*\left(Q,\frac{t}{2^m}\right)
\ge\left(\frac{2}{D^*}\right)^m\gamma^*(Q,t),
\\
\psi^*(Q,t)&\le2^{-m}\,\gamma^*(Q,2^m\,t)\le
\left(\frac{D^*}{2}\right)^m\gamma^*(Q,t),
\end{align*}
hence $\psi^*(Q,t)\approx\gamma^*(Q,t)$ uniformly in $Q\in\mathscr{D}$ and $t\ge0$.
Finally, Lemma~\ref{le:N-ineq}(c) applied to $\gamma^*$ and $\varphi^*$, together with 
$(\gamma^*)'(Q,t)=((\varphi^*)'(Q,t^{1/s}))^s$, implies
\[
\gamma^*(Q,t)\approx t\,(\gamma^*)'(Q,t)=\big(t^{1/s}\,(\varphi^*)'(Q,t^{1/s})\big)^s
\approx\big(\varphi^*(Q,t^{1/s})\big)^s
\]
and, therefore, $\psi^*(Q,t)\approx(\varphi^*(Q,t^{1/s}))^s$ uniformly in $Q\in\mathscr{D}$
and $t\ge0$.
\end{proof}

The final auxiliary Lemma~\ref{le:alpha-s-1} closely follows the proof 
of~\cite[Lemma~5.7.28]{DHHR11}, while providing additional detail for clarity.
\begin{lemma}\label{le:alpha-s-1}
Let $p(\cdot)\in\mathcal{P}(X)$ with $1<p_-\le p_+<\infty$, and
$M_{s,Q}\varphi_{p(\cdot)}\preceq(M_{s,Q}\varphi^*_{p(\cdot)})^*$ for some
$s\ge1$. Then uniformly in $Q\in\mathscr{D}$ and $t>0$, there holds 
\[
\alpha_Q(t)\approx\big(\alpha_{s,Q}(t^{1/s})\big)^s
\]
for the functions $\alpha_Q$ and $\alpha_{s,Q}$ defined as in~\eqref{eq:def-alpha}.
\end{lemma}
\begin{proof}
By Lemma~\ref{le:phi-bar-tilde}(a) we have, uniformly in $x\in X$ and $t\ge0$, that
\[
\varphi_{p(\cdot)}(x,t^s)\approx
\bar{\varphi}_{p(\cdot)}(x,t^s)=\big(\bar{\varphi}_{p(\cdot)}(x,t)\big)^s
\approx\big(\varphi_{p(\cdot)}(x,t)\big)^s,
\]
which implies 
\begin{align}\label{eq:pmean1}
(M_Q\varphi_{p(\cdot)})(t)&=
\fint_Q\varphi_{p(\cdot)}(x,(t^{1/s})^s)\,d\mu(x)
\nonumber \\
&\approx\fint_Q\big(\varphi_{p(\cdot)}(x,t^{1/s})\big)^s d\mu(x)
=\big((M_{s,Q}\varphi_{p(\cdot)})(t^{1/s})\big)^s
\end{align}
uniformly in $Q\in\mathscr{D}$ and $t\ge0$. Similarly, we apply 
Lemma~\ref{le:phi-bar-tilde}(b) to obtain, uniformly in $x\in X$ and $t\ge0$, the chain
of equivalences
\[
\varphi^*_{p(\cdot)}(x,t^s)\approx\varphi_{p'(\cdot)}(x,t^s)
\approx\big(\varphi_{p'(\cdot)}(x,t)\big)^s
\approx\big(\varphi^*_{p(\cdot)}(x,t)\big)^s,
\]
which leads to 
\begin{equation}\label{eq:pmean2}
(M_Q\varphi^*_{p(\cdot)})(t)
\approx\big((M_{s,Q}\varphi^*_{p(\cdot)})(t^{1/s})\big)^s
\end{equation}
uniformly in $Q\in\mathscr{D}$ and $t\ge0$. Note that 
$M_{s,Q}\varphi^*_{p(\cdot)}\in N(\mathscr{D})$ in view of Remark~\ref{re:N-means};
moreover, both $M_{s,Q}\varphi^*_{p(\cdot)}$ and $(M_{s,Q}\varphi^*_{p(\cdot)})^*$
satisfy the $\Delta_2$-condition, since $\varphi_{p(\cdot)}$ and $\varphi^*_{p(\cdot)}$
do so by Lemma~\ref{le:phi-bar-tilde}. Noting also that 
$M_Q\varphi_{p(\cdot)}^*\in N(\mathscr{D})$, which follows by
Lemma~\ref{le:mean-N} because $\varphi_{p(\cdot)}^*\in N(X)$ is proper and satisfies the 
$\Delta_2$-condition, we apply Lemma~\ref{le:gammas} to relation~\eqref{eq:pmean2} 
to conclude that
\[
(M_Q\varphi_{p(\cdot)}^*)^*(t)\approx
\big((M_{s,Q}\varphi^*_{p(\cdot)})^*(t^{1/s})\big)^s
\]
uniformly in $Q\in\mathscr{D}$ and $t\ge0$. This and~\eqref{eq:pmean1} imply
$\alpha_Q(t)\approx(\alpha_{s,Q}(t^{1/s}))^s$ uniformly in $Q\in\mathscr{D}$ and
$t>0$.
\end{proof}

\subsection{Equivalence of the classes $\mathscr{A^D}$ and 
$\mathscr{A}^{\mathscr{D}}_\textnormal{strong}$ for nice ``variable-Lebesgue'' functions}

Having established the auxiliary lemmas, we are now ready to prove the main result of this
section. Theorem~\ref{th:coincidence} shows that for the ``variable-Lebesgue'' functions
$\varphi_{p(\cdot)}$, membership in the class $\mathscr{A^D}$ is equivalent to membership in 
the more restrictive class $\mathscr{A}^{\mathscr{D}}_\textnormal{strong}$, provided that 
the values of the variable exponent $p(\cdot)$ remain strictly between one and infinity,
that is, if $1<p_-\le p_+<\infty$. This result is an analogue 
of~\cite[Theorem~5.7.1]{DHHR11}. 

\begin{theorem}\label{th:coincidence}
Let $p(\cdot)\in\mathcal{P}(X)$ with $1<p_-\le p_+<\infty$. Then $p(\cdot)\in\mathscr{A^D}$ 
if and only if $p(\cdot)\in\mathscr{A}^{\mathscr{D}}_\textnormal{strong}$.
\end{theorem}
\begin{proof}
Let $p(\cdot)\in\mathscr{A^D}$. Since $\varphi_{p(\cdot)}$ is proper and
the assumption $1<p_-\le p_+<\infty$ ensures that both $\varphi_{p(\cdot)}$ and
$\varphi_{p(\cdot)}^*$ satisfy the $\Delta_2$-condition, it follows from the refined 
domination property of $\mathscr{A^D}$ in Theorem~\ref{th:refinement} that
\[
M_{s,Q}\varphi_{p(\cdot)}\preceq(M_{s,Q}\varphi^*_{p(\cdot)})^*
\quad\textnormal{for some }s>1.
\]
This has the following consequences:
\begin{enumerate}
    \item[(a)] By Lemma~\ref{le:means-of-char}, there exists $\lambda\in(0,1]$ such that
for all $Q\in\mathscr{D}$ and nonnegative $t\le\lambda/\|\chi_Q\|_{p(\cdot)}$, there holds
\begin{align*}
\mu(Q)\,(M_{s,Q}\varphi_{p(\cdot)}^*)^*(t)&\le
\mu(Q)\,(M_{s,Q}\varphi_{p(\cdot)}^*)^*\left(\frac{\lambda}{\|\chi_Q\|_{p(\cdot)}}\right)
\\ 
&\le\lambda\,\mu(Q)\,
(M_{s,Q}\varphi_{p(\cdot)}^*)^*\left(\frac{1}{\|\chi_Q\|_{p(\cdot)}}\right)\le1.
\end{align*}
    \item[(b)] Further, by Lemma~\ref{le:alpha}, there exists a constant $C_1\ge1$
such that for all $Q\in\mathscr{D}$,
\[
t\in\left[\lambda\,\min\left\{1,\frac1{\|\chi_Q\|_{p(\cdot)}}\right\},\,
\max\left\{1,\frac1{\|\chi_Q\|_{p(\cdot)}}\right\}\right]
\implies\alpha_{s,Q}(t)\le C_1.
\]
    \item[(c)] Lemma~\ref{le:p-pointwise} holds with some constant $K_0\ge1$ and
function $b:\mathscr{D}\to[0,\infty)$ satisfying 
\[\|b\|_{\mathfrak{S},1}=\sup_{\mathcal{Q}\in\mathfrak{S}}
\sum_{Q\in\mathcal{Q}}\mu(Q)\,b(Q)<\infty,
\]
which correspond to the parameter $s$.
    \item[(d)] By Lemma~\ref{le:alpha-s-1}, there exists $C_2\ge1$ such that 
for all $Q\in\mathscr{D}$ and $t>0$,
\[
\alpha_Q(t)\le C_2\,\big(\alpha_{s,Q}(t^{1/s})\big)^s.
\]
\end{enumerate}
Define $K_1:=C_2\,(\max\{2K_0,C_1\})^s$ and prove the following auxiliary claim.

\medskip\noindent
\emph{Auxiliary claim:} For all $Q\in\mathscr{D}$ and $t>0$ satisfying 
$\mu(Q)\,(M_Q\varphi^*_{p(\cdot)})^*(t)\le1$, we have
\[
(M_Q\varphi_{p(\cdot)})(t)\le\left\{\begin{array}{ll}
     K_1\,(M_Q\varphi^*_{p(\cdot)})^*(t)+2b(Q)\,t^{1-1/s}
     &\textnormal{if }0<t<1,  \\
     K_1\,(M_Q\varphi^*_{p(\cdot)})^*(t) &\textnormal{if }t\ge1. 
\end{array}
\right.
\]
\emph{Proof of auxiliary claim.} By Jensen's inequality, for all $Q\in\mathscr{D}$
and $t\ge0$, there holds $(M_Q\varphi^*_{p(\cdot)})(t)\le(M_{s,Q}\varphi^*_{p(\cdot)})(t)$,
and consequently $(M_{s,Q}\varphi^*_{p(\cdot)})^*(t)\le(M_Q\varphi^*_{p(\cdot)})^*(t)$ 
in view of Lemma~\ref{le:conjugation}(a). By the assumption of the claim,
\[
\mu(Q)\,(M_{s,Q}\varphi^*_{p(\cdot)})^*(t)\le
\mu(Q)\,(M_Q\varphi^*_{p(\cdot)})(t)\le1.
\]
If $t\ge1$, we apply assertion~(c) from the list above to conclude that
\[
(M_Q\varphi_{p(\cdot)})(t)\le(M_{s,Q}\varphi_{p(\cdot)})(t)\le 
K_0\,(M_{s,Q}\varphi^*_{p(\cdot)})^*(t)\le K_1\,(M_Q\varphi^*_{p(\cdot)})^*(t).
\]
If $0<t<1$ and $\alpha_Q(t)\le K_1$, then simply by definition of $\alpha_Q$, we have
\[
(M_Q\varphi_{p(\cdot)})(t)=\alpha_Q(t)\,(M_Q\varphi^*_{p(\cdot)})^*(t)
\le K_1\,(M_Q\varphi^*_{p(\cdot)})^*(t),
\]
so the claim holds in this case.

We are left with the case when $0<t<1$ and $\alpha_Q(t)>K_1$. Note that 
$0<t^{1/s}<1$, and it follows by assertion~(d) that 
\[
\alpha_{s,Q}(t^{1/s})\ge\left(\frac{\alpha_Q(t)}{C_2}\right)^{1/s}>
\left(\frac{K_1}{C_2}\right)^{1/s}=\max\{2K_0,C_1\}\ge C_1.
\]
Then assertion~(b) implies that 
\[
0<t^{1/s}<\lambda\,\min\left\{1,\frac1{\|\chi_Q\|_{p(\cdot)}}\right\}
\le\frac{\lambda}{\|\chi_Q\|_{p(\cdot)}},
\]
therefore $\mu(Q)\,(M_{s,Q}\varphi^*_{p(\cdot)})^*(t^{1/s})\le1$ by assertion~(a).
Using~(c) once again, and noting that 
$\alpha_{s,Q}(t^{1/s})\ge\max\{2K_0,C_1\}\ge2K_0$, we obtain 
\begin{align*}
2\,(M_{s,Q}\varphi_{p(\cdot)})(t^{1/s})&\le 2K_0\,(M_{s,Q}\varphi^*_{p(\cdot)})^*(t^{1/s})
+2b(Q)
\\
&\le \alpha_{s,Q}(t^{1/s})\,(M_{s,Q}\varphi^*_{p(\cdot)})^*(t^{1/s})+2b(Q)
\\
&=(M_{s,Q}\varphi_{p(\cdot)})(t^{1/s})+2b(Q),
\end{align*}
whence $(M_{s,Q}\varphi_{p(\cdot)})(t^{1/s})\le2b(Q)$. Then this, $0<t^{1-1/s}<1$, and the
convexity of $M_{s,Q}\varphi$ in the second variable imply
\begin{align*}
(M_Q\varphi_{p(\cdot)})(t)&\le(M_{s,Q}\varphi_{p(\cdot)})(t^{1/s}t^{1-1/s})
\\
&\le(M_{s,Q}\varphi_{p(\cdot)})(t^{1/s})\,t^{1-1/s}
\\
&\le2b(Q)\,t^{1-1/s},
\end{align*}
thus the claim holds in the remaining case. The auxiliary claim is proved, and we continue 
the main proof now.$\quad\square$

Take a sequence of families $\mathcal{Q}_k\in\mathfrak{S}$, $k\in\mathbb{Z}$, such that
\[
\sum_{k=-\infty}^\infty\sum_{Q\in\mathcal{Q}_k}
\mu(Q)\,(M_Q\varphi_{p(\cdot)}^*)^*(2^k)\le1.
\]
Then by the auxiliary claim,
\begin{align*}
\sum_{k=-\infty}^\infty&\sum_{Q\in\mathcal{Q}_k}
\mu(Q)\,(M_Q\varphi_{p(\cdot)})(2^k)
\\
&\le\sum_{k=-\infty}^{-1}\sum_{Q\in\mathcal{Q}_k}\mu(Q)\,
\big(K_1\,(M_Q\varphi^*_{p(\cdot)})^*(2^k)+2b(Q)\,2^{k(1-1/s)} \big)
\\
&\quad+\sum_{k=0}^{\infty}\sum_{Q\in\mathcal{Q}_k}\mu(Q)\,
K_1\,(M_Q\varphi^*_{p(\cdot)})^*(2^k)
\\
&\le K_1\,\sum_{k=-\infty}^\infty\sum_{Q\in\mathcal{Q}_k}
\mu(Q)\,(M_Q\varphi_{p(\cdot)}^*)^*(2^k)
+2\|b\|_{\mathfrak{S},1}\sum_{k=-\infty}^{-1}2^{k(1-1/s)}
\\
&\le K_1+\frac{2\|b\|_{\mathfrak{S},1}}{2^{1-1/s}-1}<\infty.
\end{align*}
Since $M_Q\varphi_{p(\cdot)}$ satisfies the $\Delta_2$-condition together with 
$\varphi_{p(\cdot)}$, it follows by Lemma~\ref{le:strongwith1} that 
$M_Q\varphi_{p(\cdot)}\ll(M_Q\varphi^*_{p(\cdot)})^*$. Thus, 
$p(\cdot)\in\mathscr{A}^{\mathscr{D}}_\textnormal{strong}$.

Conversely, if $p(\cdot)\in\mathscr{A}^{\mathscr{D}}_\textnormal{strong}$, then by
Remark~\ref{re:strong-usual} we have $p(\cdot)\in\mathscr{A^D}$, which follows from 
the $\Delta_2$-condition for $\varphi_{p(\cdot)}$.
\end{proof}

\section{The Main Result}\label{sec:main}

After a careful reading of Sections~5--7, the reader might have noticed that 
the central theorems of these sections---Theorem~\ref{th:self-improvment} 
(the self-improvement property of the class $\mathscr{A^D}$), Theorem~\ref{th:coincidence}
(equivalence of the classes $\mathscr{A^D}$ and 
$\mathscr{A}^{\mathscr{D}}_\textnormal{strong}$ for nice ``variable-Lebesgue'' functions),
and Theorem~\ref{th:sufficient} (a sufficient condition for the boundedness of the 
maximal operator)---assemble into the following elegant line of reasoning,
modeled on the discussion preceding~\cite[Theorem~5.7.1]{DHHR11}. We
take a ``variable-Lebesgue'' function $\varphi_{p(\cdot)}$, with $1<p_-\le p_+<\infty$, 
which belongs the class $\mathscr{A^D}$. The class is self-improving, hence a 
``slowed-down'' version of $\varphi_{p(\cdot)}$---a function $\varphi_{sp(\cdot)}$ with 
some $s$ slightly below one (so that $sp_->1$)---is also in the class $\mathscr{A^D}$, or
equivalently, in the class $\mathscr{A}^{\mathscr{D}}_\textnormal{strong}$. Finally, the
presence of the ``slowed-down'' function $\varphi_{sp(\cdot)}$ in  
$\mathscr{A}^{\mathscr{D}}_\textnormal{strong}$ implies that the maximal operator 
$M$ is bounded on $L^{p(\cdot)}(X)$.

These three logical steps yield the main result of our work: it suffices for the 
boundedness of $M$ on a space $L^{p(\cdot)}(X)$ with $1<p_-\le p_+<\infty$ that 
$p(\cdot)\in\mathscr{A^D}$---i.e., that there exists $C>0$ such that 
\[
\bigg\|\sum_{Q\in\mathcal{Q}}\bigg(\fint_Q|f(x)|\,d\mu(x)\bigg)\,\chi_Q\bigg\|_{p(\cdot)}
\le C\,\|f\|_{p(\cdot)}
\]
holds for all $f\in L^{p(\cdot)}(X)$ and all families $\mathcal{Q}$ of pairwise disjoint 
cubes from the dyadic system $\mathscr{D}$ on $X$. Almost immediately, this sufficient 
condition also turns out to be necessary. Thus, we obtain a full characterization of the
boundedness of the \emph{maximal operator $M$} on the variable Lebesgue spaces 
$L^{p(\cdot)}(X)$, $1<p_-\le p_+<\infty$, in terms of the uniform boundedness
on $L^{p(\cdot)}(X)$ of the \emph{averaging operators} over all families of pairwise
disjoint dyadic cubes in $X$. 
\begin{theorem}[Main Result]\label{th:main}
Let $X$ be an unbounded space of homogeneous type equipped with a Borel-semiregular measure,
and let $\mathscr{D}$ be a dyadic system associated with $X$ by Theorem~\ref{th:HK}. Suppose 
that a variable exponent $p(\cdot)\in\mathcal{P}(X)$ is such that $1<p_-\le p_+<\infty$.
Then the maximal operator $M$ is bounded on $L^{p(\cdot)}(X)$ if and only if 
$p(\cdot)\in\mathscr{A^D}$.
\end{theorem}
\begin{proof}
The necessity part immediately follows from inequality~\eqref{eq:aver-max} and the fact that 
the dyadic maximal operator $M^\mathscr{D}$ is also bounded on $L^{p(\cdot)}(X)$ 
in view of Lemma~\ref{le:M-and-dyadicM}.

Conversely, let $p(\cdot)\in\mathscr{A^D}$. By Lemma~\ref{le:phi-bar-tilde}, the condition
$1<p_-\le p_+<\infty$ guarantees that both $\varphi_{p(\cdot)}$ and $\varphi^*_{p(\cdot)}$ 
satisfy the $\Delta_2$-condition. Then it follows from the self-improvement property of the 
class $\mathscr{A^D}$ in Theorem~\ref{th:self-improvment}, applied to the function 
$\varphi_{p(\cdot)}$, that there exists $s\in(1/p_-,1)$ such
that $q(\cdot)\in\mathscr{A^D}$ with $q(x):=sp(x)$ for all $x\in X$; for this, notice that 
\[
\varphi_{p(\cdot)}(x,t^s)\approx\bar{\varphi}_{p(\cdot)}(x,t^s)
=\bar{\varphi}_{sp(\cdot)}(x,t)\approx\varphi_{sp(\cdot)}(x,t),
\]
uniformly in $x\in X$ and $t\ge0$, by Lemma~\ref{le:phi-bar-tilde}(a).
Since $1<q_-\le q_+<\infty$, Theorem~\ref{th:coincidence} implies that 
$q(\cdot)\in\mathscr{A}^{\mathscr{D}}_\textnormal{strong}$. 
Finally, observe that for any $F\subseteq X$ with $\mu(F)=\infty$ and $\lambda>0$,
there holds 
\[
\rho_{\bar{\varphi}_{p(\cdot)}}(\lambda\chi_F)
=\int_F \lambda^{p(x)}\,d\mu(x)
\ge\min\{\lambda^{p_-},\lambda^{p_+}\}\,\mu(F)=\infty,
\]
so that $\chi_F\notin L^{p(\cdot)}(X)$. Applying Theorem~\ref{th:sufficient}, 
we derive from $q(\cdot)\in\mathscr{A}^{\mathscr{D}}_\textnormal{strong}$ that 
$M$ is bounded on $L^{p(\cdot)}(X)$.
\end{proof}

It is important to note that our Theorem~\ref{th:main} is a generalization of
the classical Diening's result in~\cite[Theorem~8.1, (i)$\Leftrightarrow$(ii)]{LD05}
from the Euclidean setting to spaces of homogeneous type. Since the Euclidean space 
$\mathbb{R}^n$, equipped with the Euclidean distance and the $n$-dimensional Lebesgue 
measure, is an example of an unbounded space of homogeneous type in the sense of 
Definition~\ref{def:SHT}, and the Lebesgue measure is Borel-semiregular (even more so, 
Borel-regular), the following corollary of Theorem~\ref{th:main} holds in view of 
Remark~\ref{re:Euclidean}: 
\begin{corollary}\label{cor:main}
Let $p(\cdot)\in\mathcal{P}(\mathbb{R}^n)$ satisfy $1<p_-\le p_+<\infty$. The maximal 
operator $M$ is bounded on $L^{p(\cdot)}(\mathbb{R}^n)$ if and only if there exists
$C>0$ such that 
\[
\bigg\|\sum_{Q\in\mathcal{Q}}\bigg(\fint_Q|f(x)|\,dx\bigg)
\,\chi_Q\bigg\|_{p(\cdot)}\le C\,\|f\|_{p(\cdot)}
\]
holds for all $f\in L^{p(\cdot)}(\mathbb{R}^n)$ and all families $\mathcal{Q}$ of 
pairwise disjoint cubes in $\mathbb{R}^n$ selected from $3^n$ adjacent dyadic grids 
\[
\mathscr{D}^{\mathbf{t}}:=\left\{
2^{-k} \big([0,1)^n+\mathbf{m}+(-1)^k\mathbf{t}\big)
\ :\ 
k\in\mathbb{Z},\,\mathbf{m}\in\mathbb{Z}^n\right\},
\quad\mathbf{t}\in\left\{0,\frac13,\frac23\right\}^n.
\] 
\end{corollary}
This provides a sharper formulation of~\cite[Theorem~8.1, (i)$\Leftrightarrow$(ii)]{LD05}
and~\cite[Theorem~5.7.2, (a)$\Leftrightarrow$(c)]{DHHR11}, allowing one to characterize
the boundedness of the maximal operator on the ``nice'' spaces $L^{p(\cdot)}(\mathbb{R}^n)$
via the averaging operators over families of pairwise disjoint \emph{dyadic} cubes 
rather than all cubes in $\mathbb{R}^n$.

With the help of Theorem~\ref{th:main}, we can also rediscover two known facts concerning
the boundedness of $M$ on $L^{p(\cdot)}(X)$ \emph{in the special case when
$1<p_-\le p_+<\infty$ and $X$ is an unbounded space of homogeneous type with a 
Borel-semiregular measure}. 
The first of these facts states that $M$ is bounded on a variable Lebesgue space 
$L^{p(\cdot)}(X)$---with $1<p_-\le p_+<\infty$, over an arbitrary space
of homogeneous type $X$---if and only if $M$ is bounded on $L^{p'(\cdot)}(X)$; this was
proved by Karlovich in~\cite[Theorem~1.1]{AK20}. The second is part of the self-improving
boundedness property of $M$, recently obtained by the author in~\cite[Corollary~4.5]{AS25}
(see also~\cite[Theorem~1.1]{AS25a}), which asserts that $M$ is bounded on a variable 
Lebesgue space $L^{p(\cdot)}(X)$---with $p_->1$, over an arbitrary space of homogeneous 
type $X$---if and only if $M$ is bounded on $L^{sp(\cdot)}(X)$ for some $s\in(1/p_-,1)$.
The following corollary of Theorem~\ref{th:main}, giving the special-case versions of
these facts, extends~\cite[Theorem~8.1, 
(ii)$\Leftrightarrow$(iv)$\Leftrightarrow$(vi)]{LD05}. 
\begin{corollary}
Under the hypotheses of Theorem~\ref{th:main}, the following are equivalent:
\begin{enumerate} 
\item[(a)] $M$ is bounded on $L^{p(\cdot)}(X)$;
\item[(b)] $M$ is bounded on $L^{p'(\cdot)}(X)$;
\item[(c)] $M$ is bounded on $L^{sp(\cdot)}(X)$ for some $s\in(1/p_-,1)$. 
\end{enumerate}
\end{corollary}
\begin{proof}
Let $M$ be bounded on $L^{p(\cdot)}(X)$. By Theorem~\ref{th:main}, this is equivalent to
$p(\cdot)\in\mathscr{A^D}$. In view of the duality property of the class
$\mathscr{A^D}$ in Lemma~\ref{le:A-duality}, $p(\cdot)\in\mathscr{A^D}$ if and only if
$\varphi_{p(\cdot)}^*\in\mathscr{A^D}$. Note that
$\varphi_{p(\cdot)}^*(x,t)\approx\varphi_{p'(\cdot)}(x,t)$, uniformly in $x\in X$ and 
$t\ge0$, by Lemma~\ref{le:phi-bar-tilde}(b). Then we apply Lemma~\ref{le:comparison}(b) 
with $s=1$ to conclude that $p'(\cdot)\in\mathscr{A^D}$, hence $M$ is bounded 
on $L^{p'(\cdot)}(X)$ by Theorem~\ref{th:main}. This proves implication
(a)$\Rightarrow$(b). The converse implication (b)$\Rightarrow$(a) follows by a similar
argument, since $\varphi_{p'(\cdot)}^*(x,t)\approx\varphi_{p(\cdot)}(x,t)$ uniformly in
$x\in X$ and $t\ge0$.

Suppose again that $M$ is bounded on $L^{p(\cdot)}(X)$, or equivalently, that 
$p(\cdot)\in\mathscr{A^D}$. Since $\varphi^*_{p(\cdot)}$ satisfies the 
$\Delta_2$-condition by Lemma~\ref{le:phi-bar-tilde}, it follows from the self-improvement
property of the class $\mathscr{A^D}$ in Lemma~\ref{th:self-improvment} that 
there exists $s\in(1/p_-,1)$ such that $sp(\cdot)\in\mathscr{A^D}$, owing to the equivalence 
$\varphi_{sp(\cdot)}(x,t)\approx\varphi_{p(\cdot)}(x,t^s)$ uniformly in $x\in X$ and 
$t\ge0$. Then $M$ is bounded on $L^{sp(\cdot)}(X)$ due to Theorem~\ref{th:main}, so 
implication (a)$\Rightarrow$(c) is proved. Conversely, assumption~(c) implies 
$sp(\cdot)\in\mathscr{A^D}$, and hence 
$sp(\cdot)\in\mathscr{A}^{\mathscr{D}}_\textnormal{strong}$ by Theorem~\ref{th:coincidence}.
Noting that $\varphi_{sp(\cdot)}$ satisfies the $\Delta_2$-condition
and $\varphi_{p(\cdot)}(x,t)\approx\varphi_{sp(\cdot)}(x,t^{1/s})$ 
holds with $1/s>1$ uniformly in $x\in X$ and
$t\ge0$, we derive from Lemma~\ref{le:strong-strong} that 
$p(\cdot)\in\mathscr{A}^{\mathscr{D}}_\textnormal{strong}$. Hence 
$p(\cdot)\in\mathscr{A^D}$, which implies by Theorem~\ref{th:coincidence} that $M$ is 
bounded on $L^{p(\cdot)}(X)$. We thus arrived at (a).
\end{proof}

\medskip\noindent
{\bf Acknowledgments} 

\medskip
This work is funded by national funds through the FCT -- 
Funda\c{c}\~{a}o para a Ci\^{e}ncia~e a Tecnologia, I.P., under the scope 
of the projects
\begin{itemize}
    \item \texttt{UIDB/00297/2020} (https://doi.org/10.54499/UIDB/00297/2020),
    \item \texttt{UIDP/00297/2020} (https://doi.org/10.54499/UIDP/00297/2020)
\end{itemize}
(Center for Mathematics and Applications) and under the scope of the PhD scholarship
\texttt{UI/BD/154284/2022}. I am grateful to Dr. Ryan Alvarado for the idea behind the proof
of Lemma~\ref{le:cones}(a), and to Dr. David Cruz-Uribe for communicating this idea to me. I 
would also like to thank my scientific advisor, Dr. Oleksiy Karlovych, for his 
constructive remarks during the preparation of this manuscript.

\bibliographystyle{plain}
\bibliography{arcubes-v2}

\end{document}